\documentclass[11pt]{article}
\usepackage[a4paper,margin=28mm]{geometry}
\usepackage{amsmath,amssymb,amsthm,mathtools}
\usepackage{newtxtext,newtxmath}
\usepackage{enumitem}
\usepackage{booktabs}
\usepackage{array}
\usepackage{bm}
\usepackage{xurl}
\usepackage{hyperref}
\usepackage[nameinlink,capitalise,noabbrev]{cleveref}
\usepackage{aliascnt}
\usepackage{needspace}
\usepackage{fancyhdr}
\usepackage{setspace}
\hypersetup{
  hidelinks,
  pdftitle={Time-Periodic Chemotaxis in Shifting Environments: Propagation and Planar Forced Waves in R N},
  pdfauthor={Hoang-Hung Vo}
}
\numberwithin{equation}{section}

\newtheorem{theorem}{Theorem}[section]
\newaliascnt{proposition}{theorem}
\newtheorem{proposition}[proposition]{Proposition}
\aliascntresetthe{proposition}
\newaliascnt{lemma}{theorem}
\newtheorem{lemma}[lemma]{Lemma}
\aliascntresetthe{lemma}
\newaliascnt{corollary}{theorem}
\newtheorem{corollary}[corollary]{Corollary}
\aliascntresetthe{corollary}
\newaliascnt{claim}{theorem}

\aliascntresetthe{claim}
\theoremstyle{definition}
\newaliascnt{definition}{theorem}
\newtheorem{definition}[definition]{Definition}
\aliascntresetthe{definition}
\newaliascnt{assumption}{theorem}
\newtheorem{assumption}[assumption]{Assumption}
\aliascntresetthe{assumption}
\theoremstyle{remark}
\newaliascnt{remark}{theorem}
\newtheorem{remark}[remark]{Remark}
\aliascntresetthe{remark}
\crefname{theorem}{Theorem}{Theorems}
\crefname{proposition}{Proposition}{Propositions}
\crefname{lemma}{Lemma}{Lemmas}
\crefname{corollary}{Corollary}{Corollaries}
\crefname{claim}{Claim}{Claims}
\crefname{definition}{Definition}{Definitions}
\crefname{assumption}{Assumption}{Assumptions}
\crefname{remark}{Remark}{Remarks}

\newcommand{\R}{\mathbb R}
\newcommand{\eps}{\varepsilon}
\newcommand{\p}{\partial}

\newcommand{\per}{\mathrm{per}}
\newcommand{\CalA}{\mathcal A}
\newcommand{\CalB}{\mathcal B}

\newcommand{\CalS}{\mathcal S}
\newcommand{\CalU}{\mathcal U}

\title{\bfseries Time-Periodic Chemotaxis in Shifting Environments: Propagation and Planar Forced Waves in $\R^N$}
\author{Hoang-Hung Vo$^{1,2}$\thanks{Corresponding author. E-mail: \texttt{vhhung@hcmiu.edu.vn}.}\\[1mm]
\small $^{1}$Department of Mathematics, International University, Ho Chi Minh City, Vietnam\\
\small $^{2}$Vietnam National University Ho Chi Minh City, Ho Chi Minh City, Vietnam}
\date{}

\begin{document}
\maketitle

\begin{abstract}
We study a Keller--Segel system in a moving habitat with time-periodic coefficients. When the signal obeys an elliptic equation ($\tau=0$), we establish global boundedness, sharp spreading bounds with intrinsic speed $c_*=2\sqrt{\bar a\,\bar r}$, and persistence inside the intersection of the KPP spreading cone and the moving favorable half-space. A localized favorable habitat is governed by a generalized periodic Floquet exponent; the extinction criterion for a negative exponent requires an additional weighted-resolvent hypothesis. For each fixed chemical relaxation parameter $\tau\ge0$, constant population diffusivity, and sufficiently weak taxis, we also obtain a locally unique branch of planar forced time-periodic waves near a pinned scalar profile, under the nondegeneracy and regularity assumptions specified below. We also determine the dimension-dependent small-slope behavior of the Yukawa first-contact coefficient, the exact quadratic dependence of the slab Floquet exponent on habitat speed for constant population diffusivity, and the uniform quadratic remainder in the weak-taxis expansion of the pinned waves. The proofs use a Yukawa first-contact estimate, a localized power-Harnack inequality, and a Fredholm argument in two-sided weighted spaces. The results are compared with the one-dimensional shifting-habitat theory of Shen--Xue \cite{SX2020JDE} and the homogeneous transition fronts of Salako--Shen \cite{SalakoShen2018M3AS}.
\end{abstract}



\Needspace{5\baselineskip}
\section{Introduction and Main Results}\label{sec:intro-main}
A population subject to a shifting environment must combine its intrinsic spatial propagation with its ability to remain inside a favorable region moving at an imposed speed.  Potapov--Lewis \cite{PotapovLewis2004} and Berestycki--Diekmann--Nagelkerke--Zegeling \cite{BDNZ2009} identified this competition in moving-habitat models.  Berestycki--Rossi \cite{BR2008,BR2009} developed persistence and generalized-eigenvalue criteria in whole-space and cylindrical geometries, while Berestycki--Fang \cite{BF2018} and Bouhours--Giletti \cite{BouhoursGiletti2019} treated forced waves and spreading/vanishing for scalar monostable equations.  Du--Wei--Zhou \cite{DuWeiZhou2018} considered a free-boundary shifting environment, De Leenheer--Shen--Zhang \cite{DeLeenheerShenZhang2020} treated nonlocal dispersal in moving habitats, and Hou--Wang--Zhao \cite{HouWangZhao2024} studied cooperative reaction--diffusion systems in periodic shifting environments.  The survey \cite{WangLiDongQiao2022} provides a broader account of scalar and system propagation in shifting media.

Time dependence changes the stationary spectral problem in the moving frame into a Floquet problem.  Fang--Peng--Zhao \cite{FPZ2021} established forced time-periodic waves together with propagation dynamics for the scalar periodic shifting equation
\begin{equation}\label{eq:FPZintro}
 u_t=u_{xx}+u\,g(t,x-ct,u),
 \; t>0,\; x\in\mathbb R.
\end{equation}
We use their scalar result when taxis vanishes.  Zhang--Zhao \cite{ZhangZhao2020} obtained periodic forced waves and spreading dynamics for nonlocal Fisher--KPP dispersal, Shen--Shen--Xue--Zhou \cite{SSXZ2022} derived persistence criteria and quantified the effect of temporal fluctuations, and Zhang--Zhao \cite{ZhangZhao2025} treated spreading and forced traveling waves under general time heterogeneity.

The chemotactic drift prevents direct application of scalar comparison arguments.  Shen--Xue \cite{SX2020JDE} considered the one-dimensional parabolic--elliptic model
\begin{equation}\label{eq:SXintro}
\begin{cases}
 u_t=u_{xx}-\chi(uv_x)_x+u\big(r(x-ct)-bu\big),\\
 0=v_{xx}-\nu v+\mu u,
\end{cases}
\; x\in\mathbb R,
\end{equation}
and established propagation and persistence results for shifting habitats; forced waves were later constructed in \cite{SX2022FW}.  That forced-wave theory is one-dimensional and parabolic--elliptic.  Fully parabolic traveling waves for spatially homogeneous logistic Keller--Segel systems were constructed in \cite{SalakoShen2020}, whereas persistence and the associated homogeneous spreading speed were studied in \cite{ShenXue2022Persistence,SX2022PP}.  These results do not concern forced time-periodic waves pinned to a shifting heterogeneous habitat. In particular, they do not treat wave continuation with temporal signal memory in that setting.  The Keller--Segel flow is not order preserving: ordered population densities need not produce ordered chemical drifts.  Scalar comparison results therefore do not extend directly to the coupled system.

The distinctions concern the geometry of the favorable region, the evolution of the signal, and the time dependence of the coefficients. Salako--Shen \cite{SalakoShen2017JDE,SalakoShen2017DCDS} established well-posedness and spreading or traveling-wave results for logistic Keller--Segel systems in spatially homogeneous environments. Their three-part study \cite{SalakoShen2018M3AS,SalakoShen2018JMAA,SalakoShen2022Fronts} treats persistence and positive entire solutions for space--time-dependent logistic sources, and constructs transition fronts in the spatially homogeneous logistic specialization on $\mathbb R^N$. The associated media are not translating favorable--unfavorable habitats with time-periodic forcing, and their transition fronts should not be confused with forced waves pinned to a moving heterogeneous resource profile. The sharp-speed criterion in one dimension of Salako--Shen--Xue \cite{SalakoShenXue2019} is a particularly important predecessor of the damping inequality \eqref{eq:Dass}; our explicit identity \eqref{eq:Theta1} recovers its kernel factor exactly in the homogeneous one-dimensional specialization.

Shen--Xue \cite{SX2020JDE,SX2022FW} treated the shifting one-dimensional parabolic--elliptic Keller--Segel problem and are the closest precedents for the propagation and forced-wave results considered here. The parabolic--parabolic results \cite{ShenXue2022Persistence,SX2022PP} address global persistence and sharp spreading in a \emph{spatially homogeneous} logistic medium; they do not establish periodic forced waves in a moving heterogeneous habitat. The multidimensional periodic monotone-semiflow results of Du--Li--Shen \cite{DuLiShen2022JFA} concern order-preserving cooperative dynamics, while the coupled Keller--Segel evolution has no such general order structure. Shen's theory for time-dependent scalar fronts \cite{Shen2011Monostable}, the space--time periodic KPP variational principle \cite{Shen2010Variational}, and the fluctuation-sensitive shifting-population criteria of Shen--Shen--Xue--Zhou \cite{SSXZ2022} supply further temporal and spectral context, but not the chemotactic first-contact estimates required below. Shen's recent preprint \cite{Shen2026Waves} establishes existence, uniqueness,
stability and monotonicity results for one-dimensional autonomous chemotaxis
traveling waves with generalized nonlinear taxis and logistic sources; neither
an externally translating resource profile nor time-periodic forcing appears in
that model. Finally, Hassan--Shen--Zhang \cite{HassanShenZhang2025} studied spreading in a different chemotaxis--consumption model with logistic growth; its chemical equation is not the signal-production equation in \eqref{eq:fullmodel}. Recent work of Chen--Ruau--Shen \cite{ChenRuauShen2026I,ChenRuauShen2026II} concerns boundedness and persistence for signal-dependent sensitivity with logistic sources on bounded domains. These results provide complementary chemotaxis estimates, but do not treat the propagation geometry or time-periodic forced waves in translating habitats considered here. The equations and hypotheses in these papers differ from those considered below.

The purpose of this paper is to study multidimensional population dynamics in a planar shifting habitat with seasonal coefficients.  We consider
\begin{equation}\label{eq:fullmodel}
\begin{cases}
 u_t=a(t)\Delta u-\eta\chi(t)\nabla\!\cdot\!\big(u\nabla v\big)
      +u f\big(t,x\cdot e-ct,u\big),\\
 \tau v_t=d(t)\Delta v-\nu(t)v+\mu(t)u,
\end{cases}
\; t>0,\; x\in\mathbb R^N.
\end{equation}
Here $\eta\ge0$ is the taxis-strength parameter and $\tau\ge0$ is the signal-relaxation parameter; $\tau=0$ corresponds to the parabolic--elliptic regime.  We assume $N\ge1$ and $e\in\mathbb S^{N-1}$.  The transport and signal coefficients are spatially homogeneous and periodic in time, whereas the demographic heterogeneity is carried by the single moving coordinate $x\cdot e-ct$.  The assumptions are stated below.  They include a favorable KPP limit behind the climate front and an unfavorable limit ahead of it.

Our first group of results concerns the parabolic--elliptic case $\tau=0$.  Under an explicit logistic damping condition, solutions issued from bounded nonnegative data are global and uniformly bounded.  For compactly supported nontrivial data, the intrinsic KPP speed is
\(c_*=2\sqrt{\bar a\,\bar r}.\)
The upper propagation bound is the Euclidean cone with speed $c_*$.  Persistence in a moving vector frame with velocity $s\in\mathbb R^N$ occurs whenever
\(|s|<c_*, \; s\cdot e<c.\)
Thus the persistent region is determined simultaneously by the intrinsic KPP cone and by the translating favorable half-space.  Under a stronger damping gap, the population and signal converge in the interior of this region to the unique positive periodic favorable state.  For one-sided occupied data we obtain the corresponding longitudinal threshold, and for a moving favorable slab we obtain persistence and extinction criteria in terms of a generalized principal Floquet exponent.  The critical zero-exponent case is not included in the statement.

Two specializations recover earlier models.  In dimension one, with constant logistic coefficients and $\tau=0$, the propagation statements reduce to those of \cite{SX2020JDE}.  When $\eta=0$, the wave problem reduces, after the stated spatial rescaling, to \eqref{eq:FPZintro}.  The remaining assumptions separate those needed for propagation from the additional regularity and pinning hypotheses used only for the continuation of forced waves. The quantitative refinements of the existence and threshold results appear in \cref{prop:YukawaSmallSlope,prop:FloquetSpeed,thm:linearresponse}: the dimensional small-slope law of the Yukawa correction, the exact climate-speed shift of the localized Floquet exponent, and the first-order weak-taxis profile response with a uniform quadratic remainder.

The hypotheses and main statements follow.
For a $T$-periodic scalar function $h$, we write
$$\langle h\rangle=\frac1T\int_0^T h(s)\,ds$$
and set
\[
 a_*:=\min a,\; a^*:=\max a,\;
 d_*:=\min d,\; d^*:=\max d,\;
 \nu_*:=\min\nu,
\]
and use analogous notation for the remaining periodic coefficients.  We always assume $e\in\mathbb S^{N-1}$.

\begin{assumption}[Periodic transport and signal coefficients]\label{A1}
For some $\alpha\in(0,1)$,
\(a,d,\chi,\nu,\mu\in C^\alpha_{\per}(\mathbb R),\)
$a,d,\nu,\mu$ are strictly positive and $\chi\ge0$.
\end{assumption}

\begin{assumption}[KPP reaction and time-dependent damping]\label{A2}
The map $f=f(t,\xi,u)$ is $T$-periodic in $t$ and belongs locally to $C^{\alpha/2,\alpha,1}$ in $(t,\xi,u)$.  In addition, for every $M>0$ its $C^{\alpha/2,\alpha}$ norm in $(t,\xi)$, its $u$-Lipschitz norm, and the corresponding local mixed H\"older seminorms are bounded uniformly on
\(\mathbb R\times\mathbb R\times[0,M].\)
These uniform bounds are used in the global parabolic estimates and in compactness arguments under spatial translations.  There is a strictly positive function
\(\beta\in C^\alpha_{\per}(\mathbb R), \; \beta_*:=\min_{t\in[0,T]}\beta(t)>0,\)
such that
\begin{equation}\label{eq:fuass}
 f_u(t,\xi,u)\le-\beta(t)
 \;(t,\xi,u)\in\mathbb R\times\mathbb R\times[0,\infty).
\end{equation}
For each $M>0$ there is $b_1(M)<\infty$ such that $-b_1(M)\le f_u(t,\xi,u)\le-\beta(t)$ for $0\le u\le M$, uniformly in $(t,\xi)$.  The linear growth rate
\(r(t,\xi):=f(t,\xi,0)\)
is bounded.  Consequently, for $0\le u\le M$,
\begin{equation}\label{eq:kppbounds}
 r(t,\xi)-b_1(M)u\le f(t,\xi,u)\le r(t,\xi)-\beta(t)u.
\end{equation}
The pointwise damping profile $\beta(t)$ enters the taxis inequalities below.
\end{assumption}

The quantitative slope hypothesis in \cref{A2} is required for the nonperturbative chemotaxis estimates.  It is not imposed on the purely scalar specialization $\eta=0$ when we state the Fang--Peng--Zhao scalar theorem in \cref{sec:recover}; that result is stated there under the original weaker KPP assumptions.

\begin{assumption}[Separated favorable and unfavorable far fields]\label{F}
There exist $T$-periodic functions $f_\pm(t,u)$ such that
\(f(t,\xi,u)\to f_\pm(t,u)\;\text{as }\xi\to\pm\infty\)
locally uniformly in $u\ge0$ and uniformly in $t$.  We set $r_\pm(t)=f_\pm(t,0)$ and assume
\(\langle r_-\rangle>0>\langle r_+\rangle\)
and
\(f_+(t,u)\le f(t,\xi,u)\le f_-(t,u)\)
for all $t,\xi,u\ge0$.
\end{assumption}

Thus the favorable region lies behind the moving hyperplane $x\cdot e=ct$.  The opposite orientation follows by replacing $e$ by $-e$ and $c$ by $-c$.

For $\kappa>0$ let $G_{N,\kappa}$ be the positive fundamental solution of $-\Delta+\kappa^2$ on $\mathbb R^N$, normalized by
\begin{equation}\label{eq:Yukawaheat}
 G_{N,\kappa}(z)=\int_0^\infty
 (4\pi s)^{-N/2}\exp\!\left(-\frac{|z|^2}{4s}-\kappa^2s\right)\,ds.
\end{equation}
We set $G_N:=G_{N,1}$.  By rotational invariance the following quantity does not depend on the choice of unit vector $e$:
\begin{equation}\label{eq:ThetaN}
 \Theta_N(\rho)
 :=\int_{\mathbb R^N}
 \big[\rho\,\p_eG_N(z)-G_N(z)\big]_+e^{\rho z\cdot e}\,dz,
 \; \rho\ge0,
\end{equation}
where $\p_e=e\cdot\nabla$ and $[q]_+=\max\{q,0\}$.  The integral is finite for every $\rho\ge0$; this is proved in \cref{lem:YukawaFacts}.  For $N=1$,
\begin{equation}\label{eq:Theta1}
 \Theta_1(\rho)=\frac12\frac{(\rho-1)_+}{\rho+1}.
\end{equation}

For $\tau=0$ set $\kappa(t)=\sqrt{\nu(t)/d(t)}$ and $\gamma_0(t)=\eta\chi(t)\mu(t)/d(t)$; for $\lambda\ge0$ set $\gamma_{N,\lambda}(t)=\gamma_0(t)[1+\Theta_N(\lambda/\kappa(t))]$.
Also put
\begin{equation}\label{eq:lambda0cstar}
 \bar a:=\langle a\rangle,\; \bar r:=\langle r_-\rangle,\;
 \lambda_*:=\sqrt{\frac{\bar r}{\bar a}},\; c_*:=2\sqrt{\bar a\bar r}.
\end{equation}

For the parabolic--elliptic problem $\tau=0$, define the basic damping gap
\(\delta_0:=\min_{t\in[0,T]}\big(\beta(t)-\gamma_0(t)\big).\)
This quantity is used already in the global boundedness and localized-slab arguments, which do not require the optimized leading-edge condition below.

\begin{assumption}[Optimized damping condition for the front problem]\label{D}
For the parabolic--elliptic problem $\tau=0$,
\begin{equation}\label{eq:Dass}
 \delta_0>0,
 \;
 \beta(t)\ge\gamma_{N,\lambda_*}(t)\;\text{for every }t.
\end{equation}
\end{assumption}

\begin{remark}[Optimization and one-dimensional exactness]\label{rem:SXcondition}
The second inequality in \eqref{eq:Dass} is imposed only at the minimizing KPP exponent $\lambda_*$.  No supremum over $0<\lambda<\lambda_*$ enters the proof: the upper cone with the exact linear threshold $c_*$ is obtained from the single exponential rate $\lambda_*$ and an auxiliary comparison speed strictly larger than $c_*$.  Hence the argument requires only the condition actually used at first contact and does not assume any monotonicity of $\Theta_N$.  We do not claim that \eqref{eq:Dass} is a necessary nonlinear damping condition in dimensions $N\ge2$; it is the optimized sufficient condition obtained by the exact first-contact estimate at the minimizing KPP exponent.  In dimension one it reduces exactly to the Shen--Xue coefficient condition, as shown below.

When $N=1$, \eqref{eq:Theta1} implies
\[
 \gamma_{1,\lambda_*}(t)=
 \frac{\eta\chi(t)\mu(t)}{d(t)}
 \left(1+\frac12
 \frac{(\lambda_*-\kappa(t))_+}{\lambda_*+\kappa(t)}\right).
\]
In the one-dimensional constant-coefficient separated-habitat case, after the reflection used in Section~\ref{sec:recover}, the favorable limit satisfies $r_-=r^*$ and hence $\lambda_*=\sqrt{r^*}$.  With $a=d=1$ and constant $\chi,\mu,\nu,\beta$, we have
\[
 \gamma_{1,\lambda_*}
 =\eta\chi\mu\left(1+\frac12
 \frac{(\sqrt{r^*}-\sqrt\nu)_+}{\sqrt{r^*}+\sqrt\nu}\right).
\]
For $\eta=1$ this coincides with the Shen--Xue coefficient condition in \eqref{eq:SXdampingexact}.
\end{remark}

The favorable scalar ordinary differential equation is
\begin{equation}\label{eq:pminus-early}
 p'=p f_-(t,p).
\end{equation}
Under \cref{A2,F}, this equation has a unique positive $T$-periodic solution, denoted by $p_-(t)$; its existence, uniqueness and uniform positive bounds are proved in \cref{lem:pminus} immediately before the favorable-medium rigidity argument where they are first needed in a proof.  In the elliptic signal regime we write
\begin{equation}\label{eq:qminus-ell}
 q_-(t)=\frac{\mu(t)}{\nu(t)}p_-(t).
\end{equation}

\begin{theorem}[Sharp geometric front propagation]\label{thm:front}
Suppose that \cref{A1,A2,F,D} hold and that $\tau=0$.  Let $u_0\ge0$ be bounded, uniformly continuous, nontrivial and compactly supported.  Then \eqref{eq:fullmodel} has a unique global bounded classical solution $(u,v)$ for $t>0$.  The following alternatives hold.
\begin{enumerate}[label=\textup{(\roman*)},leftmargin=8mm]
\item If $c<-c_*$, then
\(\lim_{t\to\infty}\lVert u(t,\cdot)\rVert_\infty=0.\)
\item If $-c_*<c<c_*$, then for every $\eps>0$ with
\(0<\eps<\frac12\min\{c+c_*,c_*-c\},\)
writing $\mathcal K_\eps(t)=\{x\in\mathbb R^N:|x|\le(c_*-\eps)t,\ x\cdot e\le(c-\eps)t\}$, we have
\begin{equation}\label{eq:frontmoderate}
 \lim_{t\to\infty}\sup_{|x|\ge(c_*+\eps)t}u(t,x)=0,\;
 \lim_{t\to\infty}\sup_{x\cdot e\ge(c+\eps)t}u(t,x)=0,\;
 \liminf_{t\to\infty}\inf_{x\in\mathcal K_\eps(t)}u(t,x)>0.
\end{equation}
If, in addition,
\begin{equation}\label{eq:strongdamp}
 \delta_2:=\min_{t\in[0,T]}\big(\beta(t)-2\gamma_0(t)\big)>0,
\end{equation}
then
\begin{equation}\label{eq:frontconvmoderate}
 \lim_{t\to\infty}
 \sup_{\substack{|x|\le(c_*-\eps)t\\x\cdot e\le(c-\eps)t}}
 \left(|u(t,x)-p_-(t)|+|v(t,x)-q_-(t)|\right)=0.
\end{equation}
\item If $c\ge c_*$, then for every $\eps\in(0,c_*)$,
\begin{equation}
 \lim_{t\to\infty}\sup_{|x|\ge(c_*+\eps)t}u(t,x)=0,\;
 \liminf_{t\to\infty}\inf_{|x|\le(c_*-\eps)t}u(t,x)>0.
 \label{eq:frontdecayfastgood}
\end{equation}
Under \eqref{eq:strongdamp},
\begin{equation}\label{eq:frontconvfastgood}
 \lim_{t\to\infty}\sup_{|x|\le(c_*-\eps)t}
 \left(|u(t,x)-p_-(t)|+|v(t,x)-q_-(t)|\right)=0.
\end{equation}
\end{enumerate}
Thus the sharp linear speed is independent of taxis within the damping regime, while the geometry of the persistent region is the intersection of the intrinsic Euclidean cone with the favorable moving half-space.  No conclusion is asserted here at the critical value $c=-c_*$ when $\eta>0$.
\end{theorem}

The critical value $c=-c_*$ for compactly supported data is therefore left open for nonzero taxis.  The one-dimensional scalar problem has a sharper critical description under additional tail hypotheses.

\begin{theorem}[One-sided occupied favorable half-space]\label{thm:frontfilled}
Suppose that \cref{A1,A2,F,D} hold and that $\tau=0$.  Let $u_0$ be bounded, nonnegative and uniformly continuous.  Suppose that for some $R_0\in\mathbb R$,
\begin{equation}\label{eq:filleddata}
 u_0(x)=0\;\text{if }x\cdot e\ge R_0,
 \;
 \liminf_{R\to\infty}\inf_{x\cdot e\le-R}u_0(x)>0.
\end{equation}
We set $s_c=\min\{c,c_*\}$.  Then for every $\eps>0$,
\begin{equation}
 \lim_{t\to\infty}\sup_{x\cdot e\ge(s_c+\eps)t}u(t,x)=0,\;
 \liminf_{t\to\infty}\inf_{x\cdot e\le(s_c-\eps)t}u(t,x)>0.
 \label{eq:filleddecay}
\end{equation}
If \eqref{eq:strongdamp} holds, then
\begin{equation}\label{eq:filledconv}
 \lim_{t\to\infty}\sup_{x\cdot e\le(s_c-\eps)t}
 \left(|u(t,x)-p_-(t)|+|v(t,x)-q_-(t)|\right)=0.
\end{equation}
\end{theorem}

Consider now a favorable slab bounded on both longitudinal sides by unfavorable regions.

\begin{assumption}[Localized favorable slab]\label{L}
The limits $r_\pm(t)=\lim_{\xi\to\pm\infty}r(t,\xi)$ exist uniformly in $t$ and
\(\langle r_-\rangle<0, \; \langle r_+\rangle<0,\)
while $\sup_{t,\xi}r(t,\xi)>0$.
\end{assumption}

For $L>0$ consider the one-dimensional periodic Dirichlet problem
\begin{equation}\label{eq:FloquetL}
\begin{cases}
 \phi_t=a(t)\phi_{\xi\xi}+c\phi_\xi+r(t,\xi)\phi-\lambda_L(c)\phi,
 &|\xi|<L,\\
 \phi(t,\pm L)=0,
 \; \phi(t+T,\xi)=\phi(t,\xi),
\end{cases}
\end{equation}
with $\phi>0$.  The principal exponent is increasing in $L$ and we set
\(\lambda_\infty(c)=\lim_{L\to\infty}\lambda_L(c).\)
Also put $r^{\mathrm{sup}}(t)=\sup_{\xi\in\mathbb R}r(t,\xi)$ and $\bar r^{\mathrm{sup}}=\langle r^{\mathrm{sup}}\rangle$.
When $\bar r^{\mathrm{sup}}>0$, define
\begin{equation}\label{eq:admissiblesup}
 \mathcal A_D=\{\lambda>0:\ \beta(t)\ge\gamma_{N,\lambda}(t)
 \text{ for every }t\in[0,T]\},
 \;
 c_D^{\rm env}=\inf_{\lambda\in\mathcal A_D}
 \left(\bar a\lambda+\frac{\bar r^{\mathrm{sup}}}{\lambda}\right),
\end{equation}
with the convention $c_D^{\rm env}=+\infty$ if $\mathcal A_D=\varnothing$.

For the negative-exponent result, let $\Phi>0$ be a generalized planar Floquet eigenfunction
\begin{equation}\label{eq:wholeFloquet}
 \Phi_t=a(t)\Phi_{\xi\xi}+c\Phi_\xi+r(t,\xi)\Phi-
 \lambda_\infty(c)\Phi,
 \; \Phi(t+T,\xi)=\Phi(t,\xi),
\end{equation}
normalized by $\Phi(0,0)=1$.  We introduce its logarithmic slope and the optional strict slope condition
\begin{equation}\label{eq:KPhi}
 K_\Phi:=\sup_{t,\xi}\frac{|\Phi_\xi(t,\xi)|}{\Phi(t,\xi)}<\kappa_*:=\min_t\kappa(t).
\end{equation}
For $0\le\rho<1$ define the weighted signed Yukawa functional
\begin{equation}\label{eq:PsiN}
 \Psi_N(\rho)=
 \sup_{|q|\le\rho}
 \int_{\mathbb R^N}
 \big[-q\,\p_eG_N(z)-G_N(z)\big]_+
 e^{\rho|z\cdot e|}\,dz.
\end{equation}
It is finite for $\rho<1$.  We define
\begin{equation}\label{eq:GammaPhi}
 \gamma_{N,\Phi}(t):=
 \gamma_0(t)
 \left(1+\Psi_N\!\left(\frac{K_\Phi}{\kappa(t)}\right)\right).
\end{equation}
For $N=1$ the same integral is finite at $\rho=1$ and is identically zero for $0\le\rho\le1$; accordingly we extend the definition by $\Psi_1(1):=0$. Indeed, $G_1(z)=e^{-|z|}/2$ and $|\partial_eG_1|=G_1$ away from the origin, so $[-q\partial_eG_1-G_1]_+=0$ whenever $|q|\le1$. The coefficient $\gamma_{N,\Phi}$ is invoked only in the admissible range $K_\Phi<\kappa_*$, or for $N=1$ in the extended endpoint range $K_\Phi\le\kappa_*$. Outside these ranges no finite weighted-Yukawa estimate is asserted.

\begin{theorem}[Moving favorable slab]\label{thm:localized}
Suppose that \cref{A1,A2,L} hold, that $\tau=0$, and that $\delta_0>0$.  Let $u_0\in BUC(\mathbb R^N)$ be nonnegative, nontrivial, and compactly supported.
\begin{enumerate}[label=\textup{(\roman*)},leftmargin=8mm]
\item For every $\sigma>0$,
\begin{equation}\label{eq:longitudinalconfinement}
 \lim_{t\to\infty}\sup_{|x\cdot e-ct|\ge\sigma t}u(t,x)=0.
\end{equation}
Thus the population is asymptotically confined, in the longitudinal direction, to every linearly expanding neighborhood of the moving favorable slab.
\item If $\bar r^{\mathrm{sup}}\le0$, then $\lVert u(t,\cdot)\rVert_\infty\to0$ for every $c$.  If $\bar r^{\mathrm{sup}}>0$ and $|c|>c_D^{\rm env}$, then again $\lVert u(t,\cdot)\rVert_\infty\to0$.
\item If $\lambda_\infty(c)>0$, then for every $R>0$,
\begin{equation}\label{eq:localpersistN}
 \liminf_{t\to\infty}
 \inf_{|x-(ct)e|\le R}u(t,x)>0.
\end{equation}
\item Suppose $\lambda_\infty(c)<0$ and either \eqref{eq:KPhi} holds, or $N=1$ and $K_\Phi\le\kappa_*$.  Suppose in addition that
\begin{equation}\label{eq:localizeddamping}
 \beta(t)\ge\gamma_{N,\Phi}(t)
 \;\text{for every }t\in[0,T].
\end{equation}
Then $\lVert u(t,\cdot)\rVert_\infty\to0$.  Moreover $u(t,x+ct\,e)$ decays exponentially on every fixed compact subset of $\mathbb R^N$.
\end{enumerate}
No conclusion is asserted when $\lambda_\infty(c)=0$.
\end{theorem}

\begin{proposition}[Exact climate-speed dependence of the periodic slab exponent]\label{prop:FloquetSpeed}
Suppose that \cref{A1,A2,L} hold and that $a(t)\equiv a_0>0$. Then for every $L>0$ and $c\in\mathbb R$,
\begin{equation}\label{eq:FloquetExactShift}
 \lambda_L(c)=\lambda_L(0)-\frac{c^2}{4a_0},\qquad
 \lambda_\infty(c)=\lambda_\infty(0)-\frac{c^2}{4a_0}.
\end{equation}
If $\lambda_\infty(0)>0$, put
\(c_{\mathrm{tr}}=2\sqrt{a_0\lambda_\infty(0)}.\)
Then $\lambda_\infty(c)>0$ exactly when $|c|<c_{\mathrm{tr}}$, and $\lambda_\infty(c)<0$ exactly when $|c|>c_{\mathrm{tr}}$. In the first regime \cref{thm:localized}(iii) implies nonlinear persistence, whereas in the second regime \cref{thm:localized}(iv) implies extinction only if its additional eigenfunction-slope and Yukawa-damping hypotheses are fulfilled. No nonlinear conclusion follows solely from \eqref{eq:FloquetExactShift} when $\lambda_\infty(c)=0$.
\end{proposition}

\begin{assumption}[Wave regularity and monotone shifting profile]\label{M}
Assumption~\ref{F} holds.  In addition, the map $f$ is nonincreasing in the shifting coordinate,
\(f_\xi(t,\xi,u)\le0.\)
\emph{For the scalar-wave theorem invoked below, we additionally require}
\(f\in C^1(\mathbb R\times\mathbb R\times[0,\infty))\)
\emph{in the joint variables $(t,\xi,u)$, with the derivatives extending continuously to $u=0$.  On each strip $0\le u\le M$, the first derivatives are bounded uniformly in $(t,\xi)$; in particular $f_t$ exists and is continuous.  This additional joint $C^1$ hypothesis is imposed only for the wave results, not for the parabolic--elliptic Cauchy propagation theorems.}
For every $M>0$, the functions $f$, $f_u$, $f_{uu}$, $f_{uuu}$ and $f_\xi$ are bounded and uniformly H\"older continuous on $\mathbb R\times\mathbb R\times[0,M]$ with the time--space exponents required by the periodic Schauder theory. In addition, the derivatives in $u$ through order three are locally Lipschitz in $u$, as maps with values in $C_{\per}^{\alpha_w/2,\alpha_w}$ in $(t,\xi)$ on the fixed bounded population range used below, with bounds uniform in spatial translations. Here $\alpha_w\in(0,\alpha)$ is the exponent fixed in the wave section; equivalently we assume the corresponding mixed H\"older composition estimates for these derivatives.  Moreover, there are $\sigma_f>0$ and constants $C_M$ such that, for $j=0,1,2$,
\[
 |\partial_u^j f(t,\xi,u)-\partial_u^j f_\pm(t,u)|
 +|f_\xi(t,\xi,u)|
 \le C_M e^{-\sigma_f|\xi|}
\]
for $\pm\xi\ge1$, uniformly in $t\in\mathbb R$ and $0\le u\le M$.  Finally, for the wave analysis only, there is $u_{\rm ext}>0$ such that $f$ admits a $T$-periodic extension to $u\in(-u_{\rm ext},\infty)$ having the same $C^3$ regularity in $u$, the same uniform H\"older and mixed-derivative bounds on bounded $u$-intervals, and $f_\xi\le0$.  The extension has limits $f_\pm(t,u)$ for $-u_{\rm ext}<u<0$ and satisfies the same exponential far-field estimate above, uniformly on every compact interval $u\in[-u_0,M]\subset(-u_{\rm ext},\infty)$.  The extension permits sign-changing perturbations near the zero end state in the implicit-function argument.  The constructed waves are positive and therefore solve the original equation.  These additional hypotheses are used only for the weighted Fredholm and implicit-function arguments in the forced-wave section.
\end{assumption}

For $\tau>0$, let $q_{-,\tau}$ be the unique positive $T$-periodic solution of
\(\tau q'=-\nu(t)q+\mu(t)p_-(t),\)
and set $q_{-,0}=q_-$.

\begin{definition}
A planar forced $T$-periodic wave of \eqref{eq:fullmodel} is a solution
\((u,v)(t,x)=(U,V)(t,x\cdot e-ct), \; (U,V)(t+T,\xi)=(U,V)(t,\xi),\)
with
\begin{equation}\label{eq:waveends}
 (U,V)(t,-\infty)=(p_-(t),q_{-,\tau}(t)),
 \;
 (U,V)(t,+\infty)=(0,0)
\end{equation}
uniformly in $t\in[0,T]$.
\end{definition}

\begin{theorem}[Weak-taxis planar forced waves]\label{thm:waves}
Suppose that \cref{A1,A2,M} hold, including the joint $C^1$ regularity in $(t,\xi,u)$ required for waves, and that $a(t)\equiv a_0>0$.  Let $\tau\ge0$ be fixed, and let $I\Subset(-\infty,c_*)$ be compact, where $c_*=2\sqrt{a_0\bar r}$.  Let $U_c$ be the zero-taxis scalar forced wave.  Suppose that the nondegenerate pinning condition
\begin{equation}\label{eq:pinning}
 -U_c(t,\xi)f_\xi(t,\xi,U_c(t,\xi))\not\equiv0
 \;\text{on }[0,T]\times\mathbb R
 \;\text{for every }c\in I.
\end{equation}
Then there exist $\eta_{I,\tau}>0$ and an admissible two-sided exponential weight such that for every $c\in I$ and $0\le\eta\le\eta_{I,\tau}$, system \eqref{eq:fullmodel} admits a planar forced $T$-periodic wave $(U_{\eta,\tau,c},V_{\eta,\tau,c})$ satisfying \eqref{eq:waveends}.  The wave is locally unique among solutions with the same end states in the specified planar weighted class, in a neighborhood of the zero-taxis family. The neighborhood and inverse estimates are uniform for $c\in I$, for the fixed value of $\tau$ under consideration.  Moreover,
\begin{equation}\label{eq:etatends}
 (U_{\eta,\tau,c},V_{\eta,\tau,c})
 \longrightarrow(U_{0,\tau,c},V_{0,\tau,c})
 \;\text{as }\eta\downarrow0
\end{equation}
uniformly for $c\in I$ in the corresponding weighted H\"older topology, and the convergence rate is $O(\eta)$ uniformly on $I$.  The statement at $\tau=0$ uses the elliptic signal space; no uniform implicit-function theorem across the singular limit $\tau\downarrow0$ is asserted.
\end{theorem}

\begin{theorem}[First-order chemotactic response of pinned periodic waves]\label{thm:linearresponse}
Assume all hypotheses of \cref{thm:waves}, fix $\tau\ge0$, and let $I\Subset(-\infty,c_*)$ be the same compact speed interval. For each $c\in I$, write $V_{\tau,c}=V_{0,\tau,c}$ for the zero-taxis signal and define $(H_{\tau,c},K_{\tau,c})$ as the unique periodic weighted solution of
\begin{equation}\label{eq:responseEq}
 \begin{cases}
 \CalA_cH_{\tau,c}=-\chi(t)\big(U_c(V_{\tau,c})_\xi\big)_\xi,\\
 \CalB_{\tau,c}K_{\tau,c}=\mu(t)H_{\tau,c},
 \end{cases}
\end{equation}
where $\CalA_c$ and $\CalB_{\tau,c}$ are the invertible scalar operators from \eqref{eq:triangularIFT}; for $\tau=0$ the second equation is elliptic. Then, uniformly in $c\in I$ as $\eta\downarrow0$,
\begin{equation}\label{eq:responseExpansion}
 \left\|
 \begin{pmatrix}U_{\eta,\tau,c}-U_c\\ V_{\eta,\tau,c}-V_{\tau,c}\end{pmatrix}
 -\eta\begin{pmatrix}H_{\tau,c}\\ K_{\tau,c}\end{pmatrix}
 \right\|_{\mathbb X_\tau}
 \le C_{I,\tau}\eta^2.
\end{equation}
In particular, the right derivative of the locally unique pinned branch with respect to $\eta$ at $\eta=0$ exists in the weighted H\"older topology and equals $(H_{\tau,c},K_{\tau,c})$. No assertion is made about the sign of either first-order correction or about stability under transverse perturbations.
\end{theorem}

We begin with the estimates used in the propagation and wave arguments.  The Keller--Segel evolution is not order preserving, so the exponential super-solution arguments used for scalar shifting equations such as \cite{BF2018,FPZ2021} cannot be applied directly to two chemotactic solutions.  The upper-bound argument therefore works only at a hypothetical first contact.  The exact resolvent identity in \cref{lem:expGreen} integrates the positive part of the kernel combination generated by the planar exponential envelope and produces the coefficient \(\gamma_{N,\lambda}(t)=\gamma_0(t)[1+\Theta_N(\lambda/\kappa(t))]\) defined above.  The coefficient is determined by the multidimensional Yukawa kernel.  Choosing only the minimizing KPP exponent \(\lambda_*\) closes the estimate under \eqref{eq:Dass}, without a supremum over exponential rates; in one dimension \eqref{eq:Theta1} recovers exactly the coefficient used in the shifting Keller--Segel analysis of Shen--Xue \cite{SX2020JDE}.  The resulting argument preserves the linearly determined speed \(c_*\) although the chemotactic drift is nonlinear and self-generated.

Second, persistence cannot be obtained by reversing the upper comparison.  Pointwise smallness of \(u\) does not by itself control \(\nabla v\) on a neighborhood.  The power-Harnack estimate in \cref{lem:Harnack} converts a small value of \(u\) into quantitative smallness on a fixed ball, and \cref{lem:localgreen,cor:smalldrift} converts that local information into smallness of the Yukawa drift.  We then truncate the drift, construct a positive entire scalar state from a periodic Dirichlet Floquet mode in \cref{lem:Dirichletpositive,lem:entireaux}, and compare only at a first contact where \eqref{eq:driftsmallchoice} proves that the truncation is inactive.  The lower bounds in \cref{thm:front,thm:frontfilled} and the positive-exponent assertion of \cref{thm:localized} follow from this comparison at first contact; a comparison principle between two arbitrary Keller--Segel solutions is not assumed.

Third, temporal signal memory creates a different obstruction in \cref{thm:waves}.  For \(\tau>0\), the signal is no longer an instantaneous elliptic resolvent, so the wave problem cannot be reduced to the parabolic--elliptic forced-wave construction of Shen--Xue \cite{SX2022FW}.  We work instead in affine two-sided weighted H\"older spaces with vanishing spatial jets.  The signal remainder is placed in these spaces by the whole-past Gaussian inverse and weighted Schauder estimate in \cref{lem:signalweighted}.  For the population linearization, the differentiated scalar profile \(W_c=-U_{c,\xi}>0\) satisfies the strict supersolution identity \eqref{eq:Wsupersolution}; the pinning condition \eqref{eq:pinning} and positivity improvement then exclude the Floquet multiplier \(1\).  Proposition~\ref{prop:fredholm} supplies an index-zero monodromy problem from the two far-field exponential gaps.  The pinning argument then excludes the multiplier $1$, and periodic extension and interior Schauder estimates upgrade monodromy invertibility to uniform invertibility of the parabolic block in \cref{lem:invert}.  The quantitative contraction in the proof of \cref{thm:waves} continues the scalar profile of Fang--Peng--Zhao \cite{FPZ2021} to weak taxis with an \(O(\eta)\) bound, for each fixed \(\tau\ge0\).  The first-contact estimate is used for propagation; the weighted periodic linearization is used separately for the wave continuation.

\subsection{Comparison with the related previous results}
The conclusions in \cref{thm:front,thm:frontfilled,thm:localized,thm:waves} address different mathematical questions, and the comparisons below are made under their respective hypotheses.

\emph{Shifting Keller--Segel habitats.} Shen--Xue \cite{SX2020JDE} considered $N=1$ with constant diffusion, signal and logistic coefficients. Here, in the elliptic-signal case, the coefficients are time-periodic and the translating favorable set is a half-space in $\mathbb R^N$. The sharp spreading region is characterized by the simultaneous inequalities $|s|<c_*$ and $s\cdot e<c$, instead of an interval bounded by two one-dimensional fronts. The anisotropic geometry is coupled to a radial propagation constraint; it is not obtained by simply replacing the one-dimensional derivative by a Laplacian. At the level of the sufficient damping hypothesis, \eqref{eq:ThetaN} uses the full $N$-dimensional Yukawa kernel and \eqref{eq:Dass} only requires the test exponent $\lambda_*$. In $N=1$ the original coefficient restriction is recovered exactly by \eqref{eq:Theta1}, so no improvement over Shen--Xue's one-dimensional constant-coefficient threshold is asserted.

\emph{Spatial dimension and planar profiles.} A planar ansatz
$U(t,x\cdot e-ct)$ reduces the spatial wave equation to a single profile variable,
regardless of the ambient dimension $N$. Thus the mere embedding of a planar wave
in $\mathbb R^N$ is not asserted as a new multidimensional effect: the
higher-dimensional contribution of the present paper lies chiefly in the
propagation and persistence geometry and the dimension-dependent
Yukawa first-contact estimate. For the forced-wave theorem, the distinction from
prior theory is the time-periodic shifting habitat, the permitted chemical memory
$\tau>0$, and the weighted continuation under a nondegenerate pinning assumption.

\emph{Homogeneous chemotaxis fronts.} Salako--Shen \cite{SalakoShen2017DCDS,SalakoShen2022Fronts} and Shen--Xue \cite{ShenXue2022Persistence,SX2022PP} treat traveling fronts, spreading or convergence without a translating unfavorable far field; the parabolic--parabolic setting has also been addressed for homogeneous coefficients. In contrast, \cref{thm:waves} constructs a \emph{forced}, temporally periodic profile pinned to a shifting resource, under $a(t)\equiv a_0$ and for each \emph{fixed} $\tau\ge0$ at sufficiently small taxis strength. It establishes local uniqueness near the zero-taxis branch rather than uniqueness among all possible forced waves; no continuation uniform as $\tau\downarrow0$ is claimed. The proof uses the far-field spectral gaps and the weighted periodic linearization at the pinned scalar profile; these features distinguish the continuation argument from the cited autonomous traveling-wave constructions.

\emph{Time-periodic scalar and cooperative propagation.} Fang--Peng--Zhao \cite{FPZ2021} established the scalar periodic forced-wave theorem. Shen \cite{Shen2010Variational,Shen2011Monostable} and Du--Li--Shen \cite{DuLiShen2022JFA} developed important time-dependent scalar or monotone-semiflow propagation theories. Their order-preserving arguments cannot be transferred unchanged to the chemotactic drift. For the parabolic--elliptic problem, we estimate the Yukawa term at first contact and obtain the lower bounds through power-Harnack localization. The localized-habitat result is governed by a periodic Dirichlet Floquet exponent, with an additional weighted-resolvent assumption in the negative-exponent case. The sharp kernel asymptotics and first-order wave response refine these comparisons quantitatively. The constant-diffusion Floquet drift identity is elementary at the spectral level and is not claimed as a novel standalone theorem. Neither the zero-exponent threshold nor general variable-$a(t)$ forced waves with $\tau>0$ are resolved here.
\Needspace{8\baselineskip}
\noindent\textbf{Organization of the paper.}  Section~\ref{sec:green} develops the multidimensional Yukawa estimates, including sharp small-slope asymptotics, underlying the comparison-free arguments.  Section~\ref{sec:global} establishes global boundedness and the power-Harnack localization.  Sections~\ref{sec:envelope} and \ref{sec:frontproof} prove the sharp geometric propagation results.  Section~\ref{sec:localized} treats moving favorable slabs and the exact spectral climate-speed law for constant diffusivity. Section~\ref{sec:waves} constructs planar forced waves for each fixed signal-relaxation parameter $\tau\ge0$ and derives their first-order chemotactic correction.  Section~\ref{sec:recover} treats the one-dimensional reductions, and Section~\ref{sec:discussion} summarizes the scope and remaining open problems.  Appendix~\ref{sec:appendix} collects the periodic principal-exponent and weighted Fredholm facts used in the main text.

\section{The elliptic signal and multidimensional resolvent estimates}\label{sec:green}
Throughout Sections~\ref{sec:green}--\ref{sec:localized} we set $\tau=0$.  At each fixed time, the signal equation is
\begin{equation}\label{eq:ellipticsignal}
 d(t)\Delta v-\nu(t)v+\mu(t)u=0
 \;\text{in }\mathbb R^N.
\end{equation}
The following elementary facts about the Yukawa kernel will be used repeatedly.

\begin{lemma}[Resolvent kernel estimates]\label{lem:YukawaFacts}
Let $G_{N,\kappa}$ be defined by \eqref{eq:Yukawaheat}.  Then
\begin{align}
 &G_{N,\kappa}(z)=\kappa^{N-2}G_N(\kappa z),
 \;
 \int_{\mathbb R^N}G_{N,\kappa}(z)\,dz=\kappa^{-2},
 \label{eq:Yukawascale}\\
 &\int_{e^\perp}G_{N,\kappa}(re+y)\,dy
 =\frac{1}{2\kappa}e^{-\kappa|r|},
 \; r\in\mathbb R.
 \label{eq:Yukawamarginal}
\end{align}
Moreover $G_N$ is radial and strictly decreasing in $|z|$, $\Theta_N$ in \eqref{eq:ThetaN} is finite and
\begin{equation}\label{eq:ThetaBound}
 0\le\Theta_N(\rho)\le\frac\rho2
 \;(\rho\ge0).
\end{equation}
For $N=1$, \eqref{eq:Theta1} holds.  The functional $\Psi_N$ in \eqref{eq:PsiN} is finite for every $0\le\rho<1$; when $N=1$ it extends to the endpoint $\rho=1$ and
\begin{equation}\label{eq:Psi1zero}
 \Psi_1(\rho)=0
 \;(0\le\rho\le1).
\end{equation}
\end{lemma}

\begin{proof}
The scaling in \eqref{eq:Yukawascale} follows from the change of variables $s=\sigma/\kappa^2$ in \eqref{eq:Yukawaheat}.  We integrate the heat kernel in $z$ and obtain
\(\int_{\mathbb R^N}G_{N,\kappa}(z)\,dz =\int_0^\infty e^{-\kappa^2s}\,ds=\kappa^{-2}.\)
Write $z=re+y$, where $y\in e^\perp$. Integration of the Gaussian over $e^\perp$ yields
\[
 \int_{e^\perp}G_{N,\kappa}(re+y)\,dy
 =\int_0^\infty (4\pi s)^{-1/2}\exp\!\left(-\frac{r^2}{4s}-\kappa^2s\right)\,ds
 =\frac{1}{2\kappa}e^{-\kappa|r|},
\]
which is \eqref{eq:Yukawamarginal}.

The heat representation shows that $G_N$ is radial.  Differentiation under the integral shows that $\p_rG_N<0$ for $r>0$.  Hence $\p_eG_N(z)>0$ can occur only when $z\cdot e<0$.  On the positive set of the integrand in \eqref{eq:ThetaN}, therefore, $e^{\rho z\cdot e}\le1$ and
\([\rho\p_eG_N-G_N]_+\le \rho(\p_eG_N)_+.\)
We use the divergence theorem on the half-space $\{z\cdot e<0\}$ together with \eqref{eq:Yukawamarginal} at $r=0$ and obtain
\begin{align*}
 \int_{z\cdot e<0}\p_eG_N(z)\,dz
 &=\int_{e^\perp}G_N(y)\,dy
 =\frac12.
\end{align*}
This proves \eqref{eq:ThetaBound}.

For $N=1$, $G_1(z)=\frac12e^{-|z|}$.  The positive part of
$\rho G_1'-G_1$ vanishes when $\rho\le1$.  If $\rho>1$, it is supported on $z<0$ and equals
\(\frac{\rho-1}{2}e^z.\)
Consequently
\[
 \Theta_1(\rho)
 =\frac{\rho-1}{2}\int_{-\infty}^0e^{(1+\rho)z}\,dz
 =\frac12\frac{\rho-1}{\rho+1},
\]
which proves \eqref{eq:Theta1}.

The weighted integrability required for $\Psi_N$ follows directly from \eqref{eq:Yukawaheat}.  Let $0\le\rho<1$ and choose $\varepsilon$ with $\rho^2<\varepsilon<1$.  Young's inequality implies
\(\rho|z|\le \varepsilon|z|^2/(4s)+\rho^2s/\varepsilon.\)
Hence the heat representation and Tonelli's theorem imply
\[
 \int_{\mathbb R^N}G_N(z)e^{\rho|z|}\,dz
 \le C_{N,\varepsilon}\int_0^\infty
 e^{-(1-\rho^2/\varepsilon)s}\,ds<\infty.
\]
Differentiating the heat kernel under the $s$-integral and using the same inequality yields
\[
 \int_{\mathbb R^N}|\nabla G_N(z)|e^{\rho|z|}\,dz
 \le C_{N,\varepsilon}\int_0^\infty
 s^{-1/2}e^{-(1-\rho^2/\varepsilon)s}\,ds<\infty.
\]
The factor $s^{-1/2}$ is integrable at the origin, and the positive exponential margin is integrable at infinity.  Since $|z\cdot e|\le|z|$, the integral in \eqref{eq:PsiN} is finite for every $\rho<1$.  Finally, for $N=1$ and $|q|\le1$,
\(-qG_1'(z)-G_1(z) =\big(q\,\operatorname{sgn}z-1\big)G_1(z)\le0,\)
so \eqref{eq:Psi1zero} follows.
\end{proof}

\begin{proposition}[Sharp small-slope Yukawa asymptotics]\label{prop:YukawaSmallSlope}
For every $N\ge2$ the functional $\Theta_N$ is strictly positive on $(0,\infty)$. Moreover, as $\rho\downarrow0$,
\begin{equation}\label{eq:YukawaSharpSmall}
 \Theta_2(\rho)=\frac{\rho^2}{8\log(1/\rho)}\big(1+o(1)\big),
 \qquad
 \Theta_N(\rho)=\frac{N-2}{4N}\rho^2\big(1+o(1)\big)
 \quad(N\ge3).
\end{equation}
By contrast, $\Theta_1(\rho)=0$ for $0\le\rho\le1$. In particular, if $N\ge2$ and $\kappa(t)\gg\lambda_*$, the coefficient obtained from the exact first-contact estimate has the expansion
\begin{equation}\label{eq:asympGamma}
 \frac{\gamma_{N,\lambda_*}(t)}{\gamma_0(t)}
 =1+\begin{cases}
 \displaystyle\frac{(\lambda_*/\kappa(t))^2}{8\log(\kappa(t)/\lambda_*)}\big(1+o(1)\big),&N=2,\\[2mm]
 \displaystyle\frac{N-2}{4N}(\lambda_*/\kappa(t))^2\big(1+o(1)\big),&N\ge3,
 \end{cases}
\end{equation}
where the expansion concerns $\lambda_*/\kappa(t)\to0$. The ratio in \eqref{eq:asympGamma} is interpreted only at times with $\gamma_0(t)>0$; equivalently its right-hand side defines the correction in $\gamma_{N,\lambda_*}=\gamma_0(1+\Theta_N)$ without division.
\end{proposition}

\begin{proof}
Write $G_N(z)=g_N(r)$, $r=|z|$, and $\omega=z/r$. Since $g_N'(r)<0$, the positive part in \eqref{eq:ThetaN} is supported on $\omega\cdot e<0$ and takes the form
\begin{equation}\label{eq:polarTheta}
 \big[\rho(-\omega\cdot e)(-g_N'(r))-g_N(r)\big]_+.
\end{equation}
For $N\ge2$, the logarithmic derivative $-g_N'(r)/g_N(r)$ tends to $+\infty$ as $r\downarrow0$. Hence, for every $\rho>0$, the expression in \eqref{eq:polarTheta} is positive on a nonempty open subset of $\{r>0,\omega\cdot e<0\}$. This proves $\Theta_N(\rho)>0$.

For $N\ge3$, set $\sigma_{N-1}=|\mathbb S^{N-1}|$. The fundamental-solution singularity of $-\Delta+1$, which also follows by splitting the heat representation \eqref{eq:Yukawaheat} at $s=1$, is
\begin{equation}\label{eq:YukawaSingN}
 g_N(r)=\frac{r^{2-N}}{(N-2)\sigma_{N-1}}\big(1+o(1)\big),\qquad
 -g_N'(r)=\frac{r^{1-N}}{\sigma_{N-1}}\big(1+o(1)\big)
 \quad(r\downarrow0).
\end{equation}
For every fixed $r_0>0$, the ratio $-g_N'/g_N$ is bounded on $[r_0,\infty)$ (the large-$r$ limit is $1$). For small $r$, \eqref{eq:YukawaSingN} shows that $-g_N'/g_N\le C/r$. Thus for all sufficiently small $\rho$ the positive part in \eqref{eq:polarTheta} is supported in $0<r<C\rho$. Put $r=\rho s$ and use polar coordinates. On this bounded interval of $s$ the integrand, after multiplication by $r^{N-1}\,dr/\rho^2$, converges pointwise to
\[
 \frac1{\sigma_{N-1}}
 \Big[-\omega\cdot e-\frac{s}{N-2}\Big]_+\,ds,
\]
with a uniform integrable bound. Also $e^{\rho r\omega\cdot e}\to1$ uniformly there. Dominated convergence yields
\begin{align*}
 \lim_{\rho\downarrow0}\frac{\Theta_N(\rho)}{\rho^2}
 &=\frac1{\sigma_{N-1}}\int_{\omega\cdot e<0}
   \int_0^{(N-2)(-\omega\cdot e)}
   \left(-\omega\cdot e-\frac{s}{N-2}\right)ds\,d\omega\\
 &=\frac{N-2}{2\sigma_{N-1}}\int_{\omega\cdot e<0}(\omega\cdot e)^2\,d\omega
 =\frac{N-2}{4N},
\end{align*}
where the last identity follows from spherical symmetry and
$\int_{\mathbb S^{N-1}}(\omega\cdot e)^2\,d\omega=\sigma_{N-1}/N$.

For $N=2$ we use the explicit radial kernel $g_2(r)=(2\pi)^{-1}K_0(r)$ and its derivative $g_2'(r)=-(2\pi)^{-1}K_1(r)$. The small-argument formulas are
\begin{equation}\label{eq:BesselSmall}
 K_0(r)=\log(1/r)+O(1),\qquad rK_1(r)=1+o(1)
 \quad(r\downarrow0).
\end{equation}
Consequently, $-g_2'/g_2\le C/[r\log(1/r)]$ for sufficiently small $r$. As in the preceding case, the ratio is bounded away from $r=0$, and therefore the positive part vanishes outside a small neighborhood of the origin for small $\rho$. Put $L=\log(1/\rho)$. For any $M>C$, the inequality $r\log(1/r)>C\rho$ whenever $M\rho/L\le r\le r_0$ and $\rho$ is sufficiently small (after increasing $M$ if necessary) shows that the positive part is supported in $0<r<M\rho/L$.

Writing $r=\rho s/L$ and using \eqref{eq:BesselSmall}, we obtain for each $s>0$
\[
 \frac{\rho(-\omega\cdot e)K_1(\rho s/L)-K_0(\rho s/L)}{L}
 \longrightarrow \frac{-\omega\cdot e}{s}-1.
\]
On $0<s<M$ the corresponding polar integral has an integrable majorant after multiplication by the Jacobian $s\,ds$; it follows directly from $rK_1(r)\le C$ and $K_0(r)\le C(1+|\log r|)$. Since $e^{\rho r\omega\cdot e}\to1$ uniformly on the support, dominated convergence yields
\begin{align*}
 \lim_{\rho\downarrow0}\frac{L\Theta_2(\rho)}{\rho^2}
 &=\frac1{2\pi}\int_{\omega\cdot e<0}\int_0^{-\omega\cdot e}
 \big[-\omega\cdot e-s\big]\,ds\,d\omega\\
 &=\frac1{4\pi}\int_{\omega\cdot e<0}(\omega\cdot e)^2\,d\omega
 =\frac18.
\end{align*}
This proves \eqref{eq:YukawaSharpSmall}. The one-dimensional statement is exactly \eqref{eq:Theta1}, and \eqref{eq:asympGamma} follows by substitution in the definition of $\gamma_{N,\lambda_*}$.
\end{proof}

The elliptic signal can now be written as
\begin{equation}\label{eq:GreenN}
 v(t,x)=\frac{\mu(t)}{d(t)}
 \int_{\mathbb R^N}G_{N,\kappa(t)}(x-y)u(t,y)\,dy.
\end{equation}
In particular,
\begin{equation}\label{eq:basicgreenN}
 0\le v(t,x)\le\frac{\mu(t)}{\nu(t)}\lVert u(t,\cdot)\rVert_\infty.
\end{equation}
If
\(C_N^\nabla:=\int_{\mathbb R^N}|\nabla G_N(z)|\,dz,\)
scaling then yields
\begin{equation}\label{eq:gradgreenN}
 |\nabla v(t,x)|
 \le C_N^\nabla\frac{\mu(t)}{\sqrt{d(t)\nu(t)}}
 \lVert u(t,\cdot)\rVert_\infty.
\end{equation}

The next lemma is the multidimensional replacement for the one-dimensional derivative-kernel identity.

\begin{lemma}[Planar exponential first-contact estimate]\label{lem:expGreen}
We fix $t$, a unit vector $\theta\in\mathbb S^{N-1}$, $\lambda>0$, and let
\(W(x)=A e^{-\lambda x\cdot\theta}, \; A>0.\)
Suppose $0\le u\le W$ and $v$ solves \eqref{eq:ellipticsignal}.  At a contact point $x_0$ at which
\(u(x_0)=W(x_0), \; \nabla u(x_0)=\nabla W(x_0)=-\lambda W(x_0)\theta,\)
one has
\begin{equation}\label{eq:sharpGreenineqN}
 -\eta\chi(t)\nabla\!\cdot(u\nabla v)(x_0)
 \le
 \frac{\eta\chi(t)\mu(t)}{d(t)}
 \left(1+\Theta_N\!\left(\frac{\lambda}{\kappa(t)}\right)\right)
 W(x_0)^2.
\end{equation}
\end{lemma}

\begin{proof}
We fix the time and abbreviate $d=d(t)$, $\nu=\nu(t)$, $\mu=\mu(t)$ and $\kappa=\sqrt{\nu/d}$.  From \eqref{eq:ellipticsignal},
\begin{equation}\label{eq:Deltavidentity}
 \Delta v=\kappa^2v-\frac\mu d u.
\end{equation}
At the contact point,
\begin{equation}\label{eq:contactexpandN}
 -\eta\chi\nabla\!\cdot(u\nabla v)
 =-\eta\chi\nabla u\cdot\nabla v-\eta\chi u\Delta v
 =\eta\chi\lambda W\,\p_\theta v-\eta\chi\kappa^2Wv+\eta\chi\frac\mu dW^2.
\end{equation}
Thus only $\lambda\p_\theta v-\kappa^2v$ has to be estimated.  By \eqref{eq:GreenN}, with $z=x_0-y$,
\begin{align}
 \lambda\p_\theta v(x_0)-\kappa^2v(x_0)
 =\frac\mu d\int_{\mathbb R^N}
 \big(\lambda\p_\theta G_{N,\kappa}(z)-\kappa^2G_{N,\kappa}(z)\big)
 u(x_0-z)\,dz.
 \label{eq:kernelcomboN}
\end{align}
Because $G_{N,\kappa}$ is radial and strictly decreasing in $|z|$, its directional derivative satisfies
\[
 \p_\theta G_{N,\kappa}(z)
 =G_{N,\kappa}'(|z|)\frac{z\cdot\theta}{|z|}\le0
 \;\hbox{when }z\cdot\theta\ge0.
\]
Consequently the positive part of the kernel in \eqref{eq:kernelcomboN} is contained in $\{z\cdot\theta<0\}$.  On that set the global exponential envelope satisfies the exact relation
\begin{equation}\label{eq:enveloperatioN}
 u(x_0-z)\le W(x_0-z)=W(x_0)e^{\lambda z\cdot\theta}.
\end{equation}
We drop the negative part of the kernel and use \eqref{eq:enveloperatioN}.  We obtain
\begin{align*}
 \lambda\p_\theta v(x_0)-\kappa^2v(x_0)
 &\le\frac\mu dW(x_0)
 \int_{\mathbb R^N}
 \big[\lambda\p_\theta G_{N,\kappa}(z)-\kappa^2G_{N,\kappa}(z)\big]_+
 e^{\lambda z\cdot\theta}\,dz.
\end{align*}
We set $w=\kappa z$.  The identities
\[
 G_{N,\kappa}(z)=\kappa^{N-2}G_N(w),
 \;
 \p_\theta G_{N,\kappa}(z)=\kappa^{N-1}\p_\theta G_N(w),
 \; dz=\kappa^{-N}dw
\]
transform the last integral into
\[
 \int_{\mathbb R^N}
 \left[\frac{\lambda}{\kappa}\p_\theta G_N(w)-G_N(w)\right]_+
 e^{(\lambda/\kappa)w\cdot\theta}\,dw
 =\Theta_N\!\left(\frac\lambda\kappa\right).
\]
Substituting into \eqref{eq:contactexpandN}, we obtain \eqref{eq:sharpGreenineqN}.
\end{proof}

We also need a local form of the resolvent estimates.

\begin{lemma}[Local signal localization]\label{lem:localgreen}
Suppose that $0\le u\le M$.  For every $R>1$ and $x_0\in\mathbb R^N$,
\begin{align}
 |\nabla v(t,x_0)|+v(t,x_0)
 &\le C_0\sup_{|y-x_0|\le R}u(t,y)+\varepsilon_RM,
 \label{eq:localsignal1N}\\
 |\Delta v(t,x_0)|
 &\le C_1\sup_{|y-x_0|\le R}u(t,y)+\varepsilon_RM,
 \label{eq:localsignal2N}
\end{align}
where $C_0,C_1$ depend only on the coefficient bounds and $N$, while $\varepsilon_R\downarrow0$ as $R\to\infty$.
\end{lemma}

\begin{proof}
From \eqref{eq:GreenN}, split the convolution into $|z|\le R$ and $|z|>R$.  Because $\kappa(t)$ ranges in a compact subset of $(0,\infty)$, the $L^1$ norms of $G_{N,\kappa(t)}$ and $\nabla G_{N,\kappa(t)}$ are uniformly bounded.  The corresponding tail integrals converge to zero uniformly in $t$ by scaling and the exponential decay of the Yukawa kernel.  Hence
\(|\nabla v(t,x_0)|+v(t,x_0) \le C_0\sup_{B_R(x_0)}u(t,\cdot)+\varepsilon_RM.\)
Finally \eqref{eq:Deltavidentity} and the lower bound for $d$ imply \eqref{eq:localsignal2N} after changing the constants.
\end{proof}

\begin{lemma}[Weighted planar resolvent estimate]\label{lem:loglipN}
Let $W(x)=w(x\cdot e)>0$ with
\(|w'(\xi)|\le K w(\xi), \; 0\le K<\kappa(t),\)
and suppose $0\le u\le MW$.  If $N=1$, the endpoint $K=\kappa(t)$ is also allowed.  At a first contact $u=MW$, $\nabla u=Mw'(x\cdot e)e$, one has
\begin{equation}\label{eq:weightedchemN}
 -\eta\chi(t)\nabla\!\cdot(u\nabla v)
 \le
 \frac{\eta\chi(t)\mu(t)}{d(t)}
 \left(1+\Psi_N\!\left(\frac{K}{\kappa(t)}\right)\right)
 M^2W^2.
\end{equation}
\end{lemma}

\begin{proof}
The logarithmic derivative bound implies
\begin{equation}\label{eq:Wratio}
 W(x-z)\le W(x)e^{K|z\cdot e|}.
\end{equation}
At the contact, we put $q=w'(x\cdot e)/w(x\cdot e)$, so $|q|\le K$.  We then use \eqref{eq:Deltavidentity} and obtain
\begin{align*}
 -\eta\chi\nabla\!\cdot(u\nabla v)
 &= -\eta\chi M w'\p_ev
    -\eta\chi MW\kappa^2v
    +\eta\chi\frac\mu dM^2W^2.
\end{align*}
The first two terms are
\begin{align*}
 \frac{\eta\chi\mu}{d}MW
 \int_{\mathbb R^N}
 \big[-q\p_eG_{N,\kappa}(z)-\kappa^2G_{N,\kappa}(z)\big]
 u(x-z)\,dz.
\end{align*}
Retain the positive part and use \eqref{eq:Wratio}.  After the scaling $y=\kappa z$, the resulting integral is bounded by
\(M W\,\Psi_N\!\left(\frac{K}{\kappa}\right).\)
We add the last quadratic term and obtain \eqref{eq:weightedchemN}.
\end{proof}

\section{Global existence and a local power estimate}\label{sec:global}
The global bound follows from the sign of the chemotactic term at spatial maxima; this argument does not depend on the dimension.

\begin{proposition}[Global boundedness for the parabolic--elliptic problem]\label{prop:global}
Suppose that \cref{A1,A2} hold and that $\delta_0>0$.  For every bounded nonnegative uniformly continuous $u_0$, the parabolic--elliptic problem has a unique global classical solution for positive time.  Moreover,
\begin{equation}\label{eq:globalbound}
 0\le u(t,x)\le M_0:=
 \max\left\{\lVert u_0\rVert_\infty,
 \frac{r^*}{\delta_0}\right\},
 \;
 r^*:=\sup_{t,\xi}r(t,\xi)_+.
\end{equation}
The signal satisfies \eqref{eq:basicgreenN}--\eqref{eq:gradgreenN} uniformly in time.
\end{proposition}

\begin{proof}
We begin with local existence and verify that the existence time depends only on an $L^\infty$ bound, as required by the continuation argument.  Since $f$ is initially specified only for $u\ge0$, extend it once and for all to $u<0$ by
\(\widetilde f(t,\xi,u)=f(t,\xi,0)+u f_u(t,\xi,0), \; u<0,\)
and keep $\widetilde f=f$ for $u\ge0$.  The one-sided derivatives agree at $u=0$, so this defines a $C^1$ extension in $u$; the uniform local bounds in \cref{A2} imply that, on every interval $[-R,R]$, the extended Nemytskii map has a Lipschitz constant independent of $(t,\xi)$.  We suppress the tilde below.  The solution constructed below is subsequently shown to remain nonnegative, so the extension never changes the original problem.  For fixed $t$ let
\(\mathcal R_t w =\frac{\mu(t)}{d(t)}G_{N,\kappa(t)}*w.\)
By \cref{lem:YukawaFacts} and the positivity of the coefficient lower bounds,
\begin{equation}\label{eq:RtBounds}
 \|\mathcal R_tw\|_\infty+\|\nabla\mathcal R_tw\|_\infty
 \le C_R\|w\|_\infty
\end{equation}
uniformly in $t$.  We let $S_a(t,s)$ be the evolution generated by $a(t)\Delta$ on $BUC(\mathbb R^N)$.  If
\(A(t,s)=\int_s^ta(\zeta)\,d\zeta,\)
then $S_a(t,s)$ is convolution with the Gaussian of variance $2A(t,s)$, and therefore
\begin{equation}\label{eq:heatgradbound}
 \|\nabla S_a(t,s)F\|_\infty
 \le C(t-s)^{-1/2}\|F\|_\infty,
 \; 0\le s<t\le s+1,
\end{equation}
where $C$ depends only on $a_*$.  A mild solution must satisfy
\begin{align}
 u(t)
 &=S_a(t,0)u_0
 -\int_0^t\nabla S_a(t,s)\!\cdot
 \big(\eta\chi(s)u(s)\nabla\mathcal R_su(s)\big)\,ds \notag\\
 &\;+\int_0^tS_a(t,s)
 \big[u(s)f(s,\,\cdot\, e-cs,u(s))\big]ds.
 \label{eq:mildlocalN}
\end{align}
Here $f(s,\,\cdot\, e-cs,u)$ denotes the function $x\mapsto f(s,x\cdot e-cs,u(x))$.

We fix $R>2\|u_0\|_\infty+1$ and let $X_T$ be the closed ball of radius $R$ in $C([0,T];BUC(\mathbb R^N))$.  On $X_T$, \eqref{eq:RtBounds} provides
\[
 \|u\nabla\mathcal R_tu-\widetilde u\nabla\mathcal R_t\widetilde u\|_\infty
 \le \|u-\widetilde u\|_\infty\|\nabla\mathcal R_tu\|_\infty
 +\|\widetilde u\|_\infty\|\nabla\mathcal R_t(u-\widetilde u)\|_\infty
 \le C_RR\|u-\widetilde u\|_\infty.
\]
The local $C^1$ regularity of the extended reaction, together with boundedness of $r$, also implies
\[
 \|u f(t,x\cdot e-ct,u)-\widetilde u f(t,x\cdot e-ct,\widetilde u)\|_\infty
 \le C_R'\|u-\widetilde u\|_\infty
\]
whenever $|u|,|\widetilde u|\le R$ pointwise.  This is the estimate required on the whole Banach ball $X_T$, which is not restricted a priori to nonnegative functions.  Hence the right-hand side of \eqref{eq:mildlocalN} defines a map $\mathcal T:X_T\to C([0,T];BUC)$ satisfying
\[
 \|\mathcal Tu-\mathcal T\widetilde u\|_{X_T}
 \le C_R''(T^{1/2}+T)\|u-\widetilde u\|_{X_T}
\]
by \eqref{eq:heatgradbound}.  The same estimates, without differences, show that $\mathcal T(X_T)\subset X_T$ when $T>0$ is small.  Banach's fixed-point theorem therefore yields a unique mild solution.  Heat-kernel smoothing and elliptic Schauder estimates for $\mathcal R_t$ imply that this solution is $C^{1,2}$ for every positive time and satisfies the system classically.

Along this solution the population equation can be written
\begin{equation}\label{eq:nondivpositiveN}
 u_t=a(t)\Delta u-\eta\chi(t)\nabla v\cdot\nabla u
 +u\big(f(t,x\cdot e-ct,u)-\eta\chi(t)\Delta v\big).
\end{equation}
The coefficients are bounded on every compact time interval.  Approximation of $u_0$ by strictly positive bounded data and the parabolic maximum principle in \eqref{eq:nondivpositiveN}, followed by passage to the limit, provide $u\ge0$; if $u_0\not\equiv0$, the strong maximum principle provides $u(t,x)>0$ for $t>0$.

Let $T_{\max}$ be the maximal existence time. We derive an $L^\infty$ bound independent of $T_{\max}$.  We let
\(M(t)=\sup_{x\in\mathbb R^N}u(t,x).\)
When the supremum is not attained, maximize $u(t,x)-\delta|x|^2$.  If $x_\delta$ is a maximizing point, boundedness of $u$ provides $\delta|x_\delta|^2\le2\|u(t,\cdot)\|_\infty$, while
\(\nabla u(t,x_\delta)=2\delta x_\delta, \; \Delta u(t,x_\delta)\le2N\delta.\)
Consequently $|\nabla u(t,x_\delta)|\to0$ and the diffusion error tends to zero as $\delta\downarrow0$; the drift error also vanishes because $\nabla v$ is uniformly bounded.  We pass to a sequence of such penalized maxima and obtain the upper Dini derivative inequality below.  The identity
\(-\eta\chi\nabla\!\cdot(u\nabla v)=-\eta\chi\nabla u\cdot\nabla v-\eta\chi u\Delta v\),
together with \eqref{eq:Deltavidentity}, implies, after letting the penalization vanish,
\begin{align*}
 D^+M(t)
 &\le M(t)\left[r^* -\beta(t)M(t)\right]
 -\eta\chi(t)M(t)
 \left[\kappa(t)^2v-\frac{\mu(t)}{d(t)}M(t)\right]\\
 &\le M(t)\left[r^*-
 \left(\beta(t)-\gamma_0(t)\right)M(t)\right]\\
 &\le M(t)\left[r^*-\delta_0M(t)\right].
\end{align*}
Comparison with the logistic ordinary differential equation proves \eqref{eq:globalbound}.  In particular the radius $R$ in the local fixed-point construction can be chosen uniformly when the solution is restarted at any time $t_1<T_{\max}$.  The contraction time then depends only on $M_0$ and the fixed coefficient bounds, not on $t_1$.  If $T_{\max}<\infty$, restarting at times $t_1\uparrow T_{\max}$ would extend the mild solution beyond $T_{\max}$, a contradiction.  Hence $T_{\max}=\infty$.  The signal bounds are direct consequences of \cref{lem:YukawaFacts}, \eqref{eq:GreenN} and \eqref{eq:RtBounds}.
\end{proof}

The persistence proof uses a quantitative power-Harnack estimate.  Only one exponent is needed.  We therefore state the form that follows directly from two-sided Aronson estimates; this avoids any hidden compatibility condition between two unrelated Gaussian constants.

\begin{lemma}[Local power estimate]\label{lem:Harnack}
We fix $R>0$ and $t_0>0$.  Fix the taxis strength $\eta\ge0$ and a bound $M_0$ as in \eqref{eq:globalbound}. There exists $p_0>1$, depending on $N,R,t_0,M_0,\eta$ and the coefficient bounds in \cref{A1,A2}, such that, for every $p\ge p_0$, there is a constant $C>0$, depending on these data and on $p$, for which every solution bounded by \eqref{eq:globalbound} satisfies
\begin{equation}\label{eq:Harnack}
 \sup_{|y-x|\le R}u(t,y)
 \le C u(t,x)^{1/p},
 \; t\ge t_0,
 \; x\in\mathbb R^N.
\end{equation}
\end{lemma}

\begin{proof}
Along the fixed solution write
\begin{equation}\label{eq:scalarrewrite}
 u_t=a(t)\Delta u+B(t,x)\cdot\nabla u+C(t,x)u,
\end{equation}
where
\[
 B(t,x)=-\eta\chi(t)\nabla v(t,x),
 \;
 C(t,x)=f(t,x\cdot e-ct,u)-\eta\chi(t)\Delta v(t,x).
\]
By \cref{prop:global,lem:localgreen}, the coefficients are bounded uniformly on $[0,\infty)\times\mathbb R^N$ for the fixed value of $\eta$. More precisely, writing $M_0$ for the global population bound, \eqref{eq:basicgreenN}--\eqref{eq:gradgreenN} and \eqref{eq:Deltavidentity} imply
\[
\begin{aligned}
 \|B\|_\infty&\le \eta\chi^* C_N^\nabla
 \max_t\frac{\mu(t)}{\sqrt{d(t)\nu(t)}}M_0,\\
 \|C\|_\infty&\le\|r\|_\infty+b_1(1+M_0)M_0
 +2\eta\chi^*\max_t\frac{\mu(t)}{d(t)}M_0.
\end{aligned}
\]
The constants in the two-sided Gaussian estimates therefore depend on $\eta$ and $M_0$; they are not asserted to be uniform as the taxis strength varies. We fix
\(0<\delta<\min\{1,t_0/2\}\)
and let $\mathcal K(t,y;t-\delta,z)$ be the positive fundamental kernel of \eqref{eq:scalarrewrite}.  The two-sided Aronson bounds hold with constants $A_0,A_1,\alpha_0,\alpha_1>0$, independent of $t,x,y,z$, such that
\begin{equation}\label{eq:AronsonN}
\begin{aligned}
 A_0\delta^{-N/2}e^{-\alpha_0|x-z|^2/\delta}
 &\le \mathcal K(t,x;t-\delta,z),\\
 \mathcal K(t,y;t-\delta,z)
 &\le A_1\delta^{-N/2}e^{-\alpha_1|y-z|^2/\delta}.
\end{aligned}
\end{equation}
We choose $\theta\in(0,1)$.  If $|y-x|\le R$, then
\(|y-z|^2\ge \theta|x-z|^2-\frac{\theta}{1-\theta}R^2.\)
Hence
\begin{equation}\label{eq:kernelcompare1N}
 \mathcal K(t,y;t-\delta,z)
 \le C_R\exp\!\left(-\frac{\alpha_1\theta}{\delta}|x-z|^2\right).
\end{equation}
Choose $p_0>\max\{1,\alpha_0/(\alpha_1\theta)\}$ and set
\(\gamma:=\frac{1}{2\delta}\big(\alpha_1\theta-\alpha_0/p_0\big)>0.\)
For every $p\ge p_0$ we then have
\(\alpha_1\theta/\delta-\alpha_0/(p\delta)\ge2\gamma.\)
The lower bound in \eqref{eq:AronsonN} and \eqref{eq:kernelcompare1N} then imply
\begin{equation}\label{eq:kernelpowerN}
 \mathcal K(t,y;t-\delta,z)
 \le C_{R,p}\,\mathcal K(t,x;t-\delta,z)^{1/p}
 e^{-\gamma|x-z|^2}.
\end{equation}
We let $q=p/(p-1)$ and put $w(z)=u(t-\delta,z)$.  Since $0\le w\le M_0$, the representation formula, \eqref{eq:kernelpowerN}, and H\"older's inequality imply
\begin{align*}
 u(t,y)
 &\le C_{R,p}\int_{\mathbb R^N}
   \mathcal K(t,x;t-\delta,z)^{1/p}w(z)e^{-\gamma|x-z|^2}\,dz\\
 &=C_{R,p}\int_{\mathbb R^N}
   [\mathcal K(t,x;t-\delta,z)w(z)]^{1/p}
   w(z)^{1-1/p}e^{-\gamma|x-z|^2}\,dz\\
 &\le C_{R,p}
   \left(\int_{\mathbb R^N}\mathcal K(t,x;t-\delta,z)w(z)\,dz\right)^{1/p}
   \left(\int_{\mathbb R^N}w(z)e^{-q\gamma|x-z|^2}\,dz\right)^{1/q}\\
 &\le C_{R,p}u(t,x)^{1/p}
   \left(M_0\int_{\mathbb R^N}e^{-q\gamma|z|^2}\,dz\right)^{1/q}.
\end{align*}
The Gaussian integral is finite and independent of $t,x,y$.  Absorbing it into the constant proves \eqref{eq:Harnack}.
\end{proof}

\begin{corollary}[Smallness of the chemical drift at a small population point]\label{cor:smalldrift}
We fix $R>1$ and choose $p\ge p_0$ as in \cref{lem:Harnack} with $t_0=1$.  There are $C_{R,p}>0$ and $\varepsilon_R>0$, with $\varepsilon_R\downarrow0$ as $R\to\infty$, such that for $t\ge1$ the sum $\eta\chi(t)|\nabla v(t,x)|+\eta\chi(t)\kappa(t)^2v(t,x)$ is bounded by $C_{R,p}u(t,x)^{1/p}+\varepsilon_RM_0$.
\end{corollary}

\begin{proof}
We apply \cref{lem:localgreen} and then \cref{lem:Harnack} to the local supremum over $B_R(x)$. For the fixed problem, $\eta$, $M_0$ and the coefficient bounds are fixed; all constants are independent of the observation point and of $t\ge1$. The Yukawa tails in \cref{lem:localgreen} produce $\varepsilon_RM_0$, whereas the local part is bounded by $C_{R,p}u(t,x)^{1/p}$. The radius is chosen first, then the Harnack exponent, and finally the population threshold; neither $p$ nor the threshold is chosen before the localization radius.
\end{proof}

\section{Planar exponential envelopes and moving-frame growth}\label{sec:envelope}
We fix once and for all a H\"older exponent $\alpha_H\in(0,\alpha)$; all drift H\"older norms in Sections~\ref{sec:envelope}--\ref{sec:localized} are measured with the parabolic exponent pair $(\alpha_H/2,\alpha_H)$.  The favorable limiting linear equation is
\(w_t=a(t)\Delta w+r_-(t)w.\)
For a unit direction $\theta\in\mathbb S^{N-1}$ and $\lambda>0$, set
\(c_\lambda=\bar a\lambda+\frac{\bar r}{\lambda}\)
and
\(P_\lambda(t)= \exp\!\left(\int_0^t [a(s)\lambda^2+r_-(s)-\lambda c_\lambda]ds\right).\)
The integrand has zero temporal mean, hence $P_\lambda$ is positive and $T$-periodic.  Direct differentiation shows
\begin{equation}\label{eq:expidentityN}
 \p_t\left(P_\lambda e^{-\lambda(x\cdot\theta-c_\lambda t)}\right)
 =a(t)\Delta\left(P_\lambda e^{-\lambda(x\cdot\theta-c_\lambda t)}\right)
 +r_-(t)P_\lambda e^{-\lambda(x\cdot\theta-c_\lambda t)}.
\end{equation}
Minimizing $c_\lambda$, we obtain $\lambda_*$ and $c_*$ in \eqref{eq:lambda0cstar}.

\begin{lemma}[Planar exponential supersolutions in arbitrary directions]\label{lem:superenv}
Suppose that \cref{A1,A2,F,D} hold and let $u_0\in BUC(\mathbb R^N)$ be nonnegative and compactly supported. For every $s>c_*$ there exist $\lambda=\lambda_*$ and a speed $\widehat c$ satisfying
\begin{equation}\label{eq:chatchoiceN}
 c_\lambda<\widehat c<s,
\end{equation}
and constants $M=M(u_0,s)>0$, $z_0=z_0(u_0,s)\in\mathbb R$, independent of $\theta\in\mathbb S^{N-1}$, such that
\(\overline U_\theta(t,x) =M P_\lambda(t)e^{-\lambda(x\cdot\theta-\widehat c t-z_0)}\)
satisfies
\begin{equation}\label{eq:envelopecontrolN}
 u(t,x)\le\overline U_\theta(t,x)
 \;(t\ge0, x\in\mathbb R^N)
\end{equation}
for the solution issued from this fixed $u_0$.  No common $M,z_0$ is asserted for all initial data.
\end{lemma}

\begin{proof}
We take $\lambda=\lambda_*$.  Then $c_\lambda=c_*<s$, so choose $\widehat c$ with
\(c_*<\widehat c<s.\)
By \cref{D},
\begin{equation}\label{eq:GammausedN}
 \gamma_{N,\lambda_*}(t)\le\beta(t)
 \;(t\in\mathbb R).
\end{equation}
We let $R_0$ satisfy $\operatorname{supp}u_0\subset B_{R_0}(0)$.  Since $P_\lambda(0)>0$ and $x\cdot\theta\le R_0$ on the support, $M$ and $z_0$ can be chosen, independently of $\theta$, so that
\(u_0(x)<M P_\lambda(0)e^{-\lambda(x\cdot\theta-z_0)} \;(x\in\mathbb R^N).\)
Outside the compact support the inequality is automatic because the right-hand side is strictly positive.

We fix $\theta$.  The identity \eqref{eq:expidentityN} and the extra speed in \eqref{eq:chatchoiceN} imply
\begin{equation}\label{eq:strictExpIdentityN}
 \overline U_{\theta,t}-a(t)\Delta\overline U_\theta-r_-(t)\overline U_\theta
 =\lambda(\widehat c-c_\lambda)\overline U_\theta>0.
\end{equation}
We justify that a first contact cannot escape to spatial infinity.  On a fixed interval $0<t\le T_1$, equation \eqref{eq:scalarrewrite} has uniformly bounded drift and zero-order coefficient.  Its positive fundamental kernel satisfies the upper Aronson estimate
\(\mathcal K(t,x;0,y) \le C_{T_1}t^{-N/2} \exp\!\left(-\frac{|x-y|^2}{C_{T_1}t}\right).\)
Since $\operatorname{supp}u_0\subset B_{R_0}$,
\begin{equation}\label{eq:finitetimeGaussianN}
 u(t,x)
 \le C_{T_1}\|u_0\|_\infty
 \exp\!\left(-\frac{(|x|-R_0)_+^2}{C_{T_1}t}\right),
 \; 0<t\le T_1.
\end{equation}
Thus, uniformly on every compact positive time interval, $u$ has Gaussian decay as $x\cdot\theta\to+\infty$; as $x\cdot\theta\to-\infty$, the exponential envelope tends to $+\infty$ while $u$ remains bounded.  We start from the strict initial domination and enlarge $M$ by a factor $1+\varepsilon$.  Then \eqref{eq:finitetimeGaussianN} shows that any first zero of $(1+\varepsilon)\overline U_\theta-u$ lies in a bounded spatial set.  After the first-contact contradiction, we let $\varepsilon\downarrow0$ and deduce the non-strict comparison.  At a putative first contact $(t_1,x_1)$,
\[
 u=\overline U_\theta,
 \;
 \nabla u=\nabla\overline U_\theta=-\lambda\overline U_\theta\theta,
 \;
 \Delta u\le\Delta\overline U_\theta,
 \;
 \p_t(u-\overline U_\theta)\ge0.
\]
By \cref{F,A2},
\(f(t,x\cdot e-ct,u) \le f_-(t,u)\le r_-(t)-\beta(t)u.\)
Since prior to the first contact $u\le\overline U_\theta$, \cref{lem:expGreen} applies at $(t_1,x_1)$ and provides
\(-\eta\chi\nabla\!\cdot(u\nabla v) \le\gamma_{N,\lambda_*}(t)\overline U_\theta^2.\)
Therefore, at the contact point,
\begin{align*}
 u_t-a(t)\Delta u
 &\le-\eta\chi\nabla\!\cdot(u\nabla v)
   +r_-(t)u-\beta(t)u^2\\
 &\le r_-(t)\overline U_\theta
   -(\beta(t)-\gamma_{N,\lambda_*}(t))\overline U_\theta^2\\
 &\le r_-(t)\overline U_\theta
 <\overline U_{\theta,t}-a(t)\Delta\overline U_\theta,
\end{align*}
where the strict last inequality is \eqref{eq:strictExpIdentityN}.  Together with $\Delta u\le\Delta\overline U_\theta$, this implies $u_t<\overline U_{\theta,t}$, contradicting the first-contact inequality.  Thus \eqref{eq:envelopecontrolN} holds for all $t\ge0$.
\end{proof}

\begin{corollary}[Sharp radial upper cone]\label{cor:uppercone}
Under the assumptions of \cref{lem:superenv}, for every $s>c_*$,
\(\lim_{t\to\infty}\sup_{|x|\ge st}u(t,x)=0.\)
\end{corollary}

\begin{proof}
We let $\lambda$ and $\widehat c<s$ be chosen as in \cref{lem:superenv}.  For $x\ne0$ put $\theta=x/|x|$.  Since $P_\lambda$ and the constants in \cref{lem:superenv} are independent of $\theta$,
\(u(t,x)\le C\exp[-\lambda(|x|-\widehat c t)].\)
Hence, uniformly for $|x|\ge st$,
\(u(t,x)\le C e^{-\lambda(s-\widehat c)t}\to0.\)
\end{proof}

\begin{lemma}[Planar exponential upper bound for half-space data]\label{lem:halfupper}
Under \cref{A1,A2,F,D}, let $\tau=0$, and suppose that $u_0$ is bounded, nonnegative and uniformly continuous with
\(u_0(x)=0\;\text{for }x\cdot e\ge R_0.\)
For every $s>c_*$ there are $\widehat c\in(c_*,s)$ and $C>0$ such that
\begin{equation}\label{eq:halfupperexp}
 u(t,x)\le C P_{\lambda_*}(t)
 e^{-\lambda_*(x\cdot e-\widehat c t)}
 \;(t\ge0,\ x\in\mathbb R^N).
\end{equation}
Consequently
\(\lim_{t\to\infty}\sup_{x\cdot e\ge st}u(t,x)=0.\)
\end{lemma}

\begin{proof}
We put $\lambda=\lambda_*$ and choose $\widehat c$ with $c_*<\widehat c<s$.  Since $P_\lambda$ is positive and periodic, let
\(P_*:=\min_{t\in[0,T]}P_\lambda(t)>0.\)
We choose $C$ so large that
\(C P_\lambda(0)e^{-\lambda R_0}\ge2\|u_0\|_\infty+1.\)
Then the function
\(\overline U(t,x)=C P_\lambda(t) e^{-\lambda(x\cdot e-\widehat c t)}\)
strictly dominates $u_0$: on $x\cdot e<R_0$ its value at $t=0$ is at least $2\|u_0\|_\infty+1$, while on $x\cdot e\ge R_0$ the datum vanishes.  As in \eqref{eq:strictExpIdentityN},
\begin{equation}\label{eq:halfupperstrict}
 \overline U_t-a(t)\Delta\overline U-r_-(t)\overline U
 =\lambda(\widehat c-c_*)\overline U>0.
\end{equation}

We verify carefully that a loss of the inequality cannot occur only at spatial infinity.  Along the already constructed bounded solution, write the population equation in the linear form \eqref{eq:scalarrewrite}.  On every finite interval $0<t\le T_1$, the corresponding fundamental kernel satisfies an upper Aronson bound
\(\mathcal K(t,x;0,y) \le C_1t^{-N/2}\exp\!\left(-\frac{|x-y|^2}{C_2t}\right),\)
where $C_1,C_2$ depend on $T_1$ and the uniform coefficient bounds but not on $x,y$.  Because $u_0(y)=0$ for $y\cdot e\ge R_0$, integration first in the transverse variables provides
\begin{align}
 u(t,x)
 &\le C_1\|u_0\|_\infty t^{-N/2}
 \int_{y\cdot e<R_0}
 e^{-|x-y|^2/(C_2t)}\,dy\notag\\
 &\le C_3(T_1)\exp\!\left(
 -\frac{((x\cdot e-R_0)_+)^2}{C_4T_1}
 \right),
 \; 0<t\le T_1.\label{eq:longitudinalGaussian}
\end{align}
Thus $u/\overline U\to0$ as $x\cdot e\to+\infty$, uniformly on compact positive time intervals.  As $x\cdot e\to-\infty$, $\overline U\to+\infty$ whereas $u$ is bounded by \cref{prop:global}.  Therefore any sequence of almost-contact points on a finite time interval has bounded longitudinal coordinate.

We suppose now that the comparison fails, and let $t_1>0$ be the first failure time after replacing $C$ by $(1+\varepsilon)C$, $\varepsilon>0$.  If a first-contact sequence has bounded transverse component, compactness provides an actual contact point.  If its transverse component diverges, write
\(x_n=(x_n\cdot e)e+x_{n,\perp},\; x_{n,\perp}\cdot e=0,\)
and translate the solution by $x_{n,\perp}$.  The coefficients of the equation and the comparison function are invariant under translations in $e^\perp$.  Positive-time parabolic estimates and elliptic estimates therefore imply a locally convergent subsequence.  The initial domination has a uniform positive gap on the longitudinal strip $x\cdot e\le R_0$, namely at least $\|u_0\|_\infty+1$, and the initial datum is zero on its complementary half-space. Consequently transverse translations preserve strict domination locally at the initial time. The longitudinal Gaussian estimate prevents a first contact escaping in the $e$-direction, while locally uniform convergence at positive times supplies an actual contact for the limiting pair at time $t_1$.  Hence in either case it is enough to rule out a finite interior contact.

At such a contact, prior domination provides the hypotheses of \cref{lem:expGreen} with $\theta=e$.  Therefore
\(-\eta\chi\nabla\!\cdot(u\nabla v) \le\gamma_{N,\lambda_*}(t_1)\overline U^2.\)
We combine $f(t,\xi,u)\le r_-(t)-\beta(t)u$ with $\gamma_{N,\lambda_*}\le\beta$.  We obtain
\(u_t-a(t)\Delta u \le r_-(t_1)\overline U <\overline U_t-a(t_1)\Delta\overline U,\)
where the strict inequality is \eqref{eq:halfupperstrict}.  At a first contact one also has $\Delta u\le\Delta\overline U$ and $u_t\ge\overline U_t$, a contradiction.  Thus $(1+\varepsilon)\overline U$ dominates for every $\varepsilon>0$; letting $\varepsilon\downarrow0$ proves \eqref{eq:halfupperexp}.  Finally, if $x\cdot e\ge st$,
\(u(t,x)\le C\max_{[0,T]}P_\lambda e^{-\lambda(s-\widehat c)t}\longrightarrow0.\)
\end{proof}

The next estimate concerns growth in an arbitrary moving frame.  Let $B_L$ denote the ball of radius $L$ centered at the origin.

\begin{lemma}[Positive periodic growth on a bounded ball]\label{lem:Dirichletpositive}
Let $K\Subset\{s\in\mathbb R^N:|s|<c_*\}$ be compact.  There exist $L>2$, $\delta_{\mathrm{sp}}>0$ and $\sigma>0$ such that for every $s\in K$ and $0<\delta\le\delta_{\mathrm{sp}}$, the periodic Dirichlet problem
\begin{equation}\label{eq:pertFloquetN}
\begin{cases}
 \phi_t=a(t)\Delta\phi+s\cdot\nabla\phi+(r_-(t)-3\delta)\phi-\lambda_s\phi,
 &x\in B_L,\\
 \phi=0,&x\in\p B_L,\\
 \phi(t+T,x)=\phi(t,x)
\end{cases}
\end{equation}
has a positive normalized eigenfunction with $\lambda_s\ge4\sigma$.  Here and in the proof we suppress the dependence of $\lambda_s$ and $\phi_s$ on $\delta$; every estimate is uniform for $(s,\delta)\in K\times(0,\delta_{\mathrm{sp}}]$.  Moreover, for every prescribed $H_A>0$ there exists $\eta_0=\eta_0(H_A)>0$ such that, if $A(t,x)$ is a H\"older vector field satisfying
\[
 \lVert A\rVert_\infty\le\eta_0,
 \;
 \lVert A\rVert_{C^{\alpha_H/2,\alpha_H}([t_0,t_0+T]\times\overline{B_L})}\le H_A
 \;\text{for every }t_0\in\mathbb R,
\]
the Dirichlet evolution of
\begin{equation}\label{eq:linearrobustN}
 Z_t=a(t)\Delta Z+(s-A(t,x))\cdot\nabla Z+(r_-(t)-3\delta)Z
\end{equation}
satisfies, for every temporal phase $t_0$,
\begin{equation}\label{eq:oneperiodampN}
 Z(t_0+T,\cdot;t_0,\phi_s(t_0,\cdot))
 \ge e^{2\sigma T}\phi_s(t_0+T,\cdot).
\end{equation}
All constants are uniform for $(s,\delta)\in K\times(0,\delta_{\mathrm{sp}}]$.
\end{lemma}

\begin{proof}
For $s\in\mathbb R^N$ define
\(\Lambda(s):=\bar r-\frac{|s|^2}{4\bar a}.\)
Consider first the principal periodic Dirichlet exponent on $B_L$ for
\(z_t=a(t)\Delta z+s\cdot\nabla z+(r_-(t)-3\delta)z\)
converges, as $L\to\infty$, to $\Lambda(s)-3\delta$, uniformly for $s$ in compact sets.  Compactness alone does not determine the limiting exponent, since the whole-space equation admits several positive exponential Floquet modes.

We set
\[
 b_s:=\frac{s}{\bar a},
 \;
 h_s'(t):=s-a(t)b_s=s\left(1-\frac{a(t)}{\bar a}\right),
 \; h_s(0)=0.
\]
Since $\int_0^T h_s'(t)\,dt=0$, $h_s$ is $T$-periodic.  We let $\phi$ solve the periodic eigenvalue equation on the fixed ball $B_L$ and set
\(y=x+h_s(t), \; \psi(t,y)=\phi(t,y-h_s(t)).\)
The spatial domain becomes the periodic moving ball
\(D_{L,s}(t)=B_L+h_s(t),\)
and a direct differentiation provides
\(\psi_t=a(t)\Delta\psi+a(t)b_s\cdot\nabla\psi +(r_-(t)-3\delta-\lambda)\psi.\)
Now put
\(\psi(t,y)=e^{-b_s\cdot y/2}\zeta(t,y).\)
Because $b_s$ is constant,
\begin{align*}
 \Delta\psi+b_s\cdot\nabla\psi
 &=e^{-b_s\cdot y/2}
 \left(\Delta\zeta-\frac{|b_s|^2}{4}\zeta\right).
\end{align*}
Hence $\zeta$ satisfies, on $D_{L,s}(t)$,
\begin{equation}\label{eq:driftremovedmoving}
 \zeta_t=a(t)\Delta\zeta+
 \left(r_-(t)-3\delta-\frac{a(t)|s|^2}{4\bar a^2}-\lambda\right)\zeta,
 \; \zeta=0\text{ on }\partial D_{L,s}(t).
\end{equation}
Both the time-dependent translation and the positive gauge preserve the principal exponent.

We let
\(H_K:=\sup_{s\in K}\sup_{t\in[0,T]}|h_s(t)|<\infty.\)
For $L>H_K$,
\(B_{L-H_K}\subset D_{L,s}(t)\subset B_{L+H_K} \;(t\in\mathbb R,\ s\in K).\)
Domain monotonicity for positive periodic Dirichlet evolution therefore sandwiches the principal exponent of \eqref{eq:driftremovedmoving} between the corresponding exponents on the two fixed balls.  On a fixed ball $B_R$ the drift-free problem separates.  If $\psi_R$ is the first Dirichlet eigenfunction of $-\Delta$ with eigenvalue $\lambda_1(B_R)$, then
\(Z(t,x)=P_{R,s}(t)\psi_R(x)\)
is periodic provided
\begin{align*}
 \frac{P_{R,s}'}{P_{R,s}}
 &=r_-(t)-3\delta-\frac{a(t)|s|^2}{4\bar a^2}
 -a(t)\lambda_1(B_R)-\Lambda_{R,s},
\end{align*}
and periodicity determines
\begin{equation}\label{eq:fixedball-exponent}
 \Lambda_{R,s}
 =\bar r-3\delta-\frac{|s|^2}{4\bar a}
 -\bar a\lambda_1(B_R).
\end{equation}
Consequently
\begin{align*}
 \Lambda_{L-H_K,s}
 \le \lambda_{L,s}(\delta)
 \le \Lambda_{L+H_K,s},
\end{align*}
where the outer quantities are defined in \eqref{eq:fixedball-exponent}.  Since
$\lambda_1(B_R)=R^{-2}\lambda_1(B_1)$, both bounds converge uniformly for $s\in K$ to
\(\Lambda(s)-3\delta.\)

Because $K\Subset\{|s|<c_*\}$,
\(\Lambda_K:=\min_{s\in K}\Lambda(s)>0.\)
We choose $0<\sigma<\Lambda_K/16$, and then we choose $0<\delta_{\mathrm{sp}}<\Lambda_K/16$.  The uniform exhaustion just proved permits a choice of $L>2$ so large that
\[
 \lambda_s:=\lambda_{L,s}(\delta)\ge4\sigma
 \;(s\in K,\ 0<\delta\le\delta_{\mathrm{sp}}).
\]
We normalize the corresponding positive periodic eigenfunctions by
$\|\phi_s\|_{L^\infty([0,T]\times B_L)}=1$.

To control drift perturbations, we work in a norm that detects behavior near the boundary. Fix $H_A>0$, as in the statement, and let
$\mathcal U_{A,s}$ denote the Dirichlet evolution.  The unperturbed Dirichlet
evolution satisfies, uniformly in $s\in K$, phase $t_0$, and
$0<t-t_0\le T$,
\[
 \|\nabla\mathcal U_{0,s}(t,t_0)\varphi\|_\infty
 \le C\|\varphi\|_{C^1(\overline{B_L})},\qquad
 \|\nabla\mathcal U_{0,s}(t,t_0)F\|_\infty
 \le C(t-t_0)^{-1/2}\|F\|_\infty
\]
for compatible zero-Dirichlet $C^1$ data $\varphi$ and bounded $F$.
These are the standard boundary gradient estimates on the fixed smooth ball;
the constants depend on the fixed geometry, ellipticity, coefficient bounds
and $T$, but not on $A$ or $t_0$.  The normalized unperturbed eigenfunctions
have uniformly bounded $C^1$ norms.  For
$u_A(t)=\mathcal U_{A,s}(t,t_0)\phi_s(t_0)$, Duhamel's formula, first for
smooth drifts and then by approximation within the common H\"older class,
implies
\[
 \|\nabla u_A(t)\|_\infty
 \le C_0+C\|A\|_\infty
 \int_{t_0}^t(t-\zeta)^{-1/2}\|\nabla u_A(\zeta)\|_\infty\,d\zeta.
\]
The fractional Gronwall inequality implies
$\sup_{t_0\le t\le t_0+T}\|\nabla u_A(t)\|_\infty\le C_T$,
uniformly for $\|A\|_\infty\le1$ and the stated common H\"older bound.
The argument uses the integrable boundary-gradient bound on a fixed smooth ball and only the zero-order compatibility $\phi_s(t_0,\cdot)|_{\partial B_L}=0$ at the starting time. No higher-order compatibility at the initial lateral corner is used. In particular, the constant in the gradient estimate is allowed to depend on $L,T,a_*,a^*$ and on the fixed H\"older bound for the drift, but not on $s,t_0$ or the particular truncated drift.
Duhamel's formula now reads
\[
 u_A(t_0+T)-\mathcal U_{0,s}(t_0+T,t_0)\phi_s(t_0)
 =-\int_{t_0}^{t_0+T}\mathcal U_{0,s}(t_0+T,\zeta)
 [A(\zeta)\cdot\nabla u_A(\zeta)]\,d\zeta.
\]
We use the $C^1$ boundary gradient estimate for the unperturbed evolution,
\[
 \|\mathcal U_{0,s}(t_0+T,\zeta)F\|_{C^1(\overline{B_L})}
 \le C\big(1+(t_0+T-\zeta)^{-1/2}\big)\|F\|_\infty,
\]
and integrate the integrable singularity $(t_0+T-\zeta)^{-1/2}$.
Since $\mathcal U_{0,s}(t_0+T,t_0)\phi_s(t_0)
=e^{\lambda_sT}\phi_s(t_0+T)$, we obtain
\begin{equation}\label{eq:evolpertN}
 \left\|\mathcal U_{A,s}(t_0+T,t_0)\phi_s(t_0)
 -e^{\lambda_sT}\phi_s(t_0+T)\right\|_{C^1(\overline{B_L})}
 \le C\|A\|_\infty.
\end{equation}
The constant is uniform in $s,t_0,\delta,A$ under the displayed hypotheses.

The family of normalized principal eigenfunctions, with $(s,\delta)$ in the closure $K\times[0,\delta_{\mathrm{sp}}]$, is compact in $C^{1,2}$ on the closed periodic cylinder by the common boundary Schauder estimates. The normal derivative is strictly negative at every boundary point for every member of this compact family. Compactness and the Hopf lemma therefore provide $\phi_s(t,x)\ge c_H\operatorname{dist}(x,\partial B_L)$ with a constant $c_H>0$ uniform in $t,x,s$.
Since both terms in \eqref{eq:evolpertN} vanish on the lateral boundary, the $C^1$ estimate implies a pointwise bound by
$C'\|A\|_\infty\phi_s(t_0+T,x)$.  We choose $\eta_0>0$ so small that
\(C'\eta_0\le e^{4\sigma T}-e^{2\sigma T}.\)
We use $\lambda_s\ge4\sigma$ in \eqref{eq:evolpertN} and obtain
\(\CalU_{A,s}(t_0+T,t_0)\phi_s(t_0,\cdot) \ge e^{2\sigma T}\phi_s(t_0+T,\cdot),\)
which is \eqref{eq:oneperiodampN}.  All choices are uniform for $s\in K$.
\end{proof}

\begin{lemma}[Small positive entire state on a bounded ball under a truncated drift]\label{lem:entireaux}
Let $K\Subset\{|s|<c_*\}$ and choose $L,\delta_{\mathrm{sp}},\sigma$ from \cref{lem:Dirichletpositive}.  Fix $H_A>0$, let $\eta_0(H_A)$ be chosen as in that lemma, fix $p>1$, $0<\delta\le\delta_{\mathrm{sp}}$, and an arbitrary prescribed constant $C_*>0$.  For every $s\in K$, every globally H\"older vector field $A$ satisfying
\[
 \lVert A\rVert_\infty\le\eta_0(H_A),
 \;
 \sup_{t_0\in\mathbb R}
 \lVert A\rVert_{C^{\alpha_H/2,\alpha_H}([t_0,t_0+T]\times\overline{B_L})}\le H_A,
\]
and every $B>0$, the problem
\begin{equation}\label{eq:auxentireN}
\begin{cases}
 w_t=a(t)\Delta w+(s-A(t,x))\cdot\nabla w\\
 \hspace{14mm}+w\big(r_-(t)-3\delta-C_*w^{1/p}-Bw\big),
 &x\in B_L,\\
 w=0,&x\in\p B_L
\end{cases}
\end{equation}
has a positive bounded entire solution $w=w_{s,A,B}$.  It may be selected so that
\(0<w(t,x)\le Q_B, \; Q_B\downarrow0\;(B\to\infty),\)
and for every $\ell\in(0,L)$,
\begin{equation}\label{eq:auxinteriorN}
 \inf_{\substack{t\in\mathbb R,\ s\in K\\|x|\le L-\ell}}
 w_{s,A,B}(t,x)>0
\end{equation}
for each fixed $B>0$, uniformly over all admissible drifts $A$.
\end{lemma}

\begin{proof}
We let $\phi_s$ be the normalized positive eigenfunction in \cref{lem:Dirichletpositive}, with $\|\phi_s\|_\infty=1$.  By \eqref{eq:oneperiodampN}, the linear Dirichlet evolution with drift $A$ and zero-order coefficient $r_-(t)-3\delta$ satisfies
\begin{equation}\label{eq:auxlinearampN}
 \CalU_{A,s}(t_0+T,t_0)\phi_s(t_0,\cdot)
 \ge e^{2\sigma T}\phi_s(t_0+T,\cdot)
\end{equation}
uniformly in $s\in K$, in the admissible drift $A$, and in $t_0\in\mathbb R$.

We fix $B>0$ for the moment.  We choose $q_{\rm nl}=q_{\rm nl}(B)>0$ so small that
\begin{equation}\label{eq:nonlinlossN}
 C_*q_{\rm nl}^{1/p}+Bq_{\rm nl}\le\sigma.
\end{equation}
We let $w$ be the nonlinear solution of \eqref{eq:auxentireN} on one period, issued at time $t_0$ from $\zeta\phi_s(t_0,\cdot)$.  We put
\(R_-:=\big(\max_{t\in[0,T]}r_-(t)-3\delta\big)_+.\)
At a spatial maximum of $w$ the drift term vanishes and the nonlinear loss is nonpositive, hence the upper Dini derivative of $\|w(t,\cdot)\|_\infty$ is at most $R_-\|w(t,\cdot)\|_\infty$.  Consequently
\begin{equation}\label{eq:auxupperfirstperiodN}
 \|w(t,\cdot)\|_\infty
 \le e^{R_-(t-t_0)}\zeta,
 \; t_0\le t\le t_0+T,
\end{equation}
as long as the solution exists.  We choose
\begin{equation}\label{eq:auxseedsmallN}
 0<\zeta\le \frac12 q_{\rm nl}e^{-R_-T}.
\end{equation}
Then \eqref{eq:auxupperfirstperiodN} provides $0\le w\le q_{\rm nl}/2$ throughout the first period; in particular the small-loss inequality \eqref{eq:nonlinlossN} is valid everywhere on that orbit.

We let $z$ denote the corresponding linear solution.  We define \(\underline z(t,x)=e^{-\sigma(t-t_0)}z(t,x).\)
The maximum estimate used in \eqref{eq:auxupperfirstperiodN} also applies to $z$ and shows that $0\le z\le q_{\rm nl}/2$ on the first period.  Hence $0\le\underline z\le q_{\rm nl}/2$, so \eqref{eq:nonlinlossN} applies to $\underline z$.  We compute
\begin{align*}
 (\underline z)_t
 &=a(t)\Delta\underline z+(s-A)\cdot\nabla\underline z
   +(r_-(t)-3\delta-\sigma)\underline z\\
 &\le a(t)\Delta\underline z+(s-A)\cdot\nabla\underline z
   +\underline z\big(r_-(t)-3\delta
   -C_*\underline z^{1/p}-B\underline z\big).
\end{align*}
Thus $\underline z$ is a subsolution of the nonlinear equation, with the same initial datum as $w$.  Scalar comparison implies $w\ge\underline z$ on the first period.  We use \eqref{eq:auxlinearampN} and obtain
\begin{equation}\label{eq:nonlinampN}
 \CalS_{s,A,B}(t_0+T,t_0)
 [\zeta\phi_s(t_0,\cdot)]
 \ge e^{\sigma T}\zeta\phi_s(t_0+T,\cdot).
\end{equation}
Every constant in \eqref{eq:auxseedsmallN}--\eqref{eq:nonlinampN} is uniform in $s$, $A$, and $t_0$ once $B$ is fixed.

We choose $Q_B>0$ as the unique positive solution of
\(BQ_B+C_*Q_B^{1/p}=(r_-^*-3\delta)_++1, \; r_-^*:=\max_t r_-(t).\)
Then $Q_B\downarrow0$ as $B\to\infty$, and
\(r_-^*-3\delta-C_*Q_B^{1/p}-BQ_B\le-1.\)
Thus the constant $Q_B$ is a strict supersolution of \eqref{eq:auxentireN}.  Notice that $Q_B$ is not required to satisfy the small-loss condition \eqref{eq:nonlinlossN}; that condition is needed only for the single orbit issued from the small eigenfunction seed.

We decrease the seed $\zeta$ in \eqref{eq:auxseedsmallN}, if necessary, so that
\(\zeta<Q_B.\)
For $n\ge1$, solve \eqref{eq:auxentireN} from time $-nT$ with initial value $\zeta\phi_s(-nT)$.  The one-period estimate \eqref{eq:nonlinampN} provides
\[
 \CalS_{s,A,B}(-(n-1)T,-nT)[\zeta\phi_s(-nT)]
 \ge e^{\sigma T}\zeta\phi_s(-(n-1)T)
 \ge\zeta\phi_s(-(n-1)T).
\]
The scalar flow is order preserving.  Iterating this inequality shows that, at every fixed later integer phase, the pullback solution obtained from time $-(n+1)T$ dominates the one obtained from time $-nT$.  The strict supersolution $Q_B$ bounds every member of the pullback sequence from above.  Interior and boundary parabolic estimates therefore provide local compactness on the full cylinder, and monotone convergence implies a bounded entire limit $w_{s,A,B}$.  The limit dominates the positive seed at every integer phase, so the strong maximum principle provides $w_{s,A,B}>0$.

We prove the uniform interior lower bound.  The seed amplitude $\zeta$ can be chosen uniformly in $s$ and in all drifts with $\lVert A\rVert_\infty\le\eta_0$, because the constants in \eqref{eq:oneperiodampN} are uniform.  Hence
\(w(kT,\cdot)\ge\zeta\phi_s(kT,\cdot) \;(k\in\mathbb Z).\)
We fix $\ell\in(0,L)$.  The normalized eigenfunctions, now with $(s,\delta)$ ranging in the compact set $K\times[0,\delta_{\mathrm{sp}}]$, satisfy
\[
 m_{\ell}:=\inf_{s\in K,\ 0\le\delta\le\delta_{\mathrm{sp}},\ t\in[0,T],\ |x|\le L-\ell}
 \phi_{s,\delta}(t,x)>0.
\]
For fixed $B$, the zero-order coefficient in \eqref{eq:auxentireN} is bounded from below on $0\le w\le Q_B$ by a constant depending only on $B,C_*,M_0$ and the fixed coefficient bounds, but not on $s,A$ or the integer phase.  The parabolic Harnack inequality and the lower Gaussian bound on the fixed ball therefore propagate the lower bound $\zeta m_{\ell/2}$ at time $kT$ throughout the strip $[kT,(k+1)T]$ on $|x|\le L-\ell$.  Thus there is $c_{\ell,B}>0$, independent of $k,s,A$, such that
\(w(t,x)\ge c_{\ell,B}\zeta m_{\ell/2} \;(kT\le t\le(k+1)T,\ |x|\le L-\ell).\)
This proves \eqref{eq:auxinteriorN}.
\end{proof}

\begin{lemma}[Persistence in admissible vector frames]\label{lem:framepersist}
Suppose that \cref{A1,A2,F,D} hold. Let
\(K\Subset\Omega_c:=\{s\in\mathbb R^N:|s|<c_*,\ s\cdot e<c\}\)
be compact. There exists $L_K>1$, depending only on $K$ and the fixed coefficients, such that for every fixed nonnegative, nontrivial, compactly supported $u_0\in BUC(\mathbb R^N)$ there exists $\rho_{K,u_0}>0$ for which the corresponding solution satisfies
\begin{equation}\label{eq:framepositiveN}
 \liminf_{t\to\infty}
 \inf_{\substack{s\in K\\|x-st|\le L_K-1}}u(t,x)
 \ge\rho_{K,u_0}.
\end{equation}
The location scale $L_K$ is uniform on $K$; the lower-bound constant in this lemma is not asserted to be uniform over the class of all nontrivial initial data.
\end{lemma}

\begin{proof}
We choose $L,\delta,\sigma$ from the unperturbed part of \cref{lem:Dirichletpositive} for the compact set $K$.  Positive-time parabolic estimates for $u$, elliptic Schauder estimates for $v$, and the global bound \eqref{eq:globalbound} provide a constant $H_A>0$, depending on $L$ and the fixed data but not on $s\in K$, such that
\[
 \sup_{t_0\ge1}\lVert \eta\chi(t)\nabla V_s\rVert_{C^{\alpha_H/2,\alpha_H}([t_0,t_0+T]\times\overline{B_L})}\le H_A.
\]
We let $\eta_0=\eta_0(C_T H_A)>0$ be the drift threshold from \cref{lem:Dirichletpositive}, where $C_T$ is a fixed Lipschitz bound for the radial truncations used below.  Choose a localization radius $R>1$ so large that the Yukawa tail term in \cref{lem:localgreen}, after multiplication by the fixed coefficient bounds, is below $\min\{\delta,\eta_0\}/4$.  Then choose $p\ge p_0(R)$ as in \cref{lem:Harnack} and apply \cref{cor:smalldrift}.  Finally, choose $q>0$ so small that the local term is also below $\min\{\delta,\eta_0\}/4$.  Then
\begin{equation}\label{eq:driftsmallchoice}
 u(t,x)\le q
 \;\Longrightarrow\;
 \begin{cases}
 \eta\chi(t)|\nabla v(t,x)|\le\dfrac{\eta_0}{2},\\[1mm]
 \eta\chi(t)\kappa(t)^2v(t,x)\le\delta.
 \end{cases}
\end{equation}
for $t\ge1$.

For $s\in K$, introduce the moving functions
\(U_s(t,y)=u(t,y+st), \; V_s(t,y)=v(t,y+st).\)
Then
\begin{align}
 (U_s)_t
 &=a(t)\Delta U_s+
 (s-\eta\chi(t)\nabla V_s)\cdot\nabla U_s\notag\\
 &\;+U_s\left[f(t,y\cdot e+(s\cdot e-c)t,U_s)
 -\eta\chi(t)\kappa(t)^2V_s
 +\frac{\eta\chi(t)\mu(t)}{d(t)}U_s\right].
 \label{eq:movingUexact}
\end{align}
Because $K\Subset\Omega_c$, there is $\delta_c>0$ such that
\(s\cdot e-c\le-\delta_c \;(s\in K).\)
Hence, uniformly for $s\in K$ and $y\in B_L$,
\(y\cdot e+(s\cdot e-c)t\to-\infty.\)
By \cref{F}, after increasing a time $t_0$ we have
\begin{equation}\label{eq:favapproxN}
 r(t,y\cdot e+(s\cdot e-c)t)\ge r_-(t)-\delta
 \;(t\ge t_0,\ s\in K,\ y\in B_L).
\end{equation}
We let $b_1=b_1(M_0)$.  We combine \eqref{eq:kppbounds}, \eqref{eq:driftsmallchoice}, and \eqref{eq:favapproxN}.  Since the last quadratic term in \eqref{eq:movingUexact} is nonnegative, we obtain, at every point where $U_s\le q$,
\begin{equation}\label{eq:UlowerN}
 (U_s)_t\ge a(t)\Delta U_s+
 (s-\eta\chi\nabla V_s)\cdot\nabla U_s
 +U_s\big(r_-(t)-2\delta-CU_s^{1/p}-b_1U_s\big).
\end{equation}

We choose a smooth radial truncation $\mathfrak T_{\eta_0}:\mathbb R^N\to\mathbb R^N$ satisfying
\[
 \mathfrak T_{\eta_0}(z)=z\;(|z|\le\eta_0/2),
 \;
 |\mathfrak T_{\eta_0}(z)|\le\eta_0\;(z\in\mathbb R^N),
\]
and define
\begin{equation}\label{eq:truncAN}
 A_s(t,y)=\mathfrak T_{\eta_0}
 \big(\eta\chi(t)\nabla V_s(t,y)\big).
\end{equation}
We choose the radial truncation with Lipschitz norm at most $C_T$, independently of $\eta_0$.  The preceding estimate implies
\[
 \sup_{t_0\ge1}\lVert A_s\rVert_{C^{\alpha_H/2,\alpha_H}([t_0,t_0+T]\times\overline{B_L})}\le C_TH_A
 \;(s\in K),
\]
and $|A_s|\le\eta_0$.  Moreover $A_s=\eta\chi\nabla V_s$ whenever the latter vector has norm at most $\eta_0/2$.  Hence \eqref{eq:driftsmallchoice} makes the truncation inactive at every point where $U_s\le q$.

The auxiliary entire-state lemma is formulated for drifts defined on the whole time axis.  We therefore extend each \(A_s\) from \([1,\infty)\times\overline{B_L}\) to \(\mathbb R\times\overline{B_L}\) by setting it equal to \(A_s(1,\cdot)\) for \(t\le1/2\) and interpolating smoothly on \(1/2<t<1\).  The extension preserves the bound \(|A_s|\le\eta_0\) after composing once more with the radial truncation, and its H\"older norm is uniform in \(s\).  Since all comparisons below start at a time \(t_0\ge1\), the arbitrary part of this extension never enters the equation for the original solution.

We let $C$ be the constant in \eqref{eq:UlowerN} and fix once and for all a number $C_*\ge C$.  We apply \cref{lem:entireaux} with this prescribed $C_*$ and with the common drift H\"older bound $C_TH_A$.  We then choose
\(B\ge b_1+1\)
so large that the selected entire solution $w_{s,A_s,B}$ satisfies $0<w<q/2$.  Since $u(t_0,\cdot)>0$ everywhere by the strong maximum principle and
\(\bigcup_{s\in K}(st_0+\overline{B_L})\)
is compact, we increase $B$ once more, if necessary, so that
\(w_{s,A_s,B}(t_0,y)<U_s(t_0,y) \;(s\in K,\ y\in B_L).\)
The value of $B$ is chosen after the initial datum has been fixed, since the preceding strict inequality involves the positive minimum of that solution at time $t_0$. Once this choice is made, the auxiliary-state lower bound \eqref{eq:auxinteriorN} is uniform for $s\in K$ and for all truncated drifts with the common bounds. It provides the constant $\rho_{K,u_0}$ in \eqref{eq:framepositiveN}; the argument makes no assertion about its independence of $u_0$. All subsequent comparisons use this same choice of $C_*$, $B$, and $w_{s,A_s,B}$.  We suppose that $U_s-w_{s,A_s,B}$ has a first zero at $(t_1,y_1)$.  It is an interior contact because $w=0$ on $\p B_L$ while $U_s>0$.  At the contact,
\(U_s=w<q, \; \nabla U_s=\nabla w, \; \Delta U_s\ge\Delta w.\)
By \eqref{eq:driftsmallchoice}, the clipping is inactive at the contact, so $A_s=\eta\chi\nabla V_s$.  We put $Z=U_s-w$.  At $(t_1,y_1)$, the equalities above and the inactive clipping provide
\(\nabla Z=0,\; \Delta Z\ge0,\; A_s=\eta\chi\nabla V_s.\)
We subtract the equation \eqref{eq:auxentireN} satisfied by this same $w$ from \eqref{eq:UlowerN} at the contact, where $U_s=w$.  We obtain
\[
 \p_t Z\ge a(t_1)\Delta Z+(s-A_s)\cdot\nabla Z
 +w\big[\delta+(C_*-C)w^{1/p}+(B-b_1)w\big]>0.
\]
On the other hand, $t_1$ is the first time at which the previously positive function $Z$ reaches zero at an interior point, so $\p_tZ(t_1,y_1)\le0$.  The contradiction proves $U_s\ge w$ for all later times.  The uniform interior lower bound \eqref{eq:auxinteriorN} proves \eqref{eq:framepositiveN}.
\end{proof}

\begin{lemma}[Persistence behind an occupied planar front]\label{lem:halffilled}
Under \cref{A1,A2,F,D}, $\tau=0$, and
\begin{equation}\label{eq:leftfilledonly}
 \liminf_{R\to\infty}\inf_{x\cdot e\le-R}u_0(x)>0.
\end{equation}
If $s<\min\{c,c_*\}$, then there is $\rho_s>0$ such that
\begin{equation}\label{eq:halfpositiveN}
 \liminf_{t\to\infty}\inf_{x\cdot e\le st}u(t,x)\ge\rho_s.
\end{equation}
The lower bound is uniform for $s$ in a compact subset of $(-\infty,\min\{c,c_*\})$.
\end{lemma}

\begin{proof}
We fix such an $s$.  We choose a scalar speed $\bar s\in(-c_*,c_*)$ and $\vartheta>0$ so that
\(s+4\vartheta<\bar s<c_*-4\vartheta.\)
This is possible because $s<c_*$.  We write every point as
\(x=(x\cdot e)e+x_\perp, \; x_\perp\cdot e=0,\)
and introduce the longitudinal moving variable
\(y=x\cdot e-st.\)
The actual population in this frame is still allowed to depend on $x_\perp$:
\(U(t,y,x_\perp)=u(t,(y+st)e+x_\perp), \; V(t,y,x_\perp)=v(t,(y+st)e+x_\perp).\)
Since $s<c$, for each fixed right boundary $R_0$,
\((s-c)t+y\to-\infty \;\text{uniformly for }y\le R_0.\)
We choose $t_0$ so large that
\(r(t,(s-c)t+y)\ge r_-(t)-\delta \;(t\ge t_0,\ y\le R_0),\)
where $\delta>0$ is sufficiently small for the one-dimensional Dirichlet operator with speed $\bar s$ to have positive periodic exponent.  We fix $\eta_0\in(0,\vartheta)$.  We choose $R>1$ so large that the Yukawa-tail contribution in \cref{lem:localgreen} is smaller than $\min\{\eta_0,\delta\}/4$, and choose $p\ge p_0(R)$ as in \cref{lem:Harnack}.  Then \cref{cor:smalldrift} allows us to choose $q>0$ so small that whenever $U\le q$,
\begin{equation}\label{eq:halfsmallN}
 \eta\chi|\nabla V|\le\eta_0<\vartheta,
 \;
 \eta\chi\kappa^2V\le\delta.
\end{equation}
Then, at such points,
\begin{equation}\label{eq:halflowerUN}
 U_t\ge a(t)\Delta U+(se-\eta\chi\nabla V)\cdot\nabla U
 +U\big(r_-(t)-2\delta-CU^{1/p}-b_1U\big).
\end{equation}

We construct a scalar function $w=w(t,y)$ on $(-\infty,R_0)$ solving
\begin{equation}\label{eq:halfauxN}
\begin{cases}
 w_t=a(t)w_{yy}+\bar s w_y
 +w\big(r_-(t)-3\delta-Cw^{1/p}-Bw\big),&y<R_0,\\
 w(t,R_0)=0,
\end{cases}
\end{equation}
where $B>b_1+1$ is large.  We choose $0<w(t_0,y)<q$, $w_y(t_0,y)<0$, and
\(\lim_{y\to-\infty}w(t_0,y)=q_0>0.\)
We take, for example,
\(w(t_0,y)=q_0\big(1-e^{\beta(y-R_0)}\big), \; y<R_0,\)
with $0<q_0<q$ and $\beta>0$.  For $B$ sufficiently large the constant $q$ is a supersolution, so the scalar half-line problem has a global solution satisfying $0<w<q$.  Its monotonicity is most cleanly obtained by sliding.  For every $h>0$, the function $w_h(t,y)=w(t,y+h)$ is defined for $y<R_0-h$, satisfies the same equation there, and obeys
\(w_h(t,R_0-h)=0<w(t,R_0-h), \; w_h(t_0,y)<w(t_0,y).\)
We apply the scalar maximum principle on the unbounded interval through an exhaustion argument.  The only point that requires care is the boundary at $y=-\infty$.  Fix $T_1>t_0$ and let $\zeta$ be the solution on $[t_0,T_1]$ of
\[
 \zeta'=\zeta\big(r_-(t)-3\delta-C\zeta^{1/p}-B\zeta\big),
 \qquad \zeta(t_0)=q_0.
\]
We claim that
\begin{equation}\label{eq:halftailODE}
 \lim_{y\to-\infty}w(t,y)=\zeta(t)
 \quad\text{uniformly for }t\in[t_0,T_1].
\end{equation}
To prove the claim, take an arbitrary sequence $y_n\to-\infty$ and set $w_n(t,z)=w(t,z+y_n)$.  On every compact cylinder $[t_0,T_1]\times[-R,R]$, the translated right boundary $R_0-y_n$ tends to $+\infty$, while $0<w_n<q$.  Interior parabolic Schauder estimates therefore yield, after extraction, convergence in $C^{1,2}$ on compact subsets of $(t_0,T_1]\times\mathbb R$ and uniformly up to $t=t_0$ to a bounded entire-space solution $W$ of the same scalar equation.  The explicit initial profile satisfies
\[
 w_n(t_0,z)=q_0\big(1-e^{\beta(z+y_n-R_0)}\big)\longrightarrow q_0
\]
uniformly for $z$ in compact sets.  Uniqueness for the scalar Cauchy problem on $\mathbb R$ with constant initial datum implies $W(t,z)\equiv\zeta(t)$.  Because every sequence $y_n\to-\infty$ has the same subsequential limit, the convergence is locally uniform in $t$; if the uniform convergence in \eqref{eq:halftailODE} failed, a sequence of violating pairs $(t_n,y_n)$ with $t_n\to t_*\in[t_0,T_1]$ would contradict the preceding translated compactness argument.  This proves \eqref{eq:halftailODE}.

Consequently $w_h(t,y)-w(t,y)\to0$ as $y\to-\infty$, uniformly for $t\in[t_0,T_1]$.  For $\varepsilon>0$, choose $J$ so large that $w_h-w\le\varepsilon$ on $[t_0,T_1]\times\{-J\}$.  On $(-J,R_0-h)$ the difference $d_h=w_h-w$ satisfies
\[
 (d_h)_t=a(t)(d_h)_{yy}+\bar s(d_h)_y+c_h(t,y)d_h,
\]
where, by the mean-value theorem and $0<w_h,w<q$, the coefficient $c_h$ is bounded above by a constant $C_0$ independent of $J$.  Since $d_h(t_0,y)<0$ and $d_h(t,R_0-h)=-w(t,R_0-h)<0$, the function $d_h-\varepsilon e^{C_0(t-t_0)}$ is nonpositive on the parabolic boundary of $[t_0,T_1]\times[-J,R_0-h]$ and satisfies the corresponding homogeneous maximum-principle inequality.  Hence
\[
 d_h(t,y)\le\varepsilon e^{C_0(t-t_0)}.
\]
First let $J\to\infty$ with the chosen boundary estimate and then let $\varepsilon\downarrow0$.  We obtain $w_h\le w$ on $[t_0,T_1]\times(-\infty,R_0-h)$.  Since $T_1$ is arbitrary and the right boundary inequality is strict, the strong maximum principle provides
\begin{equation}\label{eq:halfmonotoneN}
 w(t,y+h)<w(t,y)
 \;(t>t_0,\ y<R_0-h),
\end{equation}
so in particular $w_y\le0$ for $t>t_0$.

Because $\bar s\in(-c_*,c_*)$, the one-dimensional periodic Dirichlet operator
\(z_t=a(t)z_{yy}+\bar s z_y+(r_-(t)-3\delta)z\)
has positive principal exponent on $J_L=(R_0-2L,R_0)$ once $L$ is sufficiently large.  We keep the coefficient $B$ in the half-line equation \eqref{eq:halfauxN} fixed from now on.  We apply \cref{lem:entireaux} in one dimension with zero drift perturbation, prescribed coefficient $C_*=C$, and a larger damping coefficient $\widehat B>B$.  This provides a positive bounded entire Dirichlet solution $z_{L,\widehat B}$ on $J_L$, zero at both endpoints, and
\(m_L:=\inf_{t\in\mathbb R} z_{L,\widehat B}(t,R_0-L)>0.\)
Moreover its supersolution bound $Q_{\widehat B}$ tends to zero as $\widehat B\to\infty$.  Boundary Schauder estimates on the fixed interval provide $\sup_t\|z_{L,\widehat B}(t,\cdot)\|_{C^1(\overline{J_L})}\le C_LQ_{\widehat B}$.
The fixed half-line solution $w$ satisfies $w(t_0,y)>0$ in $J_L$ and, from the explicit initial profile,
\[
 w_y(t_0,R_0)=-\beta q_0<0.
\]
Therefore there are $\ell_0>0$ and $c_0>0$ such that
\(w(t_0,y)\ge c_0(R_0-y) \;(R_0-\ell_0\le y<R_0).\)
On the compact interval $[R_0-2L,R_0-\ell_0]$, $w(t_0,\cdot)$ has a positive minimum.  We take $\widehat B$ sufficiently large and use both the $C^1$ estimate near $R_0$ and the $C^0$ estimate on the compact remainder.  We then obtain
\(z_{L,\widehat B}(t_0,y)<w(t_0,y) \;(R_0-2L<y<R_0).\)
Because $\widehat B>B$, the function $z_{L,\widehat B}$ is a subsolution of the half-line equation \eqref{eq:halfauxN} on $J_L$.  At the two endpoints it is zero and hence does not exceed $w$.  Scalar comparison yields
\(z_{L,\widehat B}(t,y)\le w(t,y) \;(t\ge t_0,\ y\in J_L).\)
We evaluate at $R_0-L$ and use \eqref{eq:halfmonotoneN}; hence
\(w(t,y)\ge w(t,R_0-L)\ge m_L \;(t\ge t_0,\ y\le R_0-L).\)
For the remaining compact strip $[R_0-L,R_0-\ell]$, positivity, parabolic Harnack and the lower bound at $R_0-L$ provide a positive lower bound uniform for $t\ge t_0+1$.  Consequently, for every $\ell>0$,
\begin{equation}\label{eq:halfwposN}
 \liminf_{t\to\infty}\inf_{y\le R_0-\ell}w(t,y)>0.
\end{equation}

It remains to compare the planar function $w(t,x\cdot e-st)$ with the $N$-dimensional solution.  The lower bound needed at time $t_0$ is uniform also in the transverse variables.  Indeed, by \eqref{eq:leftfilledonly} there are $R_1>0$ and $m_1>0$ such that $u_0(z)\ge m_1$ whenever $z\cdot e\le-R_1$.  For a point $x$ with $x\cdot e\le R_0+st_0$, we choose
\(z_x=x-\ell_x e, \; \ell_x=\max\{0,x\cdot e+R_1+1\}.\)
Then $z_x\cdot e\le-R_1-1$.  Moreover
\(0\le\ell_x\le L_0:=\max\{0,R_0+st_0+R_1+1\},\)
so $B_{1/2}(z_x)$ lies in the uniformly occupied half-space and remains at distance at most $L_0+1/2$ from $x$, independently of the transverse component.  We apply the lower Aronson estimate over the fixed time interval $[0,t_0]$, with constants $A_*(t_0),\alpha^*(t_0)>0$, and obtain
\begin{align*}
 u(t_0,x)
 &\ge m_1\int_{B_{1/2}(z_x)}\mathcal K(t_0,x;0,z)\,dz\\
 &\ge m_1 A_*(t_0)t_0^{-N/2}|B_{1/2}|
 \exp\!\left(-\frac{\alpha^*(t_0)(L_0+1/2)^2}{t_0}\right)=:m_2>0.
\end{align*}
Thus
\(\inf_{x\cdot e\le R_0+st_0}u(t_0,x)\ge m_2>0.\)
We take $q_0$ and the initial amplitude of $w$ smaller so that
\(w(t_0,x\cdot e-st_0)<u(t_0,x) \;(x\cdot e-st_0\le R_0).\)
We justify the first-contact argument on the unbounded moving half-space.  We fix a finite time $T_1>t_0$ and write $h=t-t_0$.  The case $h=0$ follows from the uniform bound $u(t_0,x)\ge m_2$ on $x\cdot e\le R_0+st_0$ proved above.  We let $0<h\le T_1-t_0$ and let $x\cdot e\le R_0+st$.  We set
\(z_{x,h}=x-sh\,e-\sqrt h\,e.\)
For every $z\in B_{\sqrt h/2}(z_{x,h})$ we have
\(z\cdot e\le R_0+st_0-\frac12\sqrt h,\)
so the preceding uniform lower bound implies $u(t_0,z)\ge m_2$ throughout this ball.  Moreover,
\(|x-z|\le |s|h+\frac32\sqrt h.\)
The lower Aronson estimate from time $t_0$ to $t=t_0+h$ therefore provides constants $A_*,\alpha^*>0$, depending only on $T_1-t_0$ and the fixed coefficient bounds, such that
\begin{align*}
 u(t,x)
 &\ge m_2A_*h^{-N/2}\int_{B_{\sqrt h/2}(z_{x,h})}
 \exp\!\left(-\frac{\alpha^*|x-z|^2}{h}\right)dz\\
 &\ge m_2A_*|B_{1/2}|
 \exp\!\left[-\alpha^*\left(|s|\sqrt{T_1-t_0}+\frac32\right)^2\right].
\end{align*}
Together with the bound at $t=t_0$, this proves a positive lower bound $m(T_1)$ for $u$ throughout $t_0\le t\le T_1$ and $x\cdot e\le R_0+st$.
In particular an almost-contact sequence cannot approach the moving boundary $y=R_0$, because $w(t,y)\to0$ there uniformly for $t\in[t_0,T_1]$.  Set
\[
 Z(t,x):=u(t,x)-w(t,x\cdot e-st).
\]
Suppose, for contradiction, that the comparison fails before $T_1$, and let $t_1\in(t_0,T_1]$ be the first failure time.  Then $Z\ge0$ for $t_0\le t\le t_1$ in the corresponding moving half-space and
\(\inf_{x\cdot e-st_1<R_0}Z(t_1,x)=0.\)
Choose an almost-contact sequence $x_n$ at time $t_1$, so that, with
\(y_n=x_n\cdot e-st_1<R_0,\)
we have $Z(t_1,x_n)\to0$.

If $(y_n)$ stays bounded but the transverse component of $x_n$ is unbounded, translate the solution by that transverse component.  The coefficients and the planar comparison function are invariant under translations in $e^\perp$; parabolic and elliptic compactness on $[t_0,t_1]$ then produce an actual contact at a finite point in the limit.  It remains to exclude $y_n\to-\infty$.  Since $w_y\le0$ and $0<w<q$, the limit
\(w_-(t):=\lim_{y\to-\infty}w(t,y)\)
exists uniformly for $t\in[t_0,t_1]$ by \eqref{eq:halftailODE}, and it satisfies
\[
 w_-'=w_-\big(r_-(t)-3\delta-Cw_-^{1/p}-Bw_-\big).
\]
For this almost-contact sequence it is preferable to remain in the moving frame introduced above.  Put
\(\xi_n:=x_n-st_1e\), so that \(\xi_n\cdot e=y_n\), and define the spatial translates
\[
 U_n(t,z):=U(t,z+\xi_n),\qquad V_n(t,z):=V(t,z+\xi_n).
\]
For $t\in[t_0,t_1]$ and $z$ in a fixed compact set, the argument of the demographic coefficient in the $s$-frame is
\[
 z\cdot e+y_n+(s-c)t\longrightarrow-\infty
\]
uniformly.  Hence the translated demographic coefficient converges locally uniformly on the fixed finite cylinder to the favorable coefficient.  At the same time the translated planar comparator satisfies
\[
 w\bigl(t,(z+\xi_n)\cdot e\bigr)
 =w\bigl(t,z\cdot e+y_n\bigr)
 \longrightarrow w_-(t)
\]
locally uniformly in $(t,z)$.  The global bounds together with parabolic regularity for $U_n$ and elliptic regularity for $V_n$ yield, after passage to a subsequence, a limiting population--signal pair $(U_\infty,V_\infty)$ on $[t_0,t_1]$ satisfying the favorable-medium limiting equations.  Since $Z\ge0$ before the first failure, the translated comparison inequality passes to the limit, and we obtain
\[
 U_\infty(t,z)\ge w_-(t)\quad(t\le t_1),
 \qquad U_\infty(t_1,0)=w_-(t_1).
\]
Thus $Z_\infty=U_\infty-w_-$ has an actual first zero at $(t_1,0)$.  There
\(\nabla Z_\infty=0\), \(\Delta Z_\infty\ge0\), and \(\partial_tZ_\infty\le0\).  Since the comparator is spatially constant, its drift contribution vanishes.  In a neighborhood of the contact point one has $U_\infty<q$, and hence, for all sufficiently large $n$, $U_n<q$ there.  We may therefore pass \eqref{eq:halflowerUN} to the limit locally and subtract the displayed ODE for $w_-$.  At the contact point this yields
\[
 \partial_tZ_\infty
 \ge a(t_1)\Delta Z_\infty
 +w_-(t_1)\bigl[\delta+(B-b_1)w_-(t_1)\bigr]>0,
\]
a contradiction.  Therefore escape of an almost-contact sequence to $y=-\infty$ is impossible.  Consequently, if comparison fails at a finite time, it has an actual interior first contact at a finite point.

We suppose therefore that first contact occurs at $(t_1,x_1)$ with $y_1=x_1\cdot e-st_1<R_0$.  At that point,
\(\nabla U=w_y e, \; \Delta U\ge w_{yy}, \; U=w<q.\)
The transverse derivatives of the planar comparison function vanish.  By \eqref{eq:halfsmallN},
\(|\eta\chi\p_eV|\le\eta_0.\)
Since $w_y\le0$ and $\bar s>s+\eta_0$,
\((s-\eta\chi\p_eV)w_y-\bar s w_y\ge0.\)
At the first contact, \(Z=0\), \(\nabla Z=0\), \(\Delta Z\ge0\), and \(U=w\).  We subtract \eqref{eq:halfauxN} from \eqref{eq:halflowerUN} and use the displayed drift inequality.  We obtain
\[
 \p_tZ(t_1,x_1)\ge a(t_1)\Delta Z
 +\big[(s-\eta\chi\p_eV)-\bar s\big]w_y
 +w\big[\delta+(B-b_1)w\big]>0.
\]
The second term is nonnegative because
\(s-\eta\chi\p_eV-\bar s \le s+\eta_0-\bar s<0\)
and \(w_y\le0\), whereas the last term is strictly positive.  But an interior first zero of a function that was positive for earlier times satisfies \(\p_tZ(t_1,x_1)\le0\).  This contradiction proves
\(u(t,x)\ge w(t,x\cdot e-st)\)
throughout the moving half-space for all later times, and \eqref{eq:halfwposN} provides \eqref{eq:halfpositiveN}.  Uniformity on compact speed sets follows from choosing $\bar s,\vartheta,L,\delta$ uniformly.
\end{proof}

The favorable-medium rigidity theorem uses a positive scalar periodic state. We establish its uniqueness and two-sided bounds before turning to that theorem; the same state also enters the forced-wave analysis.

\begin{lemma}[Positive favorable periodic state]\label{lem:pminus}
Under \cref{A2,F}, equation \eqref{eq:pminus-early} has a unique positive $T$-periodic solution $p_-(t)$.  Moreover there are constants $0<p_*\le p^*<\infty$, depending only on $f_-$, such that
\(p_*\le p_-(t)\le p^*\;(t\in\mathbb R).\)
In the elliptic signal regime the corresponding spatially homogeneous signal is \eqref{eq:qminus-ell}.
\end{lemma}

\begin{proof}
We let $\Pi(p_0)$ denote the value at time $T$ of the solution of \eqref{eq:pminus-early} satisfying $p(0)=p_0>0$.  We verify the expansion property near zero without suppressing the dependence of the orbit on $p_0$.  Since $r_-$ is bounded, choose $R_0>0$ so that $|r_-(t)|\le R_0$ on $[0,T]$.  By continuity of $f_-$ at $u=0$, there exists $\delta>0$ such that
\(|f_-(t,u)-r_-(t)|\le\bar r/4\) for $t\in[0,T]$ and $0\le u\le\delta$.  On this strip,
\(f_-(t,u)\le R_0+\bar r/4=:R_1.\)
Hence every solution starting from
\(0<p_0\le\varepsilon_0:=\delta e^{-R_1T}\)
satisfies, as long as it remains below $\delta$,
\(p(t;p_0)\le p_0e^{R_1t}\le\delta\) for $0\le t\le T$.  If a first time $t_*\le T$ existed with $p(t_*;p_0)=\delta$, then the same differential inequality on $[0,t_*]$ would imply
\(\delta=p(t_*;p_0)\le p_0e^{R_1t_*}\le p_0e^{R_1T}\le\delta.\)
For $p_0<\varepsilon_0$ the last inequality is strict, a contradiction; the endpoint $p_0=\varepsilon_0$ follows by continuity from smaller initial values.  Thus the orbit remains in $[0,\delta]$ on $[0,T]$.  Consequently
\[
 \log\frac{\Pi(p_0)}{p_0}
 =\int_0^T f_-(t,p(t;p_0))\,dt
 \ge \int_0^T\left(r_-(t)-\frac{\bar r}{4}\right)dt
 =\frac{3T\bar r}{4}>0,
\]
and therefore $\Pi(p_0)>p_0$ for $0<p_0\le\varepsilon_0$.

For a large-data barrier, \eqref{eq:fuass} implies
\(f_-(t,p)\le r_-(t)-\beta(t)p\le r_-(t)-\beta_*p.\)
Choose
\(P_0>\max_{t\in[0,T]}r_-(t)_+/\beta_*+1.\)
At any time when a solution reaches $P_0$ its derivative is strictly negative, so $[0,P_0]$ is forward invariant.  For the solution with $p(0)=P_0$, the derivative is already strictly negative at $t=0$; uniqueness prevents a later upward crossing of $P_0$, hence $p(T;P_0)<P_0$ and $\Pi(P_0)<P_0$.  Continuity of the Poincar\'e map $\Pi$ now yields a fixed point in $(\varepsilon_0,P_0)$, whose orbit is a positive $T$-periodic solution.

For uniqueness, let $p_1,p_2$ be two positive periodic solutions.  Uniqueness for the scalar initial-value problem implies that they either coincide or are strictly ordered for every time.  If $p_1(t)>p_2(t)$ for all $t$, then
\[
 0=\int_0^T\frac{d}{dt}\log\frac{p_1(t)}{p_2(t)}\,dt
 =\int_0^T\big(f_-(t,p_1(t))-f_-(t,p_2(t))\big)\,dt
 \le-\beta_*\int_0^T(p_1(t)-p_2(t))\,dt<0,
\]
a contradiction.  The case $p_2>p_1$ is identical after exchanging the two solutions.  Thus the periodic solution is unique.  Its minimum and maximum on one period define constants $0<p_*\le p^*<\infty$ and, by periodicity, establish the asserted bounds on all of $\mathbb R$.

Finally, if $u(t,x)=p_-(t)$ is spatially homogeneous and $\tau=0$, the signal equation reduces to
\(-\nu(t)q(t)+\mu(t)p_-(t)=0\),
which is precisely \eqref{eq:qminus-ell}.
\end{proof}

\begin{lemma}[Rigidity in the favorable periodic medium]\label{lem:Liouville}
Suppose that \eqref{eq:strongdamp} holds.  Let $(U,V)$ be a bounded entire solution of
\begin{equation}\label{eq:favKSN}
\begin{cases}
 U_t=a(t)\Delta U-\eta\chi(t)\nabla\!\cdot(U\nabla V)+U f_-(t,U),\\
 d(t)\Delta V-\nu(t)V+\mu(t)U=0
\end{cases}
\end{equation}
on $\mathbb R\times\mathbb R^N$, and assume
\(\inf_{(t,x)\in\mathbb R\times\mathbb R^N}U(t,x)>0.\)
Then
\(U(t,x)=p_-(t), \; V(t,x)=q_-(t).\)
\end{lemma}

\begin{proof}
We set
\[
 M(t)=\sup_xU(t,x),
 \;
 m(t)=\inf_xU(t,x),
 \;
 \gamma(t)=\frac{\eta\chi(t)\mu(t)}{d(t)}.
\]
The positive resolvent kernel has total mass $1/\kappa(t)^2=d(t)/\nu(t)$, hence
\begin{equation}\label{eq:VbetweenN}
 \frac{\mu(t)}{\nu(t)}m(t)
 \le V(t,x)\le
 \frac{\mu(t)}{\nu(t)}M(t).
\end{equation}
We justify the extremal differential inequalities by quadratic penalization.  For a fixed time $t$, let $x_\varepsilon$ maximize $U(t,x)-\varepsilon|x|^2$.  Boundedness of $U$ implies $\varepsilon|x_\varepsilon|^2\le2\|U(t,\cdot)\|_\infty$, while
\[
 \nabla U(t,x_\varepsilon)=2\varepsilon x_\varepsilon,
 \;
 \Delta U(t,x_\varepsilon)\le2N\varepsilon.
\]
Hence $|\nabla U(t,x_\varepsilon)|\to0$ and the diffusion error tends to zero as $\varepsilon\downarrow0$.  The drift error tends to zero as well because the elliptic resolvent estimate provides a uniform bound for $\nabla V$.  We apply the same argument to $-U$, or equivalently we minimize $U(t,x)+\varepsilon|x|^2$, and obtain the corresponding lower Dini inequality for $m(t)$.  We use
\(\Delta V=\kappa^2V-\frac\mu dU\)
and \eqref{eq:VbetweenN}, and then we let the penalization vanish.  We obtain
\begin{align}
 D^+\log M(t)
 &\le f_-(t,M(t))+\gamma(t)(M(t)-m(t)),
 \label{eq:MineqN}\\
 D_-\log m(t)
 &\ge f_-(t,m(t))-\gamma(t)(M(t)-m(t)).
 \label{eq:mineqN}
\end{align}
We let $Q(t)=\log(M(t)/m(t))$.  From $f_{-,u}\le-\beta(t)$,
\begin{equation}\label{eq:QcontractN}
 D^+Q(t)\le f_-(t,M)-f_-(t,m)+2\gamma(t)(M-m)
 \le-(\beta(t)-2\gamma_0(t))(M-m).
\end{equation}
If $m_*:=\inf_{t,x}U(t,x)>0$, then
\(M-m=m(e^Q-1)\ge m_*(e^Q-1).\)
We set $\delta_L=\delta_2>0$.  Since $e^Q-1\ge Q$ for $Q\ge0$, \eqref{eq:QcontractN} implies the scalar Dini inequality
\(D^+Q(t)\le-\delta_Lm_*Q(t).\)
The Dini-version of Gronwall's lemma therefore implies, for every $s<t$,
\begin{equation}\label{eq:QGronwallN}
 Q(t)\le e^{-\delta_Lm_*(t-s)}Q(s).
\end{equation}
The entire solution is uniformly bounded above and below, so
\(0\le Q(s)\le\log\frac{\sup_{t,x}U(t,x)}{m_*}=:Q^*<\infty\)
for all $s\in\mathbb R$.  We let $s\to-\infty$ in \eqref{eq:QGronwallN} and obtain $Q(t)=0$ for every $t$.  Hence $M(t)=m(t)$ and $U(t,x)=P(t)$ is independent of $x$.

It remains to identify the positive bounded entire solution $P$ of $P'=Pf_-(t,P)$.  We let $p_-$ be the periodic solution from \cref{lem:pminus} and set
\(z(t)=\log\frac{P(t)}{p_-(t)}.\)
Then
\(z'(t)=f_-(t,P(t))-f_-(t,p_-(t)).\)
We put $m_P=\inf_tP(t)>0$ and $m_p=\min_tp_-(t)>0$.  If $z(t)>0$, then $P>p_-$ and \eqref{eq:fuass} provides
\(z'(t)\le-\beta(t)(P-p_-) =-\beta(t)p_-(e^z-1) \le-\beta_*m_p z.\)
If $z(t)<0$, then $P<p_-$ and
\(z'(t)\ge \beta(t)(p_--P) =\beta(t)P(e^{-z}-1) \ge \beta_*m_P(-z).\)
Thus the upper Dini derivative of $|z|$ satisfies
\(D^+|z(t)|\le-\beta_L|z(t)|, \; \beta_L=\beta_*\min\{m_P,m_p\}>0.\)
As above, $|z(t)|\le e^{-\beta_L(t-s)}|z(s)|$ for $s<t$.  Since $P$ and $p_-$ are bounded above and below on the whole time axis, letting $s\to-\infty$ implies $z\equiv0$.  Therefore $P=p_-$.  Finally the elliptic equation with spatially constant $P$ provides $V=(\mu/\nu)p_-=q_-$.
\end{proof}

\section{Proof of the sharp geometric propagation results}\label{sec:frontproof}
We first prove the compact-support theorem.

\begin{proof}[\upshape\bfseries Proof of \cref{thm:front}]
The global existence statement is \cref{prop:global}.

\emph{Step 1: exclusion of nonzero entire solutions in the unfavorable periodic far field.}
We fix a temporal phase $\theta\in[0,T]$.  We consider a bounded entire nonnegative solution $(U,V)$ of
\begin{equation}\label{eq:unfavlimitN}
\begin{cases}
 U_t=a(t+\theta)\Delta U-\eta\chi(t+\theta)\nabla\!\cdot(U\nabla V)
 +U f_+(t+\theta,U),\\
 d(t+\theta)\Delta V-\nu(t+\theta)V+\mu(t+\theta)U=0.
\end{cases}
\end{equation}
We let $M(t)=\sup_xU(t,x)$.  The maximum calculation from \cref{prop:global} implies
\(D^+M(t)\le r_+(t+\theta)M(t).\)
For any integer $n\ge1$,
\[
 M(t)\le M(t-nT)
 \exp\!\left(\int_{t-nT}^tr_+(s+\theta)\,ds\right)
 \le C e^{nT\langle r_+\rangle}.
\]
Because $\langle r_+\rangle<0$ and the entire solution is bounded, letting $n\to\infty$ provides $M(t)=0$.  Thus the only bounded entire nonnegative solution of \eqref{eq:unfavlimitN} is zero.

\emph{Step 2: extinction when $c<-c_*$.}
We choose $s$ with
\(c_*<s<-c.\)
By \cref{cor:uppercone}, $u(t,x)\to0$ uniformly for $|x|\ge st$.  On $|x|\le st$,
\(x\cdot e-ct\ge-|x|-ct\ge(-s-c)t\to+\infty.\)
If uniform extinction failed in the remaining cone, there would be $t_n\to\infty$ and $x_n$ with $|x_n|\le st_n$ and $u(t_n,x_n)\ge\delta>0$.  We pass to a subsequence so that $t_n$ modulo $T$ converges to $\theta$.  We define
\(u_n(t,y)=u(t+t_n,y+x_n), \; v_n(t,y)=v(t+t_n,y+x_n).\)
Uniform bounds, interior parabolic estimates for $u_n$, and elliptic estimates for $v_n$ imply local convergence to a bounded entire solution of \eqref{eq:unfavlimitN}; the convergence of the reaction follows from
\(x_n\cdot e-ct_n\to+\infty.\)
The limit has value at least $\delta$ at $(0,0)$, contradicting Step~1.  Hence extinction holds.

\emph{Step 3: persistence when $-c_*<c<c_*$.}
We fix $\eps$ as in the theorem and define the compact set of frame velocities
\begin{equation}\label{eq:Kepsmoderate}
 K_\eps=
 \{s\in\mathbb R^N:|s|\le c_*-\eps,\ s\cdot e\le c-\eps\}.
\end{equation}
The assumptions on $\eps$ ensure that $K_\eps$ is nonempty and $K_\eps\Subset\Omega_c$.  By \cref{lem:framepersist}, $u(t,st)$ has a positive lower limit uniformly for $s\in K_\eps$.  Every $x\in\mathcal K_\eps(t)$ can be written as $x=st$ with $s=x/t\in K_\eps$.  This proves the persistence assertion in \eqref{eq:frontmoderate}.

The first decay assertion in \eqref{eq:frontmoderate} is \cref{cor:uppercone}.  To prove the second decay assertion in \eqref{eq:frontmoderate}, suppose it fails.  Then there exist $t_n\to\infty$ and $x_n$ with $x_n\cdot e\ge(c+\eps)t_n$ and $u(t_n,x_n)\ge\delta>0$.
Hence $x_n\cdot e-ct_n\to+\infty$.  Translation and compactness again produce a nonzero bounded entire solution of the unfavorable limiting system, contradicting Step~1.

\emph{Step 4: persistence when $c\ge c_*$.}
For $\eps\in(0,c_*)$, set
\(K_\eps=\{s\in\mathbb R^N:|s|\le c_*-\eps\}.\)
For every $s\in K_\eps$,
\(s\cdot e\le|s|\le c_*-\eps<c,\)
with a uniform margin even at $c=c_*$.  Hence $K_\eps\Subset\Omega_c$.  We apply \cref{lem:framepersist} at $x=st$ and obtain \eqref{eq:frontdecayfastgood}.  The decay \eqref{eq:frontdecayfastgood} again follows from \cref{cor:uppercone}.

\emph{Step 5: convergence to the favorable periodic state.}
We assume \eqref{eq:strongdamp}.  We prove \eqref{eq:frontconvmoderate} and \eqref{eq:frontconvfastgood} simultaneously.  For the fixed $\eps$ in the corresponding part of the theorem, define
\[
 \mathcal K_\eps=
 \begin{cases}
 \{s\in\mathbb R^N:\ |s|\le c_*-\eps,\ s\cdot e\le c-\eps\},
 &-c_*<c<c_*,\\[1mm]
 \{s\in\mathbb R^N:\ |s|\le c_*-\eps\},
 &c\ge c_*.
 \end{cases}
\]
In the first case $\mathcal K_\eps=K_\eps$ from \eqref{eq:Kepsmoderate}.  In the second case,
\(s\cdot e-c\le |s|-c\le c_*-\eps-c\le-\eps \;(s\in\mathcal K_\eps).\)
Thus in both cases $\mathcal K_\eps\Subset\Omega_c$ and, for every $x=st$ with $s\in\mathcal K_\eps$,
\begin{equation}\label{eq:favmarginconvN}
 x\cdot e-ct=(s\cdot e-c)t\le-\eps t.
\end{equation}

We suppose that the convergence asserted in the relevant formula, \eqref{eq:frontconvmoderate} or \eqref{eq:frontconvfastgood}, fails.  Then there are $t_n\to\infty$ and $x_n=t_ns_n$ with $s_n\in\mathcal K_\eps$ such that
\begin{equation}\label{eq:convcontradN}
 |u(t_n,x_n)-p_-(t_n)|+|v(t_n,x_n)-q_-(t_n)|\ge\varepsilon_0.
\end{equation}
We pass to a subsequence such that
\(s_n\to s_\infty\in\mathcal K_\eps, \; t_n\pmod T\to\theta\in[0,T].\)
By \eqref{eq:favmarginconvN}, for every fixed $R>0$ and $|t|+|y|\le R$,
\((x_n+y)\cdot e-c(t_n+t) \le-\eps t_n+R(1+|c|)\longrightarrow-\infty.\)
The convergence in \cref{F}, positive-time parabolic estimates for $u$, and elliptic Schauder estimates for $v$ therefore imply, after diagonal extraction,
\(\big(u(t_n+t,x_n+y),v(t_n+t,x_n+y)\big) \longrightarrow(U(t,y),V(t,y))\)
locally in the classical topology on $\mathbb R\times\mathbb R^N$.  The limit is a bounded entire solution of the favorable limiting system with coefficients shifted by the temporal phase $\theta$.

It remains to show that $U$ is bounded away from zero on the whole space-time cylinder.  Since $\mathcal K_\eps$ is compactly contained in $\Omega_c$, we choose a compact set
\(\mathcal K_\eps\Subset K'\Subset\Omega_c.\)
We let $\rho_{K'}>0$ be supplied by \cref{lem:framepersist}.  We fix $R>0$.  For $|t|\le R$ and $|y|\le R$ define
\(s_{n,t,y}:=\frac{x_n+y}{t_n+t}.\)
For $n$ so large that $t_n>2R$,
\[
 |s_{n,t,y}-s_n|=\left|\frac{y-ts_n}{t_n+t}\right|
 \le\frac{R+R\sup_n|s_n|}{t_n-R}\longrightarrow0,
\]
uniformly on the cylinder.  Hence $s_{n,t,y}\in K'$ for all large $n$.  Since
\(x_n+y=s_{n,t,y}(t_n+t),\)
the definition of the liminf in \eqref{eq:framepositiveN} provides a time $T_{K'}$ such that
\(u(t_n+t,x_n+y)\ge\frac12\rho_{K'}\)
whenever $n$ is sufficiently large, uniformly for $|t|\le R$ and $|y|\le R$.  We pass to the limit and obtain
\(U(t,y)\ge\frac12\rho_{K'} \;(|t|\le R,\ |y|\le R).\)
The same constant works for every $R$, hence
\(\inf_{\mathbb R\times\mathbb R^N}U\ge\frac12\rho_{K'}>0.\)
By \cref{lem:Liouville},
\(U(t,y)=p_-(t+\theta), \; V(t,y)=q_-(t+\theta).\)
We evaluate at $(t,y)=(0,0)$ and contradict \eqref{eq:convcontradN}.  Consequently both \eqref{eq:frontconvmoderate} and \eqref{eq:frontconvfastgood} hold in their stated parameter regimes.
\end{proof}

\begin{proof}[\upshape\bfseries Proof of \cref{thm:frontfilled}]
We treat decay, persistence and convergence separately.

If $c\le c_*$, then $s_c=c$.  We suppose \eqref{eq:filleddecay} fails.  There exist $\delta>0$, $t_n\to\infty$ and $x_n$ such that
\(x_n\cdot e\ge(c+\eps)t_n, \; u(t_n,x_n)\ge\delta.\)
For every fixed $R>0$, $|t|\le R$ and $|y|\le R$,
\((x_n+y)\cdot e-c(t_n+t) \ge\eps t_n-R(1+|c|)\longrightarrow+\infty.\)
Positive-time parabolic estimates for $u$, the elliptic estimates for $v$, and periodicity of the coefficients allow us to choose a subsequence for which $t_n$ modulo $T$ tends to $\theta$ and
\(\big(u(t_n+t,x_n+y),v(t_n+t,x_n+y)\big) \longrightarrow(U_*,V_*)\)
locally in the classical topology.  The far-field convergence in \cref{F} shows that $(U_*,V_*)$ is a bounded entire solution of the unfavorable limiting system and $U_*(0,0)\ge\delta$.

We let $M_*(t)=\sup_yU_*(t,y)$.  The approximate-maximum computation from \cref{prop:global}, with $r_+(t+\theta)$ in place of $r(t,\xi)$, yields
\(D^+M_*(t) \le r_+(t+\theta)M_*(t)-\delta_0M_*(t)^2 \le r_+(t+\theta)M_*(t).\)
For every integer $m\ge1$,
\[
 M_*(t)\le M_*(t-mT)
 \exp\!\left(\int_{t-mT}^{t}r_+(s+\theta)\,ds\right)
 \le \|U_*\|_\infty e^{mT\langle r_+\rangle}.
\]
Since $\langle r_+\rangle<0$, letting $m\to\infty$ provides $M_*(t)=0$, contradicting $U_*(0,0)\ge\delta$.  Hence \eqref{eq:filleddecay} holds when $c\le c_*$.

If $c>c_*$, then $s_c=c_*$.  We choose $s$ with $c_*<s<c_*+\eps$.  The hypothesis $u_0(x)=0$ for $x\cdot e\ge R_0$ is exactly the assumption of \cref{lem:halfupper}; hence
\(\sup_{x\cdot e\ge st}u(t,x)\longrightarrow0.\)
Since $s<c_*+\eps=s_c+\eps$, the half-space
$\{x\cdot e\ge(s_c+\eps)t\}$ is contained in
$\{x\cdot e\ge st\}$ for every $t>0$.  This proves \eqref{eq:filleddecay}.

For persistence, we fix $s<s_c$.  Then $s<\min\{c,c_*\}$, so \cref{lem:halffilled} provides a positive lower limit uniformly on $x\cdot e\le st$.
We take $s=s_c-\eps$ and obtain the persistence assertion in \eqref{eq:filleddecay}.

We assume now \eqref{eq:strongdamp}.  If \eqref{eq:filledconv} failed, we choose $t_n\to\infty$ and $x_n\cdot e\le(s_c-\eps)t_n$ where the difference from $(p_-,q_-)$ stays above a fixed positive number.  Since $s_c\le c$,
\(x_n\cdot e-ct_n\le-\eps t_n.\)
The half-space persistence lemma, applied with a slightly larger speed than $s_c-\eps$, provides a positive lower bound independent of the translated centers on each fixed neighborhood of $(t_n,x_n)$.  Translation compactness therefore produces a positive bounded entire favorable limit.  The Liouville lemma identifies that limit with $(p_-,q_-)$, a contradiction.  This proves \eqref{eq:filledconv}.
\end{proof}

\section{Moving favorable slabs}\label{sec:localized}
We begin with the spectral facts needed for \cref{thm:localized}.

\begin{lemma}[Planar exhaustion and positive growth on a smooth domain]\label{lem:exhaust}
Under \cref{A1,L}, the principal exponent $\lambda_L(c)$ in
\eqref{eq:FloquetL} is strictly increasing in $L$ and converges to a finite
limit $\lambda_\infty(c)$.  If $\lambda_\infty(c)<0$, a sequence of
normalized eigenfunctions converges locally to a positive solution $\Phi$
of \eqref{eq:wholeFloquet}.

If $\lambda_\infty(c)>0$, then for every prescribed $R_0>0$ there exists a
smooth ball $\Omega=B_S(0)\Subset\mathbb R^N$ containing $\overline{B_{R_0}(0)}$
whose periodic Dirichlet principal exponent for
$a(t)\Delta+ce\cdot\nabla+r(t,y\cdot e)$ is strictly positive.
\end{lemma}

\begin{proof}
The periodic parabolic Krein--Rutman theorem on $(-L,L)$ ensures a simple principal exponent and a positive eigenfunction.  The exponent increases strictly with $L$ by domain monotonicity.  Since $r$ is bounded, the exponents are bounded above, so the limit exists.  We normalize $\phi_L(0,0)=1$ whenever the interval contains the origin.  Interior parabolic Schauder estimates and Harnack inequalities provide local compactness; a diagonal subsequence converges to a positive generalized eigenfunction solving \eqref{eq:wholeFloquet}.

We now suppose $\lambda_\infty(c)>0$ and prescribe $R_0>0$.
Choose $L>R_0+1$ and $\sigma>0$ such that $\lambda_L(c)>3\sigma$.
When $N=1$, take the smooth interval $\Omega=(-L,L)$; its exponent is
$\lambda_L(c)>0$.
For $N\ge2$, let $\psi_R$ be the positive principal Dirichlet
eigenfunction of $-\Delta$ on the transverse ball $B_R^{N-1}$,
with eigenvalue $\mu_R=R^{-2}\mu_1\downarrow0$.
Write $y=\xi e+z$, where $z\in e^\perp$, and set
\[
 Q_R(t)=\exp\!\left(-\mu_R\int_0^t[a(s)-\bar a] \,ds\right),\qquad
 \Phi_{L,R}(t,y)=Q_R(t)\phi_L(t,\xi)\psi_R(z).
\]
On the auxiliary product domain
$D_{L,R}=\{y:|y\cdot e|<L,\ |y-(y\cdot e)e|<R\}$,
this function is positive and $T$-periodic, vanishes on the boundary,
and has principal exponent $\lambda_L(c)-\bar a\mu_R$.
Choose $R>R_0+1$ large enough that $\bar a\mu_R<\sigma$.
Then the principal exponent on $D_{L,R}$ exceeds $2\sigma$.
The product domain is used here \emph{only to calculate an eigenvalue}:
its boundary has edges, so a global smooth-boundary Hopf estimate must not
be inferred from this separation of variables.
Take a sufficiently large ball $\Omega=B_S(0)$ with
$\overline{D_{L,R}}\cup\overline{B_{R_0}(0)}\subset\Omega$.
Domain monotonicity of the positive periodic Dirichlet evolution implies
$\lambda_\Omega(c)\ge\lambda_{D_{L,R}}(c)>2\sigma$.
The ball has smooth boundary, so the parabolic Hopf lemma and boundary
gradient estimates apply to its positive principal eigenfunction without
the corner difficulty of the auxiliary product domain.
\end{proof}

\begin{proof}[\upshape\bfseries Proof of \cref{thm:localized}]
Move with the climate by setting
\(U(t,y)=u(t,y+ct\,e), \; V(t,y)=v(t,y+ct\,e).\)
Then
\begin{equation}\label{eq:localizedmovingN}
\begin{cases}
 U_t=a(t)\Delta U+ce\cdot\nabla U
 -\eta\chi(t)\nabla\!\cdot(U\nabla V)+Uf(t,y\cdot e,U),\\
 d(t)\Delta V-\nu(t)V+\mu(t)U=0.
\end{cases}
\end{equation}

\emph{Part (i): longitudinal confinement.}
We fix $\sigma>0$ and suppose that \eqref{eq:longitudinalconfinement} fails.  Then there are $\varepsilon_0>0$, $t_n\to\infty$ and $y_n\in\mathbb R^N$ such that
\(|y_n\cdot e|\ge\sigma t_n, \; U(t_n,y_n)\ge\varepsilon_0.\)
We set $\xi_n=y_n\cdot e$.  We pass to a subsequence such that either $\xi_n\to+\infty$ or $\xi_n\to-\infty$, while $t_n$ modulo $T$ converges to some $\theta\in[0,T]$.  We define
\(U_n(t,z)=U(t+t_n,z+y_n), \; V_n(t,z)=V(t+t_n,z+y_n).\)
The uniform bound from \cref{prop:global}, positive-time parabolic Schauder estimates for $U$, and elliptic Schauder estimates for $V$ imply local compactness on every bounded space--time cylinder.  By the uniform local H\"older bounds in \cref{A2}, the translated reactions
\(f_n(t,z,q):=f(t+t_n,z\cdot e+\xi_n,q)\)
are locally precompact on $\mathbb R\times\mathbb R^N\times[0,M_0]$.  After a further extraction, $(U_n,V_n,f_n)$ converges locally to $(U_*,V_*,\widehat f_{\diamond})$, where $(U_*,V_*)$ is a bounded nonnegative entire solution and $U_*(0,0)\ge\varepsilon_0$.  If $\xi_n\to+\infty$, we set $r_{\diamond}=r_+$; if $\xi_n\to-\infty$, we set $r_{\diamond}=r_-$.  \cref{A2,L} imply, for every $t,z$ and $0\le q\le M_0$,
\[
 \widehat f_{\diamond}(t,z,0)=r_{\diamond}(t+\theta),
 \;
 \widehat f_{\diamond}(t,z,q)\le r_{\diamond}(t+\theta)-\beta(t+\theta)q.
\]
Thus no nonlinear far-field limit of $f$ is required in the localized-slab theorem.

We let $M_*(t)=\sup_{z\in\mathbb R^N}U_*(t,z)$.  When the supremum is not attained, we use the quadratic penalization argument detailed in \cref{prop:global}.  The maximum-point computation then passes to the limit, and we obtain
\[
 D^+M_*(t)\le r_{\diamond}(t+\theta)M_*(t)-\delta_0M_*(t)^2
 \le r_{\diamond}(t+\theta)M_*(t).
\]
For every integer $m\ge1$ and every $t\in\mathbb R$, integration from $t-mT$ to $t$ implies
\[
 M_*(t)\le M_*(t-mT)
 \exp\!\left(\int_{t-mT}^{t}r_{\diamond}(s+\theta)\,ds\right)
 =M_*(t-mT)e^{mT\langle r_{\diamond}\rangle}.
\]
The entire solution is bounded and \cref{L} provides $\langle r_{\diamond}\rangle<0$; letting $m\to\infty$ forces $M_*(t)=0$.  This contradicts $U_*(0,0)\ge\varepsilon_0$ and proves \eqref{eq:longitudinalconfinement}.

\emph{Part (ii): extinction from a nonpositive maximal mean or an excessive climate speed.}
At a spatial maximum of $u$, the computation in \cref{prop:global} provides
\begin{equation}\label{eq:rsupmaxN}
 D^+M(t)
 \le r^{\mathrm{sup}}(t)M(t)-\delta_0M(t)^2,
 \; M(t)=\lVert u(t,\cdot)\rVert_\infty.
\end{equation}
If $\bar r^{\mathrm{sup}}<0$, dropping the negative quadratic term in \eqref{eq:rsupmaxN} provides
\(D^+M(t)\le r^{\mathrm{sup}}(t)M(t).\)
Hence, for $t\ge s\ge0$,
\(M(t)\le M(s)\exp\!\left(\int_s^t r^{\mathrm{sup}}(\zeta)\,d\zeta\right).\)
We write $t-s=mT+\vartheta$ with $m\in\mathbb N$ and $0\le\vartheta<T$.  Periodicity provides
\[
 \int_s^t r^{\mathrm{sup}}(\zeta)\,d\zeta
 \le mT\bar r^{\mathrm{sup}}+T\|r^{\mathrm{sup}}\|_\infty
 \le-\omega(t-s)+C
\]
for some $\omega>0$ and $C<\infty$.  Thus $M(t)\le C_1e^{-\omega t}$ after changing $C_1$.  If $\bar r^{\mathrm{sup}}=0$ and $M$ vanishes at some time, uniqueness provides $u\equiv0$ afterwards and there is nothing to prove.  Otherwise put $Y(t)=1/M(t)$.  The upper Dini inequality implies
\(D_-Y(t)\ge-r^{\mathrm{sup}}(t)Y(t)+\delta_0.\)
We let
\(R(t)=\int_0^t r^{\mathrm{sup}}(s)\,ds.\)
Because $r^{\mathrm{sup}}$ is $T$-periodic and has zero mean, $R$ is $T$-periodic and bounded.  The product rule for Dini derivatives provides
\(D_-\big(e^{R(t)}Y(t)\big)\ge\delta_0e^{R(t)}.\)
Hence, for $t>s$,
\(e^{R(t)}Y(t)\ge e^{R(s)}Y(s) +\delta_0\int_s^t e^{R(\zeta)}\,d\zeta.\)
If $R_*:=\max_{[0,T]}|R|$, then
\(Y(t)\ge e^{-R_*}\delta_0e^{-R_*}(t-s) =\delta_0e^{-2R_*}(t-s)\)
up to the harmless positive initial term.  Thus $Y(t)\to\infty$ and $M(t)=Y(t)^{-1}\to0$.

We assume $\bar r^{\mathrm{sup}}>0$ and $|c|>c_D^{\rm env}$.  We choose $\lambda\in\mathcal A_D$ and a number $s$ such that
\(\bar a\lambda+\frac{\bar r^{\mathrm{sup}}}{\lambda}<s<|c|.\)
We apply the first-contact calculation of \cref{lem:superenv} with $r^{\mathrm{sup}}(t)$ in place of $r_-(t)$.  More explicitly, with the selected $\lambda$ define
\(c_{\lambda,\mathrm{sup}}=\bar a\lambda+\frac{\bar r^{\mathrm{sup}}}{\lambda}<s\)
and
\[
 P_{\lambda,\mathrm{sup}}(t)=
 \exp\!\left(\int_0^t
 [a(\zeta)\lambda^2+r^{\mathrm{sup}}(\zeta)-\lambda c_{\lambda,\mathrm{sup}}]\,d\zeta\right).
\]
For any $\widehat c\in(c_{\lambda,\mathrm{sup}},s)$ and every unit vector $\theta$, a sufficiently large multiple of
\(P_{\lambda,\mathrm{sup}}(t)e^{-\lambda(x\cdot\theta-\widehat c t)}\)
dominates the compactly supported datum.  At a first contact, $\lambda\in\mathcal A_D$ provides
$\gamma_{N,\lambda}(t)\le\beta(t)$, so the chemotactic quadratic term is absorbed exactly as in \eqref{eq:GammausedN}; the speed excess
$\lambda(\widehat c-c_{\lambda,\mathrm{sup}})>0$ implies the strict contradiction.  We take the supremum over directions and then we choose $\widehat c<s$.  We obtain
\begin{equation}\label{eq:locupperconeN}
 \lim_{t\to\infty}\sup_{|x|\ge st}u(t,x)=0.
\end{equation}
On $|x|\le st$, if $c>0$ then
\(x\cdot e-ct\le(s-c)t\to-\infty,\)
whereas if $c<0$ then
\(x\cdot e-ct\ge(-s-c)t=(|c|-s)t\to+\infty.\)
If extinction failed inside the cone, translation around a positive sequence and the local compactness argument used in part~\textup{(i)} would imply a nonzero bounded entire limit.  The translated reactions need not have nonlinear far-field limits under \cref{L}; nevertheless every subsequential limit $\widehat f$ satisfies $\widehat f(t,z,q)\le r_\diamond(t+\theta)-\beta(t+\theta)q$ with $\langle r_\diamond\rangle<0$.  The maximum argument of part~\textup{(i)} therefore forces the limiting population to vanish.  This contradiction and \eqref{eq:locupperconeN} prove part~\textup{(ii)}.

\emph{Part (iii): persistence when $\lambda_\infty(c)>0$.}
Fix $R_0>0$.  By \cref{lem:exhaust}, choose a smooth ball
$\Omega=B_S(0)$ that contains $\overline{B_{R_0}(0)}$ and has positive
periodic Dirichlet exponent $\lambda_\Omega>8\sigma$ for some
$\sigma>0$.  Let $\varphi(t,y)$ be its strictly positive $T$-periodic
principal eigenfunction, normalized by $\|\varphi\|_\infty=1$.  We fix $0<\delta<\sigma$.  Since subtracting the constant $3\delta$ from the zero-order coefficient subtracts exactly $3\delta$ from the periodic principal exponent, the unperturbed Dirichlet evolution generated by
\(a(t)\Delta+ce\cdot\nabla+r(t,y\cdot e)-3\delta\)
has one-period growth at least $e^{5\sigma T}$ on its positive principal Floquet mode.

We next quantify the admissible drift perturbation.  Positive-time parabolic estimates for $U$, elliptic Schauder estimates for $V$, and the global bound from \cref{prop:global} imply a constant $H_{\mathcal C}>0$ such that
\begin{equation}\label{eq:locdriftholderN}
 \sup_{t_0\ge1}
 \|\eta\chi(t)\nabla V\|_{C^{\alpha_H/2,\alpha_H}([t_0,t_0+T]\times\overline{\Omega})}
 \le H_{\mathcal C}.
\end{equation}
We choose a smooth radial truncation $\mathfrak T_{\eta_0}$ exactly as in \eqref{eq:truncAN}, with a Lipschitz constant bounded by the same universal number $C_T$.  On the fixed smooth ball $\Omega$, the Dirichlet evolution has the same $L^\infty\to C^1$ gradient bound of order $(t-s)^{-1/2}$ used in \cref{lem:Dirichletpositive}; constants depend on the fixed ball.  The fractional-Gronwall and Duhamel argument of \cref{lem:Dirichletpositive}, applied with the normalized Dirichlet eigenfunction on $\Omega$ as initial datum, provides a constant $C_{\mathcal C}=C(H_{\mathcal C},C_T)$ such that the one-period map with truncated drift $A$ satisfies
\[
 \|\mathcal U_A(t_0+T,t_0)\varphi(t_0)-
 \mathcal U_0(t_0+T,t_0)\varphi(t_0)\|_{C^1(\overline{\Omega})}
 \le C_{\mathcal C}\|A\|_\infty
\]
uniformly in $t_0$.  Because $\Omega$ has smooth boundary,
$\varphi(t,y)\ge c_H\operatorname{dist}(y,\partial\Omega)$ uniformly in
$t$ by the Hopf lemma.  Both evolutions in the preceding estimate vanish
on $\partial\Omega$, so their $C^1$ difference is bounded pointwise by
$C_{\mathcal C}\|A\|_\infty\operatorname{dist}(y,\partial\Omega)$.
This is at most $(C_{\mathcal C}/c_H)\|A\|_\infty\varphi(t,y)$.
Hence we may choose $\eta_0>0$ so small that every drift satisfying
\[
 \|A\|_\infty\le\eta_0,
 \;
 \sup_{t_0\in\mathbb R}\|A\|_{C^{\alpha_H/2,\alpha_H}([t_0,t_0+T]\times\overline{\Omega})}
 \le C_TH_{\mathcal C}
\]
preserves the quantitative amplification
\begin{equation}\label{eq:loconeperiodN}
 \mathcal U_A(t_0+T,t_0)\varphi(t_0,\cdot)
 \ge e^{3\sigma T}\varphi(t_0+T,\cdot).
\end{equation}
This is the smooth-domain analogue of the perturbation estimate proved in \cref{lem:Dirichletpositive}; the displayed $C^1$ estimate also shows explicitly why an $L^\infty$ drift bound alone would not suffice without the common H\"older control in \eqref{eq:locdriftholderN}.

Now apply \cref{cor:smalldrift}.  First choose the Yukawa localization radius so large that its tail contribution, after multiplication by the fixed coefficient bounds, is at most $\min\{\delta,\eta_0\}/4$.  Next fix the Harnack exponent $p$ furnished by \cref{lem:Harnack}, and finally choose $q>0$ so small that the local contribution is at most $\min\{\delta,\eta_0\}/4$.  Then, for $t\ge1$,
\begin{equation}\label{eq:locsmallN}
 U(t,y)\le q
 \;\Longrightarrow\;
 \begin{cases}
 \eta\chi(t)|\nabla V(t,y)|\le\dfrac{\eta_0}{2},\\[1mm]
 \eta\chi(t)\kappa(t)^2V(t,y)\le\delta.
 \end{cases}
\end{equation}
We define
\[
 A(t,y)=\mathfrak T_{\eta_0}\big(\eta\chi(t)\nabla V(t,y)\big)
 \;(t\ge1,\ y\in\Omega).
\]
Then $|A|\le\eta_0$, its H\"older norm is bounded by $C_TH_{\mathcal C}$, and, crucially,
\begin{equation}\label{eq:locinactiveN}
 A(t,y)=\eta\chi(t)\nabla V(t,y)
 \;\text{whenever }U(t,y)\le q.
\end{equation}
We extend $A$ to the whole time axis exactly as in the proof of \cref{lem:framepersist}, preserving these two uniform bounds.

We construct the auxiliary entire state rather than merely citing the ball argument.  We fix $C_*>0$ and let $b_1$ be the quadratic KPP bound on $[0,M_0]$.  We fix the comparison time $t_{\rm cmp}=1$.  The strong maximum principle provides
\(m_{\rm cmp}:= \min_{y\in\overline{\Omega}}U(t_{\rm cmp},y)>0.\)
If $r_{\max}:=\sup_{t,\xi}r(t,\xi)$ and $Q_B>0$ is determined by
\(BQ_B+C_*Q_B^{1/p}=(r_{\max}-3\delta)_++1,\)
then $Q_B\downarrow0$ as $B\to\infty$.  We choose once and for all
\begin{equation}\label{eq:locBchoiceN}
 B>b_1+1
 \;\text{so large that}\;
 Q_B<\frac12\min\{q,m_{\rm cmp}\}.
\end{equation}
Having fixed this $B$, we choose $q_1\in(0,Q_B)$ such that
\begin{equation}\label{eq:locnonlinlossN}
 C_*q_1^{1/p}+Bq_1\le\sigma.
\end{equation}
For $n\in\mathbb N$ solve on $[-nT,\infty)\times\Omega$
\begin{equation}\label{eq:locauxN}
\begin{cases}
 w_t=a(t)\Delta w+(ce-A)\cdot\nabla w\\
 \hspace{11mm}+w\big(r(t,y\cdot e)-3\delta-C_*w^{1/p}-Bw\big),\\
 w=0\;\text{on }\partial\Omega,\\
 w(-nT,\cdot)=\zeta\varphi(-nT,\cdot).
\end{cases}
\end{equation}
The constant $Q_B$ is a strict supersolution of \eqref{eq:locauxN}.  We put
\(R_{\mathcal C}:=(r_{\max}-3\delta)_+.\)
For the orbit issued at an arbitrary phase $t_0$ from $\zeta\varphi(t_0)$, the maximum principle provides
\(\|w(t,\cdot)\|_\infty \le e^{R_{\mathcal C}(t-t_0)}\zeta, \; t_0\le t\le t_0+T.\)
We choose
\(0<\zeta< \min\left\{Q_B,\frac12q_1e^{-R_{\mathcal C}T}\right\}.\)
Then the first nonlinear orbit remains below $q_1/2$ during the whole first period.  We let $z$ be the corresponding linear orbit for the operator in \eqref{eq:loconeperiodN} and set
\(\underline z(t)=e^{-\sigma(t-t_0)}z(t).\)
The same maximum estimate keeps $z$ below $q_1/2$.  By \eqref{eq:locnonlinlossN},
\begin{align*}
 (\underline z)_t
 &=a(t)\Delta\underline z+(ce-A)\cdot\nabla\underline z
 +(r(t,y\cdot e)-3\delta-\sigma)\underline z\\
 &\le a(t)\Delta\underline z+(ce-A)\cdot\nabla\underline z
 +\underline z\big(r(t,y\cdot e)-3\delta
 -C_*\underline z^{1/p}-B\underline z\big).
\end{align*}
Hence $\underline z$ is a subsolution of \eqref{eq:locauxN}.  Scalar comparison and \eqref{eq:loconeperiodN} provide the explicit one-period amplification
\begin{equation}\label{eq:locnonlinearampN}
 \mathcal S(t_0+T,t_0)[\zeta\varphi(t_0,\cdot)]
 \ge e^{2\sigma T}\zeta\varphi(t_0+T,\cdot)
 \ge\zeta\varphi(t_0+T,\cdot).
\end{equation}
Thus the pullback solutions issued from the same periodic seed at times $-nT$ form a monotone increasing sequence at every fixed later phase.  They are bounded from above by the strict supersolution $Q_B$.  Parabolic compactness therefore implies a positive bounded entire limit $w$ satisfying
\begin{equation}\label{eq:locauxboundsN}
 0<w<Q_B<q/2,
 \;
 \inf_{t\in\mathbb R,\ y\in\mathcal C'}w(t,y)>0
\end{equation}
for every compact subdomain $\mathcal C'\Subset\Omega$.  Indeed, \eqref{eq:locnonlinearampN} shows that the pullback limit dominates the fixed positive seed at every integer phase; the parabolic Harnack inequality then propagates a uniform lower bound through each period on every compact subdomain.

By \eqref{eq:locBchoiceN},
\[
 w(t_{\rm cmp},y)<Q_B<m_{\rm cmp}\le U(t_{\rm cmp},y)
 \;(y\in\overline{\Omega}).
\]
We suppose that $Z=U-w$ has a first zero at $(t_1,y_1)$.  The contact is interior because $w=0$ on the artificial lateral boundary while $U>0$ there.  At the contact,
\(U=w<q, \; \nabla Z=0, \; \Delta Z\ge0, \; Z_t\le0.\)
By \eqref{eq:locinactiveN}, the clipped drift equals the actual drift at this point.  Expanding the chemotactic divergence in \eqref{eq:localizedmovingN} and using $d\Delta V=\nu V-\mu U$ provides
\[
 U_t=a(t)\Delta U+(ce-\eta\chi\nabla V)\cdot\nabla U
 +U\left[f(t,y\cdot e,U)-\eta\chi\kappa^2V+\frac{\eta\chi\mu}{d}U\right].
\]
For $U\le q$, \eqref{eq:kppbounds} and \eqref{eq:locsmallN} imply, after discarding the nonnegative last quadratic term,
\(U_t\ge a(t)\Delta U+(ce-A)\cdot\nabla U +U\big(r(t,y\cdot e)-\delta-b_1U\big).\)
We subtract \eqref{eq:locauxN} at the contact and obtain
\[
 Z_t(t_1,y_1)\ge a(t_1)\Delta Z+(ce-A)\cdot\nabla Z
 +w\big(2\delta+C_*w^{1/p}+(B-b_1)w\big)>0,
\]
provided $B>b_1$.  This contradicts $Z_t\le0$.  Hence $U\ge w$ for all $t\ge t_{\rm cmp}$.

Because $B_{R_0}(0)\Subset\Omega$ was built into the choice of the smooth auxiliary domain, \eqref{eq:locauxboundsN} and the comparison $U\ge w$ imply
\(\liminf_{t\to\infty}\inf_{|y|\le R_0}U(t,y)>0.\)
Since $R_0>0$ was arbitrary, this proves \eqref{eq:localpersistN}.

\emph{Part (iv): extinction when $\lambda_\infty(c)<0$.}
We let $\Phi$ solve \eqref{eq:wholeFloquet}.  We fix
\(0<\varepsilon_\lambda<-\lambda_\infty(c)\)
and define the slightly strict planar supersolution
\(\overline U(t,y)=M e^{(\lambda_\infty(c)+\varepsilon_\lambda)t}\Phi(t,y\cdot e).\)
Because $u_0$ is compactly supported and $\Phi(0,\xi)>0$, one can choose $M$ so that $U(0,y)<\overline U(0,y)$ for all $y$.  By \eqref{eq:wholeFloquet},
\begin{equation}\label{eq:UbarlinearN}
 \overline U_t=a(t)\Delta\overline U
 +ce\cdot\nabla\overline U+r(t,y\cdot e)\overline U
 +\varepsilon_\lambda\overline U.
\end{equation}
The term $\varepsilon_\lambda\overline U$ ensures strictness at a first contact.

We show that failure of the whole-space comparison would produce an actual first contact.  The logarithmic-slope bound provides
\(\Phi(t,\xi)\ge \Phi(t,0)e^{-K_\Phi|\xi|} \;(t\in\mathbb R,\ \xi\in\mathbb R),\)
and periodicity provides \(m_\Phi:=\min_{t\in[0,T]}\Phi(t,0)>0\).  Hence, on every finite time interval \(0\le t\le T_1\),
\begin{equation}\label{eq:PhiLowerExp}
 \overline U(t,y)\ge
 M m_\Phi e^{(\lambda_\infty(c)+\varepsilon_\lambda)T_1}
 e^{-K_\Phi|y\cdot e|}.
\end{equation}
On the other hand, the scalar form of the population equation has bounded drift and bounded zero-order coefficient by \cref{prop:global}.  Since the initial datum is compactly supported, the upper Aronson estimate implies
\[
 U(t,y)\le C_{T_1}
 \exp\!\left(-\frac{(|y|-C_{T_1})_+^2}{C_{T_1}}\right),
 \; 0<t\le T_1.
\]
The Gaussian decay dominates the exponential lower bound in \eqref{eq:PhiLowerExp}; therefore $U/\overline U\to0$ as $|y|\to\infty$, uniformly for $t$ in compact positive time intervals.  Starting from strict domination at \(t=0\), any loss of the inequality \(U\le\overline U\) must thus occur first at a finite interior point.

We suppose such a first contact with $U$ occurs.  Then $U=\overline U$, $\nabla U=\nabla\overline U$ and $\Delta U\le\Delta\overline U$.  At the fixed contact time the planar profile $\Phi(t,\xi)$ satisfies
\(|\Phi_\xi|\le K_\Phi\Phi.\)
If $N\ge2$, the hypothesis $K_\Phi<\kappa_*\le\kappa(t)$ places the ratio strictly inside the range of \cref{lem:loglipN}.  If $N=1$, the theorem also permits the endpoint $K_\Phi=\kappa_*$; since $\kappa_*\le\kappa(t)$ and the one-dimensional endpoint in \cref{lem:loglipN} is valid for $K/\kappa\le1$, the same estimate remains valid without a limiting argument.  The weighted resolvent estimate and \eqref{eq:localizeddamping} therefore bound the chemotactic term above by $\beta(t)\overline U^2$
at the contact.  Since $f(t,\xi,U)\le r(t,\xi)-\beta(t)U$, the nonlinear right-hand side of \eqref{eq:localizedmovingN} is no larger than
\[
 a(t)\Delta\overline U+ce\cdot\nabla\overline U+r(t,y\cdot e)\overline U
 =\overline U_t-\varepsilon_\lambda\overline U.
\]
At a first contact, $Z=U-\overline U$ satisfies $Z_t\ge0$, $\nabla Z=0$ and $\Delta Z\le0$, whereas the preceding inequality provides
\(Z_t\le a(t)\Delta Z-\varepsilon_\lambda\overline U<0,\)
a contradiction.  Hence
\begin{equation}\label{eq:localexpPhiN}
 U(t,y)\le M e^{(\lambda_\infty(c)+\varepsilon_\lambda)t}\Phi(t,y\cdot e)
 \;(t\ge0,\ y\in\mathbb R^N).
\end{equation}
Because $\lambda_\infty(c)+\varepsilon_\lambda<0$ and $\Phi$ is bounded on compact longitudinal intervals, \eqref{eq:localexpPhiN} provides exponential decay on every fixed ball in the moving frame.

To obtain uniform decay, suppose instead that $U(t_n,y_n)\ge\varepsilon_1>0$ for some $t_n\to\infty$.  The local decay implies $|y_n|\to\infty$.  If the longitudinal coordinates $\xi_n=y_n\cdot e$ remain bounded, \eqref{eq:localexpPhiN} itself provides a contradiction, regardless of the transverse components.  Thus, after a subsequence, $\xi_n\to+\infty$ or $\xi_n\to-\infty$.  Translate in time and space by $(t_n,y_n)$ and pass to a subsequence such that $t_n$ modulo $T$ converges to $\theta\in[0,T]$.  Positive-time parabolic estimates for $U$, elliptic estimates for $V$, and the uniform local H\"older bounds for $f$ imply local convergence of the translated solution and, after a further extraction, local convergence of the translated reactions to a function $\widehat f_{\diamond}(t,z,q)$.  The limit $(U_*,V_*)$ is bounded and entire with $U_*(0,0)\ge\varepsilon_1$.  If $\xi_n\to+\infty$ put $r_{\diamond}=r_+$, and if $\xi_n\to-\infty$ put $r_{\diamond}=r_-$.  As in part~\textup{(i)},
\[
 \widehat f_{\diamond}(t,z,0)=r_{\diamond}(t+\theta),
 \;
 \widehat f_{\diamond}(t,z,q)\le r_{\diamond}(t+\theta)-\beta(t+\theta)q.
\]

We let $M_*(t)=\sup_yU_*(t,y)$.  At approximate maxima, exactly as in \cref{prop:global},
\(D^+M_*(t) \le r_\pm(t+\theta)M_*(t) -\delta_0M_*(t)^2 \le r_\pm(t+\theta)M_*(t).\)
Since $\langle r_\diamond\rangle<0$, for every integer $m\ge1$,
\[
 M_*(t)\le M_*(t-mT)
 \exp\!\left(\int_{t-mT}^tr_\diamond(s+\theta)\,ds\right)
 \le C e^{mT\langle r_\diamond\rangle}.
\]
The entire solution is bounded, so $m\to\infty$ provides $M_*(t)=0$, contradicting $U_*(0,0)\ge\varepsilon_1$.  Hence $\lVert u(t,\cdot)\rVert_\infty\to0$.
\end{proof}

\begin{proof}[\upshape\bfseries Proof of \cref{prop:FloquetSpeed}]
Fix $L>0$ and let $\psi_L(t,\xi)>0$ be the principal $T$-periodic Dirichlet eigenfunction for \eqref{eq:FloquetL} at $c=0$, with principal exponent $\lambda_L(0)$. Define
\(\phi_{L,c}(t,\xi)=e^{-c\xi/(2a_0)}\psi_L(t,\xi).\)
The multiplier is positive and independent of $t$, so $\phi_{L,c}$ satisfies the same periodicity and Dirichlet boundary conditions as $\psi_L$. Direct differentiation shows
\begin{align*}
 a_0(\phi_{L,c})_{\xi\xi}+c(\phi_{L,c})_\xi
 &=e^{-c\xi/(2a_0)}\left(a_0(\psi_L)_{\xi\xi}-\frac{c^2}{4a_0}\psi_L\right).
\end{align*}
Inserting the equation for $\psi_L$ yields
\[
 (\phi_{L,c})_t=a_0(\phi_{L,c})_{\xi\xi}
  +c(\phi_{L,c})_\xi+r(t,\xi)\phi_{L,c}
  -\left(\lambda_L(0)-\frac{c^2}{4a_0}\right)\phi_{L,c}.
\]
The periodic Dirichlet principal exponent is uniquely characterized by the existence of a positive eigenfunction. Hence
$\lambda_L(c)=\lambda_L(0)-c^2/(4a_0)$.
Taking the increasing Dirichlet exhaustion limit proves the second identity in \eqref{eq:FloquetExactShift}. The sign assertions are immediate, and the nonlinear consequences follow under precisely the additional conditions of \cref{thm:localized}; they do not require a positive whole-line eigenfunction for the spectral identity itself.
\end{proof}

\section{Planar forced periodic waves in \texorpdfstring{$\mathbb R^N$}{RN}}\label{sec:waves}
We return to \eqref{eq:fullmodel} with a fixed $\tau\ge0$.  We fix once and for all a H\"older exponent $\alpha_w\in(0,\alpha)$; all weighted Schauder spaces in this section and in the weighted Fredholm appendix use $\alpha_w$.  A planar wave depends on
\(\xi=x\cdot e-ct.\)
Since $|e|=1$,
\(\Delta U(t,\xi)=U_{\xi\xi}(t,\xi), \; \nabla\!\cdot(U\nabla V)=(UV_\xi)_\xi.\)
Consequently the profile equations are exactly
\begin{equation}\label{eq:movingfull}
\begin{cases}
 U_t=a(t)U_{\xi\xi}+cU_\xi-
 \eta\chi(t)(UV_\xi)_\xi+U f(t,\xi,U),\\
 \tau V_t=d(t)V_{\xi\xi}+\tau cV_\xi-\nu(t)V+\mu(t)U.
\end{cases}
\end{equation}
The ambient dimension therefore does not enter the planar profile equation; it enters the Cauchy problem through the multidimensional perturbations and the Yukawa resolvent.

We first isolate the scalar wave at $\eta=0$.

\begin{proposition}[Scalar forced periodic wave]\label{prop:scalarwave}
Under \cref{A1,A2,M}, set $\eta=0$, and suppose $a(t)\equiv a_0>0$.  Then for each $c<c_*=2\sqrt{a_0\bar r}$ there exists a unique positive $T$-periodic solution $U_c$ of
\begin{equation}\label{eq:scalarwaveeq}
 U_t=a_0U_{\xi\xi}+cU_\xi+U f(t,\xi,U)
\end{equation}
satisfying
\begin{equation}\label{eq:scalarwaveends}
 U_c(t,-\infty)=p_-(t),
 \;
 U_c(t,+\infty)=0
\end{equation}
uniformly in $t$.  Moreover $U_{c,\xi}(t,\xi)<0$ for all $(t,\xi)$.  No such KPP-type forced wave exists for $c\ge c_*$.  If $I\Subset(-\infty,c_*)$, the profiles have two-sided exponential bounds whose admissible weights can be chosen uniformly for $c\in I$.
\end{proposition}

\begin{proof}
We reduce the assertion exactly to the constant-diffusion theorem of Fang--Peng--Zhao \cite{FPZ2021}.  The rescaling fixes the sign of the drift and the speed threshold.  We put
\[
 y=\frac{\xi}{\sqrt{a_0}},
 \;
 C=\frac{c}{\sqrt{a_0}},
 \;
 \widetilde U(t,y)=U(t,\sqrt{a_0}y),
 \;
 \widetilde f(t,y,u)=f(t,\sqrt{a_0}y,u).
\]
Then \eqref{eq:scalarwaveeq} is equivalent to
\begin{equation}\label{eq:scalarrescaled}
 \widetilde U_t=\widetilde U_{yy}+C\widetilde U_y
 +\widetilde U\,\widetilde f(t,y,\widetilde U).
\end{equation}
We check the hypotheses of Fang--Peng--Zhao \cite[Theorem~1.1]{FPZ2021}. Their (G1) requires joint $C^1$ regularity in $(t,y,u)$ on $u\ge0$ and $T$-periodicity in $t$. The joint $C^1$ clause in \cref{M} verifies (G1), and this property is preserved by the constant spatial rescaling. Their (G2) requires monotonicity in $y$ and $u$, with strict $u$-monotonicity in the favorable limit: these follow respectively from $f_\xi\le0$ and $f_u\le-\beta_*<0$. Their (G3) requires positive averaged favorable linear growth and a positive $M_-$ satisfying $f_-(t,M_-)\le0$ for all $t$; the first part is \cref{F}, while the second holds by choosing
\[
 M_- > \frac{\sup_{t\in[0,T]}r_-(t)_+}{\beta_*},
 \qquad f_-(t,M_-)\le r_-(t)-\beta_*M_-<0.
\]
Finally, (G4) is precisely $\langle r_+\rangle<0$, as required by \cref{F}.  Thus all four hypotheses (G1)--(G4) are verified under \cref{A1,A2,F,M}. The joint $C^1$ assumption in \cref{M} is used here to satisfy (G1); H\"older continuity in time alone does not imply (G1). The favorable state of \eqref{eq:scalarrescaled} is the same periodic solution $p_-$ of \eqref{eq:pminus-early}.  With our orientation, \cite[Theorem~1.1]{FPZ2021} establishes that the forced KPP profile connecting $p_-$ at $y=-\infty$ to $0$ at $y=+\infty$ exists uniquely precisely when
\(C<2\sqrt{\bar r}.\)
Returning to $\xi$ provides
\(c<2\sqrt{a_0\bar r}=c_*.\)
Spatial monotonicity and ordering by climate speed can be proved without assuming continuity of the profiles in $c$. This order will provide the uniform tail smallness needed for the exponential barriers. Let $z^c(t,\xi)$ solve the scalar Cauchy problem in the moving frame with the common constant initial value $z^c(0,\xi)=p_-(0)$. The scalar comparison principle implies
\(0<z^c(t,\xi)\le p_-(t).\)
For $h>0$, the translate $z^c(t,\xi+h)$ has the same initial data and is a subsolution of the unshifted equation, since $f(t,\xi+h,u)\le f(t,\xi,u)$ for $u\ge0$. Therefore
\(z^c(t,\xi+h)\le z^c(t,\xi)\), and $z^c_\xi\le0$ at positive times. Now fix $c_1<c_2<c_*$. Direct subtraction yields
\[
 (z^{c_1})_t-a_0(z^{c_1})_{\xi\xi}-c_2(z^{c_1})_\xi
 -z^{c_1}f(t,\xi,z^{c_1})
 =(c_1-c_2)(z^{c_1})_\xi\ge0.
\]
As $z^{c_1}(0,\cdot)=z^{c_2}(0,\cdot)$, another application of scalar comparison shows
\(z^{c_2}(t,\xi)\le z^{c_1}(t,\xi)\) for all $t\ge0$.

The proof of the scalar existence theorem of Fang--Peng--Zhao, specifically the Poincar\'e iteration in their Lemma~3.3 \cite{FPZ2021}, identifies the forced wave $U_c$ as the limit of the iterates from $p_-(0)$ at integer multiples of $T$. The iterates are decreasing and their pointwise limit is the unique nonzero periodic KPP profile for $c<c_*$. Periodic parabolic compactness upgrades this convergence on every fixed cylinder $[0,T]\times[-R,R]$; by scalar uniqueness the limit is $U_c$. Passing to the limit in the preceding comparisons proves
\begin{equation}\label{eq:speedorder}
 0<U_{c_2}(t,\xi)\le U_{c_1}(t,\xi)\le p_-(t),
 \qquad c_1<c_2<c_*,
\end{equation}
and $U_c(t,\xi+h)\le U_c(t,\xi)$ for every $h>0$.

We put $Z=-U_{c,\xi}\ge0$. Differentiating \eqref{eq:scalarwaveeq}, we find
\[
 Z_t=a_0Z_{\xi\xi}+cZ_\xi+
 \big(f(t,\xi,U_c)+U_cf_u(t,\xi,U_c)\big)Z
 -U_cf_\xi(t,\xi,U_c).
\]
Thus $Z$ is a nonnegative supersolution of the homogeneous linearized equation. The distinct limits in \eqref{eq:scalarwaveends} exclude $Z\equiv0$, and the strong parabolic minimum principle implies $Z(t,\xi)>0$ at every finite point. Hence $U_{c,\xi}<0$.

The tail thresholds can be selected uniformly on compact speed intervals. For an arbitrary compact interval $I=[c_-,c_+]\Subset(-\infty,c_*)$, equation \eqref{eq:speedorder} implies
\[
 U_{c_+}(t,\xi)\le U_c(t,\xi)\le U_{c_-}(t,\xi)
 \qquad(t\in[0,T],\ \xi\in\mathbb R,\ c\in I).
\]
Consequently, using the two endpoint profiles and their pointwise uniform-in-time limits, we obtain
\begin{equation}\label{eq:uniform-scalar-zerotails}
 \lim_{R\to\infty}\sup_{\substack{c\in I,\ t\in[0,T]\\\xi\ge R}}U_c(t,\xi)=0,
 \qquad
 \lim_{R\to\infty}\sup_{\substack{c\in I,\ t\in[0,T]\\\xi\le-R}}
 \big(p_-(t)-U_c(t,\xi)\big)=0.
\end{equation}
In particular, the cutoff $R$ at which the far-field nonlinear coefficients are close to their limits can be chosen before, and independently of, the individual speed $c$. Thus the tail thresholds are uniform on $I$.  Neither continuity in $c$ nor the later Fredholm argument is used.

We prove uniform two-sided exponential bounds using \eqref{eq:uniform-scalar-zerotails} and the exponential coefficient convergence in \cref{M}.  On the unfavorable side, the required bound follows from an exponential supersolution.  If $I\Subset(-\infty,c_*)$, then $\bar r_+<0$ allows one to choose $\rho_+>0$ and $\delta_+>0$, independent of $c\in I$, such that
\(a_0\rho_+^2-c\rho_++\bar r_+\le-2\delta_+ \;(c\in I).\)
We let
\(P_+(t)=\exp\!\left(\int_0^t[r_+(s)-\bar r_+]\,ds\right).\)
By \cref{M} and the uniform far-field control in \eqref{eq:uniform-scalar-zerotails}, after choosing $R$ sufficiently large, independently of $c$, and then choosing $C$ to dominate the bounded traces at $\xi=R$, the function $CP_+(t)e^{-\rho_+\xi}$ is a strict supersolution of \eqref{eq:scalarwaveeq} on $\xi\ge R$ and dominates $U_c$ at $\xi=R$, uniformly for $c\in I$.  To justify the half-cylinder comparison even when $r_+(t)$ changes sign, take $R$ so that $r(t,\xi)\le r_+(t)+\varepsilon$ for $\xi\ge R$ and $\bar r_++\varepsilon<0$. The positive part of $(U_c-CP_+e^{-\rho_+\xi})/P_+$ would attain a strictly positive maximum at a finite spatial point (both summands vanish at $+\infty$). At such a maximum the periodic parabolic inequality has the strictly negative zero-order coefficient $\bar r_++\varepsilon$, a contradiction. The comparison is therefore uniform in $c$.  Hence
\begin{equation}\label{eq:scalarroughrighttail}
 0<U_c(t,\xi)\le C e^{-\rho_+\xi}
 \;(\xi\ge R,\ c\in I).
\end{equation}
At $-\infty$, write $H_c=p_--U_c\ge0$ and $g_-(t,u)=u f_-(t,u)$.  The linearization at the positive periodic state is
\(\mathfrak q_-(t):=\partial_u g_-(t,p_-(t)) =f_-(t,p_-(t))+p_-(t)f_{-,u}(t,p_-(t)),\)
and
\(\langle \mathfrak q_-\rangle =\left\langle p_-f_{-,u}(t,p_-)\right\rangle<0.\)
Since $H_c\to0$ uniformly in both $t$ and $c\in I$ as $\xi\to-\infty$ by \eqref{eq:uniform-scalar-zerotails}, the $C^2$ regularity in \cref{M} provides, after increasing the common $R$,
\[
 (H_c)_t\le a_0(H_c)_{\xi\xi}+c(H_c)_\xi
 +(\mathfrak q_-(t)+\varepsilon)H_c+C_f e^{\sigma_f\xi}
 \;(\xi\le-R)
\]
with one $\varepsilon>0$ that can be chosen uniformly for $c\in I$.  We choose $0<\rho_-<\sigma_f$ so small that
\(a_0\rho_-^2+c\rho_-+\langle \mathfrak q_-\rangle+\varepsilon<0 \;(c\in I).\)
With
\(P_-(t)=\exp\!\left(\int_0^t[\mathfrak q_-(s)-\langle \mathfrak q_-\rangle]ds\right),\)
a sufficiently large multiple of $P_-(t)e^{\rho_-\xi}$ is a strict supersolution of the preceding inhomogeneous inequality, because $\rho_-<\sigma_f$ makes its tail dominate $e^{\sigma_f\xi}$.  We choose one multiple to dominate $H_c$ at $\xi=-R$ for all $c\in I$, using $0\le H_c\le p_-$, and to dominate the inhomogeneous term on $\xi\le-R$. More precisely, because the barrier is a strict supersolution of the inhomogeneous inequality, the function $(H_c-CP_-e^{\rho_-\xi})/P_-$ satisfies a scalar parabolic inequality whose zero-order coefficient is the negative constant $\langle\mathfrak q_-\rangle+\varepsilon<0$. The separate inequality $a_0\rho_-^2+c\rho_-+\langle\mathfrak q_-\rangle+\varepsilon<0$ is precisely what makes the exponential barrier a strict supersolution. Since $H_c$ and the barrier vanish at $-\infty$, the maximum principle on the half-cylinder applies.  We obtain
\begin{equation}\label{eq:scalarroughlefttail}
 0<p_-(t)-U_c(t,\xi)\le C e^{\rho_-\xi}
 \;(\xi\le-R,\ c\in I).
\end{equation}
Periodic Schauder estimates on unit strips, applied to the equations for $U_c$ and $H_c$, provide the same exponential bounds for the derivatives entering the weighted H\"older norms.  The weights $\rho_\pm$ may be decreased without destroying these estimates.

The nonexistence for $c\ge c_*$ and uniqueness are preserved by the bijective change of variables, while \eqref{eq:scalarroughrighttail}--\eqref{eq:scalarroughlefttail} provide the uniform two-sided weights asserted in the proposition.  This proves the result without invoking any extension of the scalar theorem to time-dependent population diffusion.
\end{proof}

\Needspace{6\baselineskip}
For each scalar profile $U_c$, we denote by $V_{\tau,c}$ the bounded $T$-periodic solution of
\(\tau V_t=d(t)V_{\xi\xi}+\tau cV_\xi-\nu(t)V+\mu(t)U_c.\)
When $\tau>0$, the positive damping $\nu_*>0$ makes the whole-past variation-of-constants integral converge and produces a bounded periodic solution; the periodic parabolic maximum principle implies uniqueness.  The weighted version of this construction is proved in \cref{lem:signalweighted}.  When $\tau=0$,
\begin{equation}\label{eq:planarellipticresolvent}
 V(t,\xi)=\frac{\mu(t)}{2d(t)\kappa(t)}
 \int_{\mathbb R}e^{-\kappa(t)|\xi-y|}U_c(t,y)\,dy.
\end{equation}
We choose a smooth two-sided exponential weight $\omega_\rho(\xi)$ satisfying $\omega_\rho(\xi)=e^{\rho_+\xi}$ for $\xi\ge2$ and $\omega_\rho(\xi)=e^{-\rho_-\xi}$ for $\xi\le-2$, with $\omega_\rho\ge1$.  We choose the interpolation on $[-2,2]$ so that $\|\partial_\xi\log\omega_\rho\|_\infty+\|\partial_{\xi\xi}\log\omega_\rho\|_\infty\le C_\omega\max\{\rho_+,\rho_-\}$, where $C_\omega$ is independent of sufficiently small $\rho_\pm$.  A fixed smooth partition of unity interpolates the affine functions $\rho_+\xi$ and $-\rho_-\xi$ and produces such a weight after exponentiation.  Hence all logarithmic derivatives needed below are bounded, and the first two can be made arbitrarily small as $\rho_\pm\downarrow0$.

For the standard (large) periodic parabolic H\"older class $C_{\per}^{\sigma_t,\sigma_x}$, set
$\|h\|_{\rho;\sigma_t,\sigma_x}=\|\omega_\rho h\|_{C_{\per}^{\sigma_t,\sigma_x}}$.
We define $X_\rho^P$ to consist of all $T$-periodic functions $h$ for which
$\omega_\rho h\in C_{\per}^{1+\alpha_w/2,2+\alpha_w}$ and
\[
 \lim_{R\to\infty}\sup_{\substack{t\in[0,T]\\|\xi|\ge R}}
 \big(|\partial_t(\omega_\rho h)|+\sum_{j=0}^{2}|\partial_\xi^j(\omega_\rho h)|\big)=0.
\]
The norm is the full weighted parabolic H\"older norm.  We define $Y_\rho^P$ by
$\omega_\rho g\in C_{\per}^{\alpha_w/2,\alpha_w}$ and
$\sup_{t\in[0,T]}|\omega_\rho(\xi)g(t,\xi)|\to0$ as $|\xi|\to\infty$,
with its usual weighted H\"older norm.  These are closed subspaces of the standard
(large) H\"older spaces and are therefore Banach spaces.  We do \emph{not}
identify them with the closures of compactly supported smooth functions in the
same H\"older norm: such closures are the strictly smaller little H\"older
spaces and require additional small-scale seminorm conditions.  In particular,
spatial decay of derivatives alone does not characterize a little H\"older space.

For the elliptic signal equation at $\tau=0$, let $X_\rho^E$ consist of periodic
functions with $\omega_\rho h\in C_{\per}^{\alpha_w/2}([0,T];C^{2+\alpha_w}(\mathbb R))$
whose weighted spatial derivatives through order two vanish at both ends,
uniformly in $t$.  Let $Y_\rho^E$ consist of periodic functions with
$\omega_\rho g\in C_{\per}^{\alpha_w/2}([0,T];C^{\alpha_w}(\mathbb R))$
and $\omega_\rho g\to0$ at both ends uniformly in $t$.
We use the corresponding full norms; all four spaces are Banach.  The periodic
H\"older norms are taken on the time circle, so matching values at $0$ and $T$
is part of the definition.  No small-scale H\"older vanishing is assumed.
We also use the weighted Banach lattice
\[
 E_{\rho,0}:=\{\phi\in C(\mathbb R):\omega_\rho\phi\in C_0(\mathbb R)\},
 \qquad \|\phi\|_{E_{\rho,0}}=\|\omega_\rho\phi\|_\infty.
\]
Multiplication by $\omega_\rho$ identifies this lattice isometrically and
order-preservingly with $C_0(\mathbb R)$.  The periodic spaces above and this
trace lattice are different Banach spaces; the periodic regularity will be
recovered by interior estimates on the complete time cylinder, not by an
unjustified endpoint smoothing assertion.  The numbers $\rho_\pm$ are chosen
strictly inside the far-field spatial gaps.

We fix a cutoff $\mathfrak h_{\rm L}\in C^\infty(\mathbb R;[0,1])$ such that $\mathfrak h_{\rm L}(\xi)=1$ for $\xi\le-2$ and $\mathfrak h_{\rm L}(\xi)=0$ for $\xi\ge2$.  We define the end-state lifts by
\begin{equation}\label{eq:affine-lifts}
 \mathscr U_-(t,\xi):=\mathfrak h_{\rm L}(\xi)p_-(t),
 \;
 \mathscr V_{-,\tau}(t,\xi):=\mathfrak h_{\rm L}(\xi)q_{-,\tau}(t).
\end{equation}
The scalar tail estimates in \cref{prop:scalarwave} imply $U_c-\mathscr U_-\in X_\rho^P$ whenever $\rho_\pm$ lie strictly below the corresponding scalar decay rates.  The next lemma establishes the analogous assertion for the signal profile rather than assuming it.

The scalar and signal linearizations at $(U_c,V_{\tau,c})$ are
\begin{align}
 \CalA_c h
 &=h_t-a_0h_{\xi\xi}-ch_\xi
 -\big(f(t,\xi,U_c)+U_cf_u(t,\xi,U_c)\big)h,
 \label{eq:Ac}\\
 \CalB_{\tau,c}k
 &=\tau k_t-d(t)k_{\xi\xi}-\tau ck_\xi+\nu(t)k.
 \label{eq:Btau}
\end{align}

\begin{lemma}[Weighted signal inverse and zero-taxis signal tails]\label{lem:signalweighted}
Let $I\Subset(-\infty,c_*)$ be compact and let $\tau\ge0$ be fixed.  We choose $\rho_\pm>0$ strictly below the uniform scalar decay rates supplied by \cref{prop:scalarwave} and decrease them further if necessary.  Then the following assertions hold uniformly for $c\in I$.

If $\tau>0$, the operator $\CalB_{\tau,c}:X_\rho^P\to Y_\rho^P$ is an isomorphism and
\begin{equation}\label{eq:Binverseuniform}
 \sup_{c\in I}\|\CalB_{\tau,c}^{-1}\|_{\mathcal L(Y_\rho^P,X_\rho^P)}<\infty.
\end{equation}
If $\tau=0$, the elliptic operator $\CalB_{0,c}=-d(t)\partial_{\xi\xi}+\nu(t)$ is an isomorphism from $X_\rho^E$ onto $Y_\rho^E$, with an inverse bound independent of $c$.  In either case,
\begin{equation}\label{eq:signal-affine-tail}
 \sup_{c\in I}\|V_{\tau,c}-\mathscr V_{-,\tau}\|_{X_\rho^{\tau}}<\infty,
\end{equation}
where $X_\rho^{\tau}=X_\rho^P$ for $\tau>0$ and $X_\rho^{\tau}=X_\rho^E$ for $\tau=0$.
\end{lemma}

\begin{proof}
We first suppose $\tau>0$ and put $C_I=\sup_{c\in I}|c|$.  We write $\|h\|_{C^0_\rho}:=\|\omega_\rho h\|_\infty$.  We decrease $\rho_\pm$ so that
\(\varrho_\rho:=\|\partial_\xi\log\omega_\rho\|_\infty\)
satisfies
\begin{equation}\label{eq:signal-weight-margin}
 m_{\tau,I}:=\nu_*-d^*\varrho_\rho^2-\tau C_I\varrho_\rho>0.
\end{equation}
After division by $\tau$, the homogeneous equation $\CalB_{\tau,c}k=0$ has the evolution family
\begin{equation}\label{eq:signal-evolution}
 (\mathscr S_c(t,s)\phi)(\xi)
 =e^{-A(t,s)}\int_{\mathbb R}
 \frac{\exp\!\left(-\frac{(\xi-y+c(t-s))^2}{4D(t,s)}\right)}{\sqrt{4\pi D(t,s)}}\phi(y)\,dy,
\end{equation}
where
\(A(t,s)=\tau^{-1}\int_s^t\nu(r)\,dr\)
and
\(D(t,s)=\tau^{-1}\int_s^t d(r)\,dr.\)
Differentiation of the Gaussian verifies that \eqref{eq:signal-evolution} solves the homogeneous equation.  Since $\log\omega_\rho$ is $\varrho_\rho$-Lipschitz, we have
\(\omega_\rho(\xi)/\omega_\rho(y)\le e^{\varrho_\rho|\xi-y|}.\)
We combine $|\xi-y|\le |\xi-y+c(t-s)|+C_I(t-s)$ with the Gaussian moment estimate $\int_{\mathbb R}(4\pi D)^{-1/2}e^{-z^2/(4D)}e^{\varrho_\rho|z|}\,dz\le2e^{\varrho_\rho^2D}$.  Equation \eqref{eq:signal-weight-margin} then implies
\begin{equation}\label{eq:signal-semigroup-weighted}
 \|\mathscr S_c(t,s)\phi\|_{C^0_\rho}
 \le 2e^{-m_{\tau,I}(t-s)/\tau}\|\phi\|_{C^0_\rho},
 \; t\ge s,
\end{equation}
uniformly for $c\in I$.

For a $T$-periodic $h\in Y_\rho^P$, extended periodically to all $t\in\mathbb R$, we define
\begin{equation}\label{eq:signal-periodic-inverse}
 k(t)=\frac1\tau\int_{-\infty}^t\mathscr S_c(t,s)h(s)\,ds.
\end{equation}
The estimate \eqref{eq:signal-semigroup-weighted} shows that the integral converges in the weighted $C^0$ norm and is $T$-periodic.  It also implies
\(\|k\|_{C^0_\rho}\le 2m_{\tau,I}^{-1}\|h\|_{C^0_\rho}.\)
The construction is valid directly for every $h\in Y_\rho^P$; no density in the full H\"older norm is used. We justify the equation and the regularity for such an arbitrary datum. On every bounded time interval, the whole-past integral splits into an integral over an earlier time slice and the ordinary variation-of-constants integral starting at that slice. Strong continuity on $E_{\rho,0}$ makes the resulting function a mild solution. Testing the latter identity against smooth compactly supported space--time functions, or using an approximation of $h$ in the weaker weighted $C^0$ norm, shows that
\(\tau k_t-d(t)k_{\xi\xi}-\tau ck_\xi+\nu(t)k=h\)
in the distributional sense. Here only convergence in $C^0$ is used to pass to the weak equation; no claim of H\"older-norm density is made.

Put $z=\omega_\rho k$ and $H=\omega_\rho h$. The conjugated equation is
\begin{equation}\label{eq:signal-conjugated-V6}
 \tau z_t=d(t)z_{\xi\xi}
 +[\tau c-2d(t)\ell'(\xi)]z_\xi
 +[d(t)(\ell'(\xi)^2-\ell''(\xi))-\tau c\ell'(\xi)-\nu(t)]z+H,
 \qquad \ell=\log\omega_\rho.
\end{equation}
The coefficients have uniformly bounded parabolic H\"older norms on translated cylinders, with constants independent of $c\in I$. The standard interior parabolic regularity theorem for bounded distributional solutions with H\"older forcing applies on each cylinder $[t-2,t+2]\times[\xi-2,\xi+2]$. Since the mild solution has been constructed on the \emph{whole} time axis, no initial-time compatibility condition is needed. Taking the supremum of the resulting uniform interior Schauder estimates, including cylinders across the period joins, yields
\begin{equation}\label{eq:signal-Schauder-bound}
 \|k\|_{\rho;1+\alpha_w/2,2+\alpha_w}
 \le C_{\tau,I}(\|k\|_{C^0_\rho}+\|h\|_{Y_\rho^P})
 \le C'_{\tau,I}\|h\|_{Y_\rho^P}.
\end{equation}
The estimate holds directly for $h\in Y_\rho^P$; it does not use density of smooth functions in the large H\"older norm.

It remains to check the vanishing-jet condition in the definition of $X_\rho^P$. The whole-past formula belongs to $C_{\per}(\mathbb R;E_{\rho,0})$, since the forcing belongs to that space and the evolution is strongly continuous there. Compactness of the time circle implies \(\sup_{t\in[0,T]}|z(t,\xi)|\to0\) as $|\xi|\to\infty$. To check decay of the spatial derivatives, put
\(\epsilon_R:=\sup_{t\in[0,T],\,|\xi|\ge R-2}|z(t,\xi)|.\)
Then $\epsilon_R\to0$, and the uniform $C^{2+\alpha_w}$ bounds on unit strips, together with local interpolation, imply
\[
 \sup_{t\in[0,T],\,|\xi|\ge R}
 (|z_\xi|+|z_{\xi\xi}|)
 \le C\big(\epsilon_R^{(1+\alpha_w)/(2+\alpha_w)}
 +\epsilon_R^{\alpha_w/(2+\alpha_w)}\big).
\]
The constant is independent of $R$ and $c\in I$; when $\epsilon_R=0$ the inequality follows by approximation. Since $H(t,\xi)\to0$ uniformly and the coefficients in \eqref{eq:signal-conjugated-V6} are bounded, the equation also implies $z_t(t,\xi)\to0$ uniformly. Hence $k\in X_\rho^P$, with the stated bound, for every $h\in Y_\rho^P$.
If $h=0$ and $k$ is a periodic homogeneous solution, then periodicity and \eqref{eq:signal-semigroup-weighted} imply, for every integer $n\ge1$,
\(\|k(t)\|_{C^0_\rho}\le2e^{-m_{\tau,I}nT/\tau}\|k(t)\|_{C^0_\rho}.\)
We let $n\to\infty$ and obtain $k=0$.  This proves the parabolic isomorphism and \eqref{eq:Binverseuniform}.

Set $W_{\tau,c}=V_{\tau,c}-\mathscr V_{-,\tau}$.  The definition of $q_{-,\tau}$ and \eqref{eq:affine-lifts} imply the exact identity
\begin{equation}\label{eq:signal-lift-equation}
 \CalB_{\tau,c}W_{\tau,c}
 =\mu(t)(U_c-\mathscr U_-)
 +d(t)q_{-,\tau}(t)\mathfrak h_{\rm L}''
 +\tau c q_{-,\tau}(t)\mathfrak h_{\rm L}'.
\end{equation}
The first term on the right belongs to $Y_\rho^P$ by the scalar tail estimates, while the last two terms have compact spatial support.  Their norms are uniform for $c\in I$.  The inverse estimate \eqref{eq:Binverseuniform} therefore proves \eqref{eq:signal-affine-tail} when $\tau>0$.

We finally suppose $\tau=0$.  The inverse of $\CalB_{0,c}$ is
\begin{equation}\label{eq:signal-elliptic-inverse}
 (\CalB_{0,c}^{-1}h)(t,\xi)
 =\frac1{2d(t)\kappa(t)}\int_{\mathbb R}e^{-\kappa(t)|\xi-y|}h(t,y)\,dy.
\end{equation}
We decrease the weights so that $\rho_\pm<\kappa_*$.  Weighted convolution estimates, followed by differentiation of \eqref{eq:signal-elliptic-inverse} and use of the elliptic equation, show that $\CalB_{0,c}^{-1}:Y_\rho^E\to X_\rho^E$ is bounded.  Formula \eqref{eq:signal-lift-equation} remains valid with its last term absent.  The same source estimate proves \eqref{eq:signal-affine-tail} for $\tau=0$ and completes the proof.
\end{proof}

\begin{lemma}[Continuity of the zero-taxis profile family]\label{lem:profilecontinuity}
Let $I\Subset(-\infty,c_*)$.  For every pair of two-sided weights strictly inside the uniform scalar far-field gaps supplied by \cref{prop:scalarwave}, the map
\(c\longmapsto U_c-\mathscr U_-\)
is continuous from $I$ into the corresponding weighted periodic $C^{1+\alpha_w/2,2+\alpha_w}$ space.  For each fixed $\tau\ge0$, we may decrease those weights, depending on $\tau$, so that \cref{lem:signalweighted} applies; with this choice, the map
\(c\longmapsto V_{\tau,c}-\mathscr V_{-,\tau}\)
is continuous in the associated weighted signal space.
\end{lemma}

\begin{proof}
We let $c_n\to c\in I$.  The scalar profiles satisfy $0<U_{c_n}<p_-$ and the coefficients in their equations are uniformly bounded on compact cylinders.  By periodic interior Schauder estimates, after passing to a subsequence,
\[
 U_{c_n}\longrightarrow U_*
 \;\text{in }C^{1+\alpha'/2,2+\alpha'}_{\rm loc}(\mathbb R\times\mathbb R)
\]
for every $0<\alpha'<\alpha$.  We choose here and below one $\alpha'$ with $\alpha_w<\alpha'<\alpha$.  We pass to the equation and obtain
\((U_*)_t=a_0(U_*)_{\xi\xi}+c(U_*)_\xi+U_*f(t,\xi,U_*).\)
The two-sided exponential estimates in \cref{prop:scalarwave} are uniform for $c\in I$.  They imply, before passing to the limit, that
\(\sup_{t\in[0,T]}|U_{c_n}(t,\xi)|\longrightarrow0 \;(\xi\to+\infty)\)
and
\(\sup_{t\in[0,T]}|U_{c_n}(t,\xi)-p_-(t)|\longrightarrow0 \;(\xi\to-\infty)\)
uniformly in $n$.  Hence $U_*$ has the same end states.  Uniqueness in \cref{prop:scalarwave} implies $U_*=U_c$.  Since every subsequence has the same limit, the full sequence converges locally.

We fix weights $\rho_\pm$ strictly smaller than the uniform decay rates on $I$.  We choose intermediate exponents $\widehat\rho_\pm>\rho_\pm$, still inside the uniform gaps, and let $\omega_{\widehat\rho}$ be a two-sided weight constructed as above with these stronger exponents.  The exponential tail estimates in \cref{prop:scalarwave}, followed by periodic Schauder estimates on unit strips, provide
\begin{equation}\label{eq:profile-strong-tail}
 \sup_n
 \big\|\omega_{\widehat\rho}(U_{c_n}-U_c)\big\|_{C_{\per}^{1+\alpha'/2,2+\alpha'}}
 \le C.
\end{equation}
We put $\delta_\rho=\min\{\widehat\rho_+-\rho_+,\widehat\rho_--\rho_-\}>0$ and set $q_\rho=\omega_\rho/\omega_{\widehat\rho}$.  On each tail $\xi\ge R$ and $\xi\le-R$, once $R\ge2$, the weights are exact exponentials; hence, for $j=0,1,2$,
\[
 \|\partial_\xi^j q_\rho\|_{L^\infty(|\xi|\ge R)}
 +[\partial_\xi^2 q_\rho]_{C^{\alpha_w}((R,\infty))}
 +[\partial_\xi^2 q_\rho]_{C^{\alpha_w}(( -\infty,-R))}
 \le C e^{-\delta_\rho R}.
\]
After multiplication by fixed smooth left- and right-tail cutoffs whose first two derivatives are bounded independently of $R$, the same $C^{2+\alpha_w}$ bound holds for the cutoff ratios.  The product estimate in parabolic H\"older spaces and \eqref{eq:profile-strong-tail} therefore make the contribution of both tails to the weaker weighted $C^{1+\alpha_w/2,2+\alpha_w}$ norm arbitrarily small, uniformly in $n$, by first taking $R$ large.  On the middle strip $|\xi|\le R$, convergence in the stronger local exponent $\alpha'>\alpha_w$ implies convergence in the $C^{1+\alpha_w/2,2+\alpha_w}$ norm by compact embedding.  We first take $R$ large and then $n$ large; we deduce that $U_{c_n}-U_c$ converges to zero in the weighted space.

For $\tau>0$, \cref{lem:signalweighted} first shows that every $V_{\tau,c}-\mathscr V_{-,\tau}$, and hence $(V_{\tau,c})_\xi$, belongs to the required weighted space.  We subtract the two signal equations for $V_{\tau,c_n}$ and $V_{\tau,c}$.  The difference $K_n$ satisfies
\begin{align*}
 \CalB_{\tau,c_n}K_n
 =\mu(t)(U_{c_n}-U_c)+\tau(c_n-c)(V_{\tau,c})_\xi.
\end{align*}
The uniform inverse estimate \eqref{eq:Binverseuniform} therefore implies
\[
 \|K_n\|_{X_\rho^P}
 \le C\left(
 \|U_{c_n}-U_c\|_{Y_\rho^P}
 +|c_n-c|\,\|(V_{\tau,c})_\xi\|_{Y_\rho^P}
 \right),
\]
where $C$ is independent of $n$.  The right-hand side tends to zero by the population convergence just proved.  When $\tau=0$, the operator is independent of $c$ and the difference satisfies $\CalB_{0,c}K_n=\mu(t)(U_{c_n}-U_c)$.  The bounded elliptic inverse in \cref{lem:signalweighted} therefore provides the same convergence.  Thus the signal profiles are continuous in the stated affine weighted spaces.
\end{proof}

\begin{lemma}[Invertibility of the zero-taxis Jacobian]\label{lem:invert}
Suppose that \cref{A1,A2,M} hold and that $a(t)\equiv a_0>0$.  Let $I\Subset(-\infty,c_*)$, and let $U_c$ be the scalar profiles from \cref{prop:scalarwave}.  Suppose in addition that the pinning condition \eqref{eq:pinning} holds for every $c\in I$.  There are positive weights $\rho_\pm$, independent of $c\in I$, such that $\CalA_c:X_\rho^P\to Y_\rho^P$ is an isomorphism and
\begin{equation}\label{eq:Ainverseuniform}
 \sup_{c\in I}\|\CalA_c^{-1}\|_{\mathcal L(Y_\rho^P,X_\rho^P)}<\infty.
\end{equation}
For every fixed $\tau>0$, after reducing $\rho_\pm$ if necessary, $\CalB_{\tau,c}:X_\rho^P\to Y_\rho^P$ is also an isomorphism with an inverse bound uniform for $c\in I$.  For $\tau=0$, the corresponding elliptic signal operator is an isomorphism.  These signal assertions use the same reduced weights as \cref{lem:signalweighted}.
\end{lemma}

\begin{proof}
We first locate the far-field spatial gaps of the scalar operator.  Since $a(t)\equiv a_0$, the right far-field linearization is
\(h_t=a_0h_{\xi\xi}+ch_\xi+r_+(t)h.\)
A mode $h(t,\xi)=P(t)e^{-\lambda\xi}$ is $T$-periodic precisely when
\begin{equation}\label{eq:rightroot}
 a_0\lambda^2-c\lambda+\bar r_+=0,
 \; \bar r_+:=\langle r_+\rangle<0.
\end{equation}
Thus \eqref{eq:rightroot} has the unique positive root $\lambda_+(c)=(c+\sqrt{c^2-4a_0\bar r_+})/(2a_0)$.
At the favorable state $p_-$, we put $\mathfrak q_-=f_-(t,p_-)+p_-f_{-,u}(t,p_-)$.  Since $p_-'/p_-=f_-(t,p_-)$, periodicity implies $\langle f_-(t,p_-)\rangle=0$.  Hence $\langle\mathfrak q_-\rangle=\langle p_-f_{-,u}(t,p_-)\rangle<0$.  A mode $P(t)e^{\lambda\xi}$ at $-\infty$ is periodic precisely when $a_0\lambda^2+c\lambda+\langle\mathfrak q_-\rangle=0$, and this equation has exactly one positive root
\(\lambda_-(c)=\frac{-c+\sqrt{c^2-4a_0\langle \mathfrak q_-\rangle}}{2a_0}>0.\)
Because $I$ is compact, we may choose
\begin{equation}\label{eq:weightchoice}
 0<\rho_+<\inf_{c\in I}\lambda_+(c),
 \;
 0<\rho_-<\inf_{c\in I}\lambda_-(c).
\end{equation}
The exponential convergence in \cref{M}, together with \eqref{eq:scalarroughrighttail}--\eqref{eq:scalarroughlefttail}, implies that the zero-order coefficient of $\CalA_c$ converges exponentially, uniformly for $c\in I$, to $r_+(t)$ at $+\infty$ and to $\mathfrak q_-(t)$ at $-\infty$.  The weighted far-field growth rates are therefore
\[
 \sigma_+(c,\rho_+)=a_0\rho_+^2-c\rho_++\bar r_+<0,
 \;
 \sigma_-(c,\rho_-)=a_0\rho_-^2+c\rho_-+\langle \mathfrak q_-\rangle<0,
\]
with a uniform negative margin on $I$.  Proposition~\ref{prop:fredholm} applies to this family.  In particular, on the conjugated weighted $C_0$ lattice the one-period maps $\mathcal P_c$ satisfy
\begin{equation}\label{eq:PcessentialgapN}
 \sup_{c\in I}r_{\rm ess}(\mathcal P_c)<1,
\end{equation}
and $\CalA_c:X_\rho^P\to Y_\rho^P$ is Fredholm of index zero.

It remains to exclude the Floquet multiplier $1$.  We let $\mathcal P_c$ be the derivative of the scalar one-period map at $U_c$, acting on the weighted lattice $E_{\rho,0}$ corresponding to \eqref{eq:weightchoice}.  Under the isometry $h\mapsto\omega_\rho h$, this operator is conjugate to the $C_0$ monodromy used in Proposition~\ref{prop:fredholm}.  The strictly positive parabolic kernel makes $\mathcal P_c$ positivity improving on $E_{\rho,0}$. By \eqref{eq:PcessentialgapN}, its essential spectral radius is smaller than one. If its spectral radius were at least one, it would therefore be a positive quasi-compact operator with a leading eigenfunctional.  By the monotonicity in \cref{prop:scalarwave},
\(W_c(t,\xi):=-U_{c,\xi}(t,\xi)>0.\)
We differentiate \eqref{eq:scalarwaveeq} in $\xi$ and use the definition \eqref{eq:Ac}.  We obtain the exact identity
\begin{equation}\label{eq:Wsupersolution}
 \CalA_c W_c=-U_c f_\xi(t,\xi,U_c)\ge0.
\end{equation}
The right-hand side is nonnegative and is not identically zero by the pinning condition \eqref{eq:pinning}.  The exponential tail estimates for $U_c$ and the corresponding differentiated equation imply
\(\omega_\rho(\xi)W_c(0,\xi)\longrightarrow0 \;(|\xi|\to\infty),\)
so $W_c(0,\cdot)$ belongs to the positive cone of the weighted $C_0$ space.  We let $\mathcal U_c(t,s)$ denote the positive evolution family of the homogeneous linearized equation.  Variation of constants applied to \eqref{eq:Wsupersolution} provides
\begin{align}
 W_c(T,\cdot)-\mathcal P_cW_c(0,\cdot)
 &=\int_0^T\mathcal U_c(T,s)
 \big[-U_cf_\xi(s,\cdot,U_c)\big](s,\cdot)\,ds.\label{eq:Wduhamel}
\end{align}
The forcing is continuous, nonnegative and not identically zero on one time period by \eqref{eq:pinning}. Consequently it is positive on some open space--time cylinder with time strictly less than $T$ (by periodicity, if necessary). The fundamental solution of the scalar linearized parabolic equation is strictly positive for every positive time separation.  Hence the right-hand side of \eqref{eq:Wduhamel} is strictly positive at every spatial point.  Since $W_c$ is $T$-periodic,
\begin{equation}\label{eq:Pcw}
 0<\mathcal P_cW_c(0,\cdot)<W_c(T,\cdot)=W_c(0,\cdot)
 \;\text{pointwise on }\mathbb R.
\end{equation}

We claim that the spectral radius $r(\mathcal P_c)$ is strictly smaller than one.  We suppose, to the contrary, that $r(\mathcal P_c)\ge1$.  Then \eqref{eq:PcessentialgapN} provides $r(\mathcal P_c)>r_{\rm ess}(\mathcal P_c)$.  The Perron--Frobenius theorem for positive quasi-compact operators on Banach lattices \cite{Schaefer1974} therefore implies a nonzero positive eigenfunctional $\ell_c\in E_{\rho,0}^*$ such that
\begin{equation}\label{eq:adjointPF}
 \mathcal P_c^*\ell_c=r(\mathcal P_c)\ell_c.
\end{equation}
In the present setting $\ell_c$ is in fact strictly positive on the positive cone of $E_{\rho,0}$.  To see this without changing ambient spaces, let
\(J_\rho:E_{\rho,0}\to C_0(\mathbb R),\ J_\rho\varphi=\omega_\rho\varphi,\)
be the order-preserving isometric isomorphism introduced above and define
\(\widetilde\ell_c:=\ell_c\circ J_\rho^{-1}\in(C_0)^*.\)
Then $\widetilde\ell_c$ is a nonzero positive functional, so the Riesz representation theorem identifies it with a nonzero positive Radon measure.  If $0\ne\varphi\in E_{\rho,0}$ and $\varphi\ge0$, positivity improvement implies $\mathcal P_c\varphi>0$ pointwise.  Hence $J_\rho(\mathcal P_c\varphi)>0$ pointwise and
\[
 \langle\ell_c,\mathcal P_c\varphi\rangle
 =\langle\widetilde\ell_c,J_\rho(\mathcal P_c\varphi)\rangle>0.
\]
Equation \eqref{eq:adjointPF} and $r(\mathcal P_c)>0$ then imply $\langle\ell_c,\varphi\rangle>0$.

We set $G_c=W_c(0,\cdot)-\mathcal P_cW_c(0,\cdot)$.  By \eqref{eq:Pcw}, $G_c>0$ pointwise and $G_c\ne0$, so strict positivity of $\ell_c$ provides $\langle\ell_c,G_c\rangle>0$.  On the other hand, \eqref{eq:adjointPF} implies
\[
 \langle\ell_c,G_c\rangle
 =\langle\ell_c,W_c\rangle-\langle\mathcal P_c^*\ell_c,W_c\rangle
 =\big(1-r(\mathcal P_c)\big)\langle\ell_c,W_c\rangle\le0,
\]
a contradiction.  Thus $r(\mathcal P_c)<1$.  Therefore $I-\mathcal P_c$ is invertible by the convergent Neumann series.  The periodic trace and whole-cylinder Schauder construction \eqref{eq:periodicTraceEquationFredholm}--\eqref{eq:periodicSchauderFredholm} in Proposition~\ref{prop:fredholm} then yields, on the same underlying space, the unique periodic solution $h\in X_\rho^P$ for every $g\in Y_\rho^P$.  Hence $\CalA_c$ is an isomorphism.  By \cref{lem:profilecontinuity}, $c\mapsto\CalA_c$ is continuous in the operator norm from $X_\rho^P$ to $Y_\rho^P$.  Since every $\CalA_c$ is invertible, the inverse-map identity
\[
 \CalA_{c_2}^{-1}-\CalA_{c_1}^{-1}
 =\CalA_{c_2}^{-1}(\CalA_{c_1}-\CalA_{c_2})\CalA_{c_1}^{-1}
\]
shows local continuity of the inverse wherever it exists.  For each $c_0\in I$, operator-norm continuity ensures a neighborhood on which $\|(\CalA_c-\CalA_{c_0})\CalA_{c_0}^{-1}\|<1/2$. The Neumann perturbation formula then yields $\|\CalA_c^{-1}\|\le2\|\CalA_{c_0}^{-1}\|$ in that neighborhood. A finite subcover of $I$ proves \eqref{eq:Ainverseuniform}.

The signal assertions follow from \cref{lem:signalweighted}, after decreasing the weights in \eqref{eq:weightchoice} if necessary.  For $\tau>0$, the identity
\(\CalB_{\tau,c_2}-\CalB_{\tau,c_1}=-\tau(c_2-c_1)\partial_\xi\)
shows that $c\mapsto\CalB_{\tau,c}$ is continuous from $I$ into $\mathcal L(X_\rho^P,Y_\rho^P)$; for $\tau=0$ the signal operator is independent of $c$.  The inverse-map identity and the uniform inverse bound in \cref{lem:signalweighted} therefore show that the signal inverse family is continuous in $c$.  Together with the scalar inverse continuity proved above, the triangular zero-taxis Jacobian and its inverse depend continuously on $c\in I$.  This completes the proof.
\end{proof}

\begin{lemma}[Positivity of the weak-taxis branch]\label{lem:tailsign}
Let $I\Subset(-\infty,c_*)$ and fix $\tau\ge0$.  There is a weighted neighborhood of the zero-taxis family such that every solution $(U,V)$ of \eqref{eq:movingfull} in this neighborhood, with the end states \eqref{eq:waveends}, satisfies $U>0$ and $V>0$ on $\mathbb R^2$, uniformly for $c\in I$.
\end{lemma}

\begin{proof}
The left tail and the compact middle region are controlled directly by the weighted topology.  Since $U_c(t,\xi)\to p_-(t)$ uniformly in $t$ and $\min_t p_-(t)>0$, there are $R>1$ and a neighborhood radius, both independent of $c\in I$, such that
\begin{equation}\label{eq:positive-left-middle}
 U(t,\xi)>0\;(\xi\le R).
\end{equation}
Here we first choose $R$ so that $U_c\ge\frac12\min_t p_-$ on $\xi\le-R$ uniformly in $c$, then use compactness of $[0,T]\times[-R,R]\times I$ and continuity of $c\mapsto U_c$ to keep a strictly positive lower bound on the middle strip.

It remains to exclude a negative tail for $\xi\ge R$.  Expanding the taxis divergence in the first profile equation provides
\begin{equation}\label{eq:U-linear-with-V}
 U_t=a_0U_{\xi\xi}+b(t,\xi)U_\xi+q(t,\xi)U,
 \;
 b=c-\eta\chi(t)V_\xi,
 \;
 q=f(t,\xi,U)-\eta\chi(t)V_{\xi\xi}.
\end{equation}
We let $m_+:=-\langle r_+\rangle>0$ and define the positive $T$-periodic function
\begin{equation}\label{eq:Pplus-positivity}
 P_+(t)=\exp\!\left(\int_0^t[r_+(s)-\langle r_+\rangle]ds\right).
\end{equation}
Because the extended reaction in \cref{M} converges uniformly to $f_+(t,u)$ for $u$ in a fixed neighborhood of zero, because $f_+(t,u)=r_+(t)+O(u)$ there, and because both $U$ and $V_{\xi\xi}$ tend to zero in the weighted topology at $+\infty$, we may increase $R$ and shrink the neighborhood so that
\begin{equation}\label{eq:qplus-small}
 |q(t,\xi)-r_+(t)|\le\frac{m_+}{2}
 \;(t\in\mathbb R,\ \xi\ge R,\ c\in I).
\end{equation}
This step uses only boundedness of $\eta\chi$ and the uniform weighted bounds; no sign of $U$ is assumed.

We set $Z=U/P_+$.  Since $P_+'/P_+=r_+-\langle r_+\rangle$, equations \eqref{eq:U-linear-with-V} and \eqref{eq:Pplus-positivity} provide
\begin{equation}\label{eq:Z-positive-tail}
 Z_t=a_0Z_{\xi\xi}+b(t,\xi)Z_\xi+Q(t,\xi)Z,
 \;
 Q=q-r_++\langle r_+\rangle\le-\frac{m_+}{2}
\end{equation}
for $\xi\ge R$, by \eqref{eq:qplus-small}.  We have $Z(t,R)>0$ by \eqref{eq:positive-left-middle} and $Z(t,\xi)\to0$ uniformly in $t$ as $\xi\to+\infty$.  If $Z$ were negative somewhere, its infimum on the periodic half-cylinder $[0,T]\times[R,\infty)$ would be a negative number attained at a finite interior point $(t_0,\xi_0)$.  At that point, after identifying $t=0$ and $t=T$, one has $Z_t=0$, $Z_\xi=0$ and $Z_{\xi\xi}\ge0$.  Equation \eqref{eq:Z-positive-tail} would then imply
\(0=a_0Z_{\xi\xi}(t_0,\xi_0)+Q(t_0,\xi_0)Z(t_0,\xi_0)>0,\)
a contradiction.  Hence $Z\ge0$, and the strong parabolic maximum principle together with $Z(t,R)>0$ implies $Z>0$ for $\xi>R$.  Thus $U>0$ everywhere.

For $\tau>0$, the signal equation can be written
\(\tau V_t-d(t)V_{\xi\xi}-\tau cV_\xi+\nu(t)V=\mu(t)U>0.\)
If $V$ were negative somewhere, its negative global infimum would be attained at a finite point because $V(t,-\infty)=q_{-,\tau}(t)>0$ and $V(t,+\infty)=0$ uniformly in $t$; the parabolic minimum principle provides a contradiction.  Hence $V\ge0$.  Since the right-hand side $\mu(t)U$ is strictly positive, the strong parabolic maximum principle excludes an interior zero and provides $V>0$.  For $\tau=0$, positivity follows immediately from the positive Yukawa representation \eqref{eq:planarellipticresolvent}.  This proves the lemma.
\end{proof}

\begin{proof}[\upshape\bfseries Proof of \cref{thm:waves}]
For fixed $\tau\ge0$ define the profile and target spaces
\[
 \mathbb X_\tau=
 \begin{cases}
 X_\rho^P\times X_\rho^P,&\tau>0,\\
 X_\rho^P\times X_\rho^E,&\tau=0,
 \end{cases}
 \;
 \mathbb Y_\tau=
 \begin{cases}
 Y_\rho^P\times Y_\rho^P,&\tau>0,\\
 Y_\rho^P\times Y_\rho^E,&\tau=0.
 \end{cases}
\]
The weights are those supplied by \cref{lem:invert}, reduced if necessary so that all products below remain inside the same weighted spaces.  By \eqref{eq:affine-lifts}, the base profiles lie in the corresponding affine weighted spaces, while the perturbations belong to the weighted spaces with vanishing spatial jets.  Thus, for $z=(h,k)\in\mathbb X_\tau$, put
\(U=U_c+h, \; V=V_{\tau,c}+k,\)
and define
\begin{align}
 \widehat{\mathcal F}_{\tau,1}(\eta,c;h,k)
 &=U_t-a_0U_{\xi\xi}-cU_\xi
 +\eta\chi(t)(UV_\xi)_\xi-Uf(t,\xi,U),
 \label{eq:Fhat1}\\
 \widehat{\mathcal F}_{\tau,2}(\eta,c;h,k)
 &=\tau V_t-d(t)V_{\xi\xi}-\tau cV_\xi
 +\nu(t)V-\mu(t)U.
 \label{eq:Fhat2}
\end{align}
We choose the neighborhood radius $r_0>0$ smaller than $u_{\rm ext}/2$.  Since $0\le U_c\le p^*$, every perturbation with $\|h\|_{X_\rho^P}\le r_0$ takes values in $(-u_{\rm ext},p^*+u_{\rm ext})$, so the extended reaction from \cref{M} makes the nonlinear profile map well defined on an open Banach neighborhood of the zero-taxis family.  For $\tau>0$ the map takes this neighborhood in $\mathbb X_\tau$ into $\mathbb Y_\tau$.  For $\tau=0$, the second equation in \eqref{eq:Fhat2} is interpreted as the elliptic equation and belongs to $Y_\rho^E$.  The final profiles are positive by \cref{lem:tailsign}, hence they solve the original, unextended system.

The nonlinear terms satisfy the following estimates.  The weighted spaces are stable under multiplication by bounded unweighted functions whose derivatives have the regularity prescribed in \cref{M}.  On a fixed ball $\|z\|_{\mathbb X_\tau}\le r_0$, the product rule provides, with a constant independent of $c\in I$,
\begin{equation}\label{eq:taxisLipIFT}
 \|(UV_\xi)_\xi-(\widetilde U\widetilde V_\xi)_\xi\|_{Y_\rho^P}
 \le C\big(1+\|z\|_{\mathbb X_\tau}+\|\widetilde z\|_{\mathbb X_\tau}\big)
 \|z-\widetilde z\|_{\mathbb X_\tau}.
\end{equation}
Here, when $\tau=0$, the $V$-component is estimated in $X_\rho^E$.  For the reaction set
\(\mathcal N(U)=Uf(t,\xi,U).\)
Its first derivative is
\(D\mathcal N(U)h= [f(t,\xi,U)+Uf_u(t,\xi,U)]h.\)
By the strengthened mixed H\"older bounds in \cref{M}, the coefficient
\(D^2\mathcal N(U)=2f_u(t,\xi,U)+Uf_{uu}(t,\xi,U)\)
 is uniformly bounded in the unweighted parabolic $C^{\alpha_w/2,\alpha_w}$ norm on the relevant bounded ranges and is locally Lipschitz there as a function of $U$ in the corresponding unweighted H\"older norm. The integral identity
\[
 R_c(h)=\int_0^1(1-s)D^2\mathcal N(U_c+sh)[h,h]\,ds
\]
 and the weighted product inequality \(\|h_1h_2\|_{Y_\rho^P}\le C\|h_1\|_{X_\rho^P}\|h_2\|_{X_\rho^P}\) (using $\omega_\rho\ge1$) therefore imply
\begin{align}
 &\|\mathcal N(U_c+h)-\mathcal N(U_c)
   -D\mathcal N(U_c)h\|_{Y_\rho^P}
 \le C\|h\|_{X_\rho^P}^2,
 \label{eq:NemyQuad}\\
 &\|R_c(h)-R_c(\widetilde h)\|_{Y_\rho^P}
 \le C(\|h\|_{X_\rho^P}+\|\widetilde h\|_{X_\rho^P})
       \|h-\widetilde h\|_{X_\rho^P},
 \label{eq:NemyLip}
\end{align}
where $R_c(h)$ denotes the remainder in \eqref{eq:NemyQuad}.  The constants are uniform for $c\in I$ because of \cref{lem:profilecontinuity} and the uniform two-sided tail bounds.

At $\eta=0$ and $z=0$ both equations vanish.  The derivative with respect to $z=(h,k)$ is
\begin{equation}\label{eq:triangularIFT}
 \mathbb L_{\tau,c}:=
 D_z\widehat{\mathcal F}_\tau(0,c;0)
 =
 \begin{pmatrix}
 \CalA_c&0\\
 -\mu(t)&\CalB_{\tau,c}
 \end{pmatrix}.
\end{equation}
By \cref{lem:invert}, $\mathbb L_{\tau,c}:\mathbb X_\tau\to\mathbb Y_\tau$ is an isomorphism, and its inverse norms are uniformly bounded on $I$; denote their supremum by $M_I$.
This follows directly from the triangular inverse: first solve $\CalA_ch=g_1$, then solve
$\CalB_{\tau,c}k=g_2+\mu h$.

We use a quantitative contraction argument rather than requiring differentiability of the base family with respect to $c$.  We write
\(\mathcal R_{\eta,c}:= \widehat{\mathcal F}_\tau(\eta,c;0).\)
Since $(U_c,V_{\tau,c})$ solves the zero-taxis population equation and the signal equation,
\begin{equation}\label{eq:IFTresidual}
 \mathcal R_{\eta,c}
 =\big(\eta\chi(t)(U_c(V_{\tau,c})_\xi)_\xi,0\big),
 \;
 \sup_{c\in I}\|\mathcal R_{\eta,c}\|_{\mathbb Y_\tau}
 \le C_0\eta.
\end{equation}
We define the nonlinear remainder
\[
 \mathcal Q_{\eta,c}(z)=
 \widehat{\mathcal F}_\tau(\eta,c;z)
 -\widehat{\mathcal F}_\tau(\eta,c;0)
 -\mathbb L_{\tau,c}z.
\]
Equations \eqref{eq:taxisLipIFT}--\eqref{eq:NemyLip} provide, for
$\|z\|,\|\widetilde z\|\le r_0$,
\begin{align}
 \|\mathcal Q_{\eta,c}(z)-\mathcal Q_{\eta,c}(\widetilde z)\|_{\mathbb Y_\tau}
 &\le C_1\big(\eta+\|z\|_{\mathbb X_\tau}+\|\widetilde z\|_{\mathbb X_\tau}\big)
 \|z-\widetilde z\|_{\mathbb X_\tau},
 \label{eq:QlipIFT}\\
 \mathcal Q_{\eta,c}(0)&=0.
 \label{eq:QzeroIFT}
\end{align}
with $C_1$ uniform for $c\in I$.

For each $(\eta,c)$ consider
\begin{equation}\label{eq:fixedpointIFT}
 \mathfrak T_{\eta,c}(z)
 =-\mathbb L_{\tau,c}^{-1}
  [\mathcal R_{\eta,c}+\mathcal Q_{\eta,c}(z)].
\end{equation}
We set
\(r_\eta=2M_IC_0\eta.\)
We choose $\eta_{I,\tau}>0$ so small that $r_\eta\le r_0$ and
\begin{equation}\label{eq:etasmallIFT}
 M_IC_1(\eta+2r_\eta)\le\frac14
 \;(0\le\eta\le\eta_{I,\tau}).
\end{equation}
For $\|z\|\le r_\eta$, \eqref{eq:QlipIFT}--\eqref{eq:QzeroIFT} and
\eqref{eq:etasmallIFT} imply
\[
 \|\mathfrak T_{\eta,c}(z)\|
 \le M_IC_0\eta+M_IC_1(\eta+\|z\|)\|z\|
 \le\frac{r_\eta}{2}+\frac{r_\eta}{4}<r_\eta,
\]
and the same ball has contraction factor at most $1/4$.  Thus Banach's fixed-point theorem produces a unique correction
$z_{\eta,\tau,c}=(h_{\eta,\tau,c},k_{\eta,\tau,c})$ in the ball of radius $r_\eta$ for every $c\in I$, and
\(\sup_{c\in I}\|z_{\eta,\tau,c}\|_{\mathbb X_\tau} \le2M_IC_0\eta.\)
In particular, \eqref{eq:etatends} follows with the quantitative rate $O(\eta)$.

The fixed-point map depends continuously on $c$: the base profiles are continuous by \cref{lem:profilecontinuity}, the coefficients in \eqref{eq:Fhat1}--\eqref{eq:Fhat2} are continuous in $c$, and the inverse family is continuous by \cref{lem:invert} $\mathbb L_{\tau,c}^{-1}$.  Since the contraction constant is uniformly smaller than one, the fixed point depends continuously on $c$.  The weighted perturbations vanish at both spatial ends, so the end states \eqref{eq:waveends} are unchanged.  We finally shrink $\eta_{I,\tau}$ if necessary.  The $O(\eta)$ bound places every corrected profile in the uniform neighborhood of \cref{lem:tailsign}; consequently $U_{\eta,\tau,c}>0$ and $V_{\eta,\tau,c}>0$.

The self-map estimate used above actually places the fixed point $z_{\eta,\tau,c}$ in the concentric ball of radius $3r_\eta/4$.  If $\eta>0$, every other planar solution with the same end states and weighted distance less than $r_\eta/4$ from $z_{\eta,\tau,c}$ therefore lies in the contraction ball of radius $r_\eta$ centered at the zero-taxis profile.  The contraction property then forces the two solutions to coincide.  When $\eta=0$, local uniqueness follows from the scalar uniqueness in \cref{prop:scalarwave} together with uniqueness of the associated signal profile.  Thus the wave is locally unique in the planar weighted class stated in the theorem.  The construction is carried out separately for each fixed $\tau>0$ and for $\tau=0$.  In particular, no estimate regards $\tau\partial_t$ as a regular perturbation of the elliptic signal operator when $\tau\downarrow0$.
\end{proof}

\begin{proof}[\upshape\bfseries Proof of \cref{thm:linearresponse}]
For fixed $\tau$ write $z_{\eta,c}=(U_{\eta,\tau,c}-U_c,V_{\eta,\tau,c}-V_{\tau,c})\in\mathbb X_\tau$. In the notation of \eqref{eq:triangularIFT}--\eqref{eq:fixedpointIFT}, set
\[
 b_c=\big(\chi(t)(U_c(V_{\tau,c})_\xi)_\xi,0\big),\qquad
 z_{1,c}=-\mathbb L_{\tau,c}^{-1}b_c.
\]
The weighted resolvent inverses in \cref{lem:invert} make these expressions well defined, and the inverse bound is uniform for $c\in I$. The triangular form of $\mathbb L_{\tau,c}$ shows that the two components of $z_{1,c}$ solve \eqref{eq:responseEq}, establishing existence and uniqueness of $(H_{\tau,c},K_{\tau,c})$.

The nonlinear profile identity and \eqref{eq:IFTresidual} now read
\begin{equation}\label{eq:responseExact}
 \mathbb L_{\tau,c}z_{\eta,c}+\eta b_c
 +\mathcal Q_{\eta,c}(z_{\eta,c})=0.
\end{equation}
The contraction proof of \cref{thm:waves} establishes
$\sup_{c\in I}\|z_{\eta,c}\|_{\mathbb X_\tau}\le2M_IC_0\eta$.
Using \eqref{eq:QlipIFT} with $\widetilde z=0$ and \eqref{eq:QzeroIFT}, we infer
\[
 \|\mathcal Q_{\eta,c}(z_{\eta,c})\|_{\mathbb Y_\tau}
 \le C_1\big(\eta+\|z_{\eta,c}\|_{\mathbb X_\tau}\big)
 \|z_{\eta,c}\|_{\mathbb X_\tau}
 \le C_1(1+2M_IC_0)(2M_IC_0)\eta^2.
\]
Subtract $\eta\mathbb L_{\tau,c}z_{1,c}=-\eta b_c$ from \eqref{eq:responseExact} and apply the inverse of $\mathbb L_{\tau,c}$. We obtain
\[
 \|z_{\eta,c}-\eta z_{1,c}\|_{\mathbb X_\tau}
 \le M_IC_1(1+2M_IC_0)(2M_IC_0)\eta^2,
\]
uniformly for $c\in I$. This is \eqref{eq:responseExpansion}. Division by $\eta>0$ proves the one-sided differentiability at zero. The proof uses only the previously established inverse and quadratic remainder estimates; no differentiability in the climate speed or uniformity as $\tau\downarrow0$ is needed.
\end{proof}

\begin{remark}
The planar wave construction is dimension independent because the profile depends only on $x\cdot e-ct$.  The theorem does not assert nonlinear stability with respect to arbitrary transverse perturbations.  Such a result requires a separate uniform analysis of the transverse Fourier modes of the coupled linearization.
\end{remark}

\section{One-dimensional specializations}\label{sec:recover}
The following specializations recover two earlier models.

\subsection{The scalar time-periodic shifting model}
We take $N=1$, $e=1$, $\eta=0$, $a\equiv1$, and omit the decoupled signal equation.  Then \eqref{eq:fullmodel} becomes
\begin{equation}\label{eq:FPZexact}
 u_t=u_{xx}+u f(t,x-ct,u).
\end{equation}
To state the scalar result we now drop the quantitative chemotaxis slope hypothesis \cref{A2}, which is unnecessary when $\eta=0$, and impose the scalar hypotheses used by Fang--Peng--Zhao: $f\in C^1$, it is $T$-periodic in $t$, nonincreasing in $\xi$ and $u$, the limit $f_-(t,u)=f(t,-\infty,u)$ exists and is strictly decreasing in $u$, 
\(\int_0^T f_-(t,0)\,dt>0,\)
there is $M>0$ with $f_-(t,M)\le0$, and the right limit $f_+(t,u)=f(t,+\infty,u)$ exists with
\(\int_0^T f_+(t,0)\,dt<0.\)
Under these assumptions the favorable scalar periodic state is the unique positive periodic solution of $p'=pf_-(t,p)$ and
\begin{equation}\label{eq:FPZspeed}
 c_*=2\sqrt{\frac1T\int_0^Tf_-(t,0)\,dt}.
\end{equation}

\begin{proposition}[Scalar time-periodic shifting model]\label{prop:FPZ}
Under the preceding scalar specialization and the Fang--Peng--Zhao assumptions, let $p_-(t)$ be the positive $T$-periodic solution of $p'=pf_-(t,p)$ and let $c_*$ be defined by \eqref{eq:FPZspeed}.  The proposition below restates the forced-wave, pulse-wave, and Cauchy results of Fang--Peng--Zhao.  For the logarithmic outer estimates we use the weighted-tail formulation appearing in their Lemma~5.4, which is the form needed in the proof and the only form used here.

\begin{enumerate}[label=\textup{(\roman*)},leftmargin=8mm]
\item A forced $T$-periodic KPP wave $U(t,x-ct)$ satisfying
\(U(t,-\infty)=p_-(t), \; U(t,+\infty)=0\)
exists if and only if $c<c_*$, and it is unique whenever it exists.

\item A forced pulse is a positive $T$-periodic profile satisfying $U(t,\pm\infty)=0$.  If $c>-c_*$, no forced pulse exists.  Suppose, in addition, that for some $r_0>0$ and $m\in\{1,2\}$,
\begin{equation}\label{eq:FPZtail}
 \sup_{t\in[0,T]}|f(t,\xi,0)-f_-(t,0)|
 =o(|\xi|^{-r_0-m})
 \;(\xi\to-\infty).
\end{equation}
If $c<-c_*$ and \eqref{eq:FPZtail} holds with $m=1$, or if $c=-c_*$ and \eqref{eq:FPZtail} holds with $m=2$, then there exist infinitely many forced pulse waves.

\item Let $\phi\in C(\mathbb R,\mathbb R_+)\setminus\{0\}$ be bounded and let $u(t,x;\phi)$ denote the solution with initial datum $\phi$.
\begin{enumerate}[label=\textup{(\alph*)},leftmargin=7mm]
\item If $c\le-c_*$ and
\begin{equation}\label{eq:FPZweightedleft}
 \int_{\mathbb R}e^{-c_*x/2}\phi(x)\,dx<\infty,
\end{equation}
then
\(\lim_{t\to\infty}\sup_{x\in\mathbb R}u(t,x;\phi)=0.\)
\item If $c\ge c_*$, then for every $\mu\in(0,c_*)$,
\(\lim_{t\to\infty}\sup_{|x|\le\mu t}|u(t,x;\phi)-p_-(t)|=0.\)
If, in addition, \eqref{eq:FPZweightedleft} holds, then for every $\gamma\in(0,2/c_*)$,
\(\lim_{t\to\infty}\sup_{|x|\ge c_*t-\gamma\ln t}u(t,x;\phi)=0.\)
\item If $-c_*<c<c_*$ and $U$ is the unique forced KPP profile from part \textup{(i)}, then for every $\mu\in(0,c_*)$,
\(\lim_{t\to\infty}\sup_{x\ge-\mu t} |u(t,x;\phi)-U(t,x-ct)|=0.\)
If \eqref{eq:FPZweightedleft} holds, then for every $\gamma\in(0,2/c_*)$,
\begin{equation}\label{eq:FPZleftlog}
 \lim_{t\to\infty}\sup_{x\le-c_*t+\gamma\ln t}u(t,x;\phi)=0.
\end{equation}
If instead $\liminf_{x\to-\infty}\phi(x)>0$ and there exists $\delta>0$ such that
\(f_u(t,\xi,u)<-\delta\)
for $t\in[0,T]$, $\xi\in\mathbb R$ and
$u\in[-\delta,\delta+\max_{t\in[0,T]}p_-(t)]$, then for some $\sigma>0$,
\begin{equation}\label{eq:FPZexpattract}
 \lim_{t\to\infty}e^{\sigma t}
 \sup_{x\in\mathbb R}|u(t,x;\phi)-U(t,x-ct)|=0.
\end{equation}
\end{enumerate}
\end{enumerate}
\end{proposition}

\begin{proof}
After the specialization, \eqref{eq:FPZexact} is exactly equation (1.1) of Fang--Peng--Zhao \cite{FPZ2021}, with the same favorable and unfavorable orientations.  Their Theorem~1.1 establishes part \textup{(i)} and the speed \eqref{eq:FPZspeed}.  Their condition (1.13) is precisely \eqref{eq:FPZtail}, and their Theorem~1.2 provides part \textup{(ii)}.  The Cauchy dynamics in part \textup{(iii)} follow from their Theorem~1.3; the logarithmic outer estimate is used here with the weighted-tail hypothesis stated explicitly in their Lemma~5.4.  The logarithmic corrections and the global exponential attraction in \eqref{eq:FPZleftlog}--\eqref{eq:FPZexpattract} are scalar conclusions and are not claimed for nonzero taxis elsewhere in this paper.  No quantitative lower bound on $-f_u$ is needed for parts \textup{(i)}--\textup{(iii)(b)} beyond the original Fang--Peng--Zhao hypotheses; the extra derivative condition in part \textup{(iii)(c)} is exactly the one used for their exponential convergence statement.
\end{proof}

\subsection{The one-dimensional logistic chemotaxis model}
We take $N=1$, $e=1$, $a=d=1$, $\eta=1$, $\tau=0$, and set
\(f(t,\xi,u)=r(\xi)-bu\)
with $\chi,\mu,\nu,b$ constant.  Then \eqref{eq:fullmodel} is precisely \eqref{eq:SXintro}.  We use the notation of Shen--Xue: $r$ is bounded and H\"older continuous, the limits $r(\pm\infty)$ exist, and
\(r^*:=\sup_{\xi\in\mathbb R}r(\xi)>0.\)
Their separated-habitat case satisfies
\(r(-\infty)<0<r(+\infty)\) and
\(r(-\infty)\le r(\xi)\le r(+\infty)\), so $r^*=r(+\infty)$; their localized-habitat case satisfies $r(\pm\infty)<0$ and
\(\min\{r(-\infty),r(+\infty)\}\le r(\xi)\).

The coefficient appearing in the Shen--Xue damping condition is, by \eqref{eq:Theta1},
\begin{equation}\label{eq:SXdampingexact}
 \gamma_{1,\sqrt{r^*}}
 =\chi\mu\left(1+\frac12
 \frac{(\sqrt{r^*}-\sqrt\nu)_+}
 {\sqrt{r^*}+\sqrt\nu}\right),
\end{equation}
while $\gamma_0=\chi\mu$.  Their standing hypothesis is therefore
\[
 b>\chi\mu,
 \;
 b\ge\chi\mu\left(1+\frac12
 \frac{(\sqrt{r^*}-\sqrt\nu)_+}
 {\sqrt{r^*}+\sqrt\nu}\right).
\]
In the separated-habitat case, Shen--Xue place the favorable far field at $+\infty$.  After the reflection $x\mapsto-x$, $c\mapsto-c$, one has $\bar r=r^*$ and $\lambda_*=\sqrt{r^*}$, so \eqref{eq:SXdampingexact} and $b>\chi\mu$ are exactly \cref{D}.  In the localized-habitat case, the same standing hypothesis makes the test exponent $\lambda=\sqrt{r^*}$ admissible in \eqref{eq:admissiblesup}.

\begin{corollary}[Global well-posedness in one dimension]\label{cor:SXglobal}
Under the one-dimensional constant-coefficient specialization above, \cref{prop:global} reduces to Shen--Xue Proposition~1.1: if $b>\chi\mu$, every bounded nonnegative uniformly continuous initial datum generates a unique global bounded classical solution.
\end{corollary}

\begin{proof}
In this specialization $\gamma_0=\chi\mu$ and $\beta=b$, so $\delta_0=b-\chi\mu>0$.  The assertion is therefore exactly \cref{prop:global}.
\end{proof}

\begin{corollary}[The separated-habitat case]\label{cor:SXfront}
Assume that $r$ satisfies the separated-habitat conditions above.  Under the preceding specialization and reflection, \cref{thm:front} provides the three propagation regimes of Shen--Xue Theorem~1.1.  In the original orientation, extinction holds for $c>2\sqrt{r^*}$; for $-2\sqrt{r^*}\le c<2\sqrt{r^*}$ the population persists between the moving climate edge and the intrinsic KPP edge; and for $c<-2\sqrt{r^*}$ it spreads in both directions with intrinsic speed $2\sqrt{r^*}$.  If $b>2\chi\mu$, the interior limit is $r^*/b$.
\end{corollary}

\begin{proof}
Here $\bar a=1$, $\bar r=r^*$, hence $c_*=2\sqrt{r^*}$.  The damping assumptions coincide by \eqref{eq:SXdampingexact}.  In one dimension the Euclidean condition $|x|\le(c_*-\eps)t$ and the half-space condition $x\le(C-\eps)t$ reduce to an interval.  With the reflected speed $C=-c$, the three cases $C<-c_*$, $-c_*<C<c_*$, and $C\ge c_*$ translate exactly into the three cases in the statement.  The endpoint $C=c_*$ corresponds to $c=-c_*$ and is covered by part (iii) of \cref{thm:front}; the excluded endpoint $C=-c_*$ corresponds to $c=c_*$, which is not part of the Shen--Xue theorem.  Finally $p_-=r^*/b$ and $\delta_2>0$ becomes $b>2\chi\mu$.
\end{proof}

\begin{corollary}[One-sided initial data]\label{cor:SXfilled}
Assume that $r$ satisfies the separated-habitat conditions above.  Under the same specialization, suppose the original Shen--Xue initial datum vanishes on a left half-line and is bounded away from zero at $+\infty$.  After reflection, it satisfies \eqref{eq:filleddata}.  Theorem~\ref{thm:frontfilled} then provides exactly the two regimes of Shen--Xue Theorem~1.2, including convergence to $r^*/b$ when $b>2\chi\mu$.
\end{corollary}

\begin{proof}
After reflection the climate speed is $C=-c$ and the occupied side is $x\to-\infty$.  The boundary speed in \cref{thm:frontfilled} is
\(s_C=\min\{C,c_*\}.\)
If $C\le c_*$, equivalently $c\ge-c_*$, the occupied population tracks the habitat with speed $C$.  If $C>c_*$, equivalently $c<-c_*$, the boundary moves with the intrinsic speed $c_*$.  Reflecting the two half-line inequalities back to the original coordinate provides the decay and persistence regions of Shen--Xue Theorem~1.2.
\end{proof}

We now assume that $r$ satisfies the localized-habitat conditions above.  The autonomous specialization of \eqref{eq:FloquetL} is
\(\phi''+c\phi'+r(\xi)\phi=\lambda\phi.\)
Shen--Xue establish for their positive generalized eigenfunction the estimate
\begin{equation}\label{eq:SXslope}
 \frac{|\phi_\infty'(\xi)|}{\phi_\infty(\xi)}
 \le K_{SX}:=\frac{\sqrt{8r^*+c^2}+|c|}{2}.
\end{equation}
Their signal condition
\begin{equation}\label{eq:SXnu}
 \nu\ge\nu^*:=
 \frac{(\sqrt{8r^*+c^2}+|c|)^2}{4}
\end{equation}
implies $K_{SX}\le\sqrt\nu$.  In dimension one, \eqref{eq:Psi1zero} shows that the weighted Yukawa correction in \eqref{eq:GammaPhi} vanishes up to this endpoint.

\begin{corollary}[The localized-habitat case]\label{cor:SXlocal}
Under the localized-habitat specialization above, \cref{thm:localized} provides all conclusions of Shen--Xue Theorem~1.3: extinction when $|c|>2\sqrt{r^*}$; for every $\sigma>0$,
\(\lim_{t\to\infty}\sup_{|x\cdot e-ct|\ge\sigma t}u(t,x)=0;\)
persistence on every fixed neighborhood of the moving favorable habitat when $\lambda_\infty(c)>0$; and extinction when $\lambda_\infty(c)<0$ under \eqref{eq:SXnu}.  The case $\lambda_\infty(c)=0$ remains open.
\end{corollary}

\begin{proof}
The choice $\lambda=\sqrt{r^*}$ belongs to $\mathcal A_D$ by \eqref{eq:SXdampingexact}.  Since $\bar r^{\mathrm{sup}}=r^*$ in this autonomous specialization, the arithmetic--geometric mean inequality and \eqref{eq:admissiblesup} imply $c_D^{\rm env}=2\sqrt{r^*}$.  The principal exponent in \eqref{eq:FloquetL} is exactly the generalized principal eigenvalue used by Shen--Xue.  The longitudinal confinement assertion is part (i) of \cref{thm:localized}.  If the exponent is positive, part (iii) provides persistence on every fixed interval in the moving coordinate.  If it is negative, \eqref{eq:SXslope} and \eqref{eq:SXnu} provide $K_\Phi\le\kappa=\sqrt\nu$.  For $N=1$, $\Psi_1=0$ on $[0,1]$, so \eqref{eq:localizeddamping} reduces to $b\ge\chi\mu$, already implied by the standing condition.  The endpoint $K_\Phi=\kappa$ is permitted directly by the one-dimensional endpoint statement in \cref{lem:loglipN}; no limiting argument is required.  Hence part (iv) reproduces their extinction mechanism.
\end{proof}

\section{Discussion}\label{sec:discussion}
The multidimensional formulation reveals a geometric structure that is hidden in one space dimension.  The intrinsic invasion set of the favorable limiting equation is the Euclidean ball of velocities $|s|<c_*$.  A frame remains asymptotically favorable precisely when $s\cdot e<c$.  The persistent velocities are therefore
\(\{s\in\mathbb R^N:|s|<c_*,\ s\cdot e<c\}.\)
This intersection is the correct higher-dimensional replacement of the interval between the two one-dimensional fronts.

The Yukawa functional $\Theta_N$ depends nontrivially on the spatial dimension.  In one dimension the resolvent kernel satisfies a pointwise derivative identity, and the sharp coefficient can be computed in closed form.  In higher dimension the logarithmic derivative of the Yukawa kernel is singular at the pole, so no inequality of the form $|\nabla v|\le\kappa v$ can hold for arbitrary nonnegative sources.  At first contact, we retain the positive part of the kernel combination arising from the planar exponential envelope.  This retains a quadratic bound at the leading edge and is the reason the KPP speed remains linearly determined under \cref{D}.

A further quantitative distinction follows from \cref{prop:YukawaSmallSlope}: the first-contact kernel correction vanishes for all sufficiently small slopes in dimension one, but is positive for every slope in dimensions two and higher. Its leading order is logarithmically reduced in dimension two and quadratic in dimensions at least three. These explicit asymptotics refine the hypothesis \cref{D}; no necessity of that nonlinear sufficient condition is asserted.

The moving-slab problem separates linear growth from nonlinear chemotactic control.  The sign of $\lambda_\infty(c)$ is a linear tracking criterion.  Positive exponent survives the transverse Dirichlet cost on sufficiently wide cylinders and produces persistence near the moving slab.  For a negative exponent, the planar Floquet mode is a natural supersolution only after the weighted Yukawa interaction is controlled.  The functional $\Psi_N$ measures the difference from the one-dimensional cancellation.

For constant population diffusivity, \cref{prop:FloquetSpeed} computes the moving-slab spectral threshold exactly by a positive gauge transformation on bounded Dirichlet intervals. This identity is the classical drift-removal identity and should not be confused with an unconditional nonlinear extinction theorem. Under the additional wave-only $C^1$-in-time and mixed H\"older regularity in \cref{M}, \cref{thm:linearresponse} identifies the first-order effect of taxis through an invertible forced periodic linear system and establishes a quadratic approximation error. These formulas make the time-periodic regime quantitatively accessible, while the sign of the shape response remains problem-dependent.

The fully parabolic problem remains open beyond weak taxis at the level of sharp Cauchy propagation.  A complete extension of \cref{thm:front} to arbitrary fixed $\tau>0$ would require a space--time analogue of \cref{lem:expGreen} that is uniform at the minimizing KPP exponent.  The wave theorem avoids this unresolved point by restricting to constant population diffusion, compact speed ranges below $c_*$, fixed $\tau$, and nondegenerate pinning of the scalar profile; in that regime the far-field spectral gaps and the pinning inequality make the coupled profile Jacobian invertible, after which a quantitative inverse argument applies.

\appendix
\section{Periodic principal exponents and weighted linear operators}\label{sec:appendix}
We shall use the following two facts about periodic principal exponents and weighted linear operators.

\begin{proposition}[Continuity of periodic principal exponents]\label{prop:eigcont}
Let $\Omega\subset\mathbb R^N$ be a bounded $C^{2,\vartheta}$ domain for some $\vartheta\in(0,1)$ and consider
\(\mathcal Q_{A,q}\phi =\phi_t-a(t)\Delta\phi-A(t,x)\cdot\nabla\phi-q(t,x)\phi\)
with $T$-periodicity in time and zero Dirichlet boundary data.  Suppose $A_n,A,q_n,q$ belong to one bounded H\"older set with a common exponent and that $A_n\to A$, $q_n\to q$ uniformly.  Here the principal exponent $\lambda$ is defined by $\mathcal Q_{A,q}\phi=-\lambda\phi$ with $\phi>0$.  Then the corresponding principal exponents converge.  The convergence is locally uniform on such bounded H\"older families of coefficients.
\end{proposition}

\begin{proof}
We let $\lambda_n$ denote the principal exponent and let $\mathcal P_n$ be the time-$T$ Dirichlet map of the equation without the spectral term.  With the present sign convention,
\(r(\mathcal P_n)=e^{\lambda_nT}.\)
Uniform Gaussian upper bounds for the Dirichlet kernels provide
\(\sup_n\|\mathcal P_n\|_{C_0(\Omega)\to C_0(\Omega)}<\infty.\)
To obtain a lower bound for the spectral radii, we fix a ball $B\Subset\Omega$ and a nonzero function $0\le\psi\in C_c(B)$.  The uniform interior lower bound for the positive Dirichlet kernels over one period implies a constant $m_B>0$, independent of $n$, such that
\(\mathcal P_n\psi\ge m_B\psi.\)
Hence $r(\mathcal P_n)\ge m_B$, while the preceding operator bound provides a uniform upper bound.  Therefore $(\lambda_n)$ is bounded.

We normalize the positive periodic eigenfunctions by $\lVert\phi_n\rVert_\infty=1$.  The common H\"older bounds on the coefficients, the uniform bound for $\lambda_n$ and parabolic boundary estimates imply compactness in a fixed parabolic H\"older class; in particular, after extraction,
\(\phi_n\longrightarrow\phi\)
locally in $C^{1,2}$ and uniformly on the closed periodic cylinder.  The limit is nonnegative, nontrivial and satisfies the zero lateral boundary condition.  We pass to the equation and obtain the limiting eigenvalue problem for $(A,q)$.  The strong maximum principle and the parabolic Hopf lemma imply that $\phi>0$ in the interior and has strictly negative outward normal derivative.  Thus the limiting exponent is the principal exponent of $(A,q)$.  Its simplicity shows that every convergent subsequence has the same limit, and consequently the whole sequence converges.

For local uniformity, argue by contradiction.  A violating sequence in a bounded H\"older family has, by compact embedding, a uniformly convergent subsequence of coefficients.  The sequential argument just proved then forces convergence of the associated principal exponents, contradicting the violation.  This proves the proposition.
\end{proof}

\begin{proposition}[Essential gap and periodic trace for the planar scalar linearization]\label{prop:fredholm}
Let $I\Subset(-\infty,c_*)$ and suppose that the zero-order coefficients $q_c(t,\xi)$ of
\(\mathcal A_ch=h_t-a_0h_{\xi\xi}-ch_\xi-q_c(t,\xi)h\)
are $T$-periodic in $t$, uniformly H\"older, and converge exponentially, uniformly for $c\in I$, to $\mathfrak q_-(t)$ as $\xi\to-\infty$ and to $\mathfrak q_+(t)$ as $\xi\to+\infty$.  Choose two-sided weights $\rho_\pm>0$ such that
\begin{equation}\label{eq:appendixfarfieldmargin}
 \sup_{c\in I}\{a_0\rho_+^2-c\rho_++\langle \mathfrak q_+\rangle\}<0,\;
 \sup_{c\in I}\{a_0\rho_-^2+c\rho_-+\langle \mathfrak q_-\rangle\}<0.
\end{equation}
Let $\mathcal P_c$ denote the time-$T$ map on $E_{\rho,0}$, equivalently the time-$T$ map of the conjugated equation on $C_0(\mathbb R)$.  Then
\begin{equation}\label{eq:appendixessentialgap}
 \sup_{c\in I}r_{\rm ess}(\mathcal P_c)<1.
\end{equation}
In particular $I-\mathcal P_c$ is Fredholm of index zero on $E_{\rho,0}$, uniformly for $c\in I$.  If, in addition, $1\notin\sigma(\mathcal P_c)$, then the periodic operator
\(\mathcal A_c:X_\rho^P\longrightarrow Y_\rho^P\)
is an isomorphism.  More precisely, its inverse is obtained from the periodic trace equation and obeys the parabolic Schauder estimate with a constant controlled by $\|(I-\mathcal P_c)^{-1}\|_{\mathcal L(E_{\rho,0})}$ and the fixed coefficient bounds.
\end{proposition}

\begin{proof}
We write $z=\omega_\rho h$.  Outside the fixed transition region of the smooth weight, the conjugated equation is
\(z_t=a_0z_{\xi\xi}+b_{\pm,c}z_\xi+Q_{\pm,c}(t,\xi)z,\)
where
\(b_{+,c}=c-2a_0\rho_+, \; b_{-,c}=c+2a_0\rho_-,\)
and
\[
 Q_{+,c}(t,\xi)=a_0\rho_+^2-c\rho_++q_c(t,\xi),\;
 Q_{-,c}(t,\xi)=a_0\rho_-^2+c\rho_-+q_c(t,\xi).
\]
We let $Q_{\pm,c}^{\infty}(t)$ denote the limits obtained by replacing $q_c$ with $q_\pm$.  By \eqref{eq:appendixfarfieldmargin}, there is $\delta>0$ such that
\(\int_0^TQ_{\pm,c}^{\infty}(t)\,dt\le-4\delta T \;(c\in I).\)
The whole-line evolution of either limiting equation is a Gaussian convolution followed by the scalar factor
\(
 \exp(\int_s^tQ_{\pm,c}^{\infty}(\zeta)d\zeta)
\).
Hence, for every integer $m\ge1$,
\(\|\mathcal P_{\pm,c}^{m}\|_{C_0\to C_0} \le e^{-4\delta mT} \;(c\in I).\)

We estimate the essential norm of the full map.  We fix an integer $m$ so large that
\(e^{-2\delta mT}<\frac18\)
and put $S=mT$.  By exponential convergence of $q_c$, there is $R>1$, independent of $c\in I$, such that
\begin{equation}\label{eq:coeffcloseFredholm}
 |Q_{+,c}(t,\xi)-Q_{+,c}^{\infty}(t)|\le\delta\ (\xi\ge R),\;
 |Q_{-,c}(t,\xi)-Q_{-,c}^{\infty}(t)|\le\delta\ (\xi\le-R).
\end{equation}
All conjugated drifts and zero-order coefficients are bounded by a common constant $M_F$.

We use the following tail estimate.  There are constants $C_1,C_2>0$, depending only on $a_0$, $M_F$, $S$ and the fixed H\"older bounds, such that, if $f\in C_0(\mathbb R)$ is supported in $[R+L,\infty)$ and $\|f\|_\infty\le1$, then
\begin{equation}\label{eq:righttailFredholm}
 \|\mathcal U_c(S,0)f\|_\infty
 \le e^{-2\delta S}+C_1e^{C_1S}e^{-L^2/(C_2S)}.
\end{equation}
To prove this estimate, we compare the full evolution with the Dirichlet evolution on a half-line.  We let $u(t)=\mathcal U_c(t,0)f$ and let $v$ solve the same conjugated equation on the half-line $(R,\infty)$ with zero Dirichlet data at $R$ and initial value $f$.  On $\xi\ge R$, \eqref{eq:coeffcloseFredholm} implies
\(Q_{+,c}(t,\xi)\le Q_{+,c}^{\infty}(t)+\delta.\)
After multiplication by the integrating factor associated with $Q_{+,c}^{\infty}+\delta$, the parabolic maximum principle on $(R,\infty)$ provides
\[
 \|v(S,\cdot)\|_{L^\infty(R,\infty)}
 \le\exp\!\left(\int_0^S[Q_{+,c}^{\infty}(t)+\delta]dt\right)\|f\|_\infty
 \le e^{-3\delta S}.
\]
The Aronson kernel $K_c$ of the full conjugated equation satisfies, uniformly for $c\in I$,
\[
 |K_c(t,\xi;s,y)|\le C_1(t-s)^{-1/2}e^{C_1(t-s)}
 \exp\!\left(-\frac{|\xi-y|^2}{C_2(t-s)}\right).
\]
Since $y\ge R+L$ on the support of $f$, integration of this estimate provides
\begin{equation}\label{eq:FredholmBoundaryLeak}
 \sup_{0<t\le S}\sup_{\xi\le R}|u(t,\xi)|
 \le C_1e^{C_1S}e^{-L^2/(C_2S)}.
\end{equation}
On $(R,\infty)$ the difference $w=u-v$ has zero initial value and boundary trace $w(t,R)=u(t,R)$.  The coefficients are uniformly bounded, so the parabolic maximum principle, after multiplication by $e^{-M_Ft}$, implies
\(\|w(S,\cdot)\|_{L^\infty(R,\infty)} \le e^{M_FS}\sup_{0<t\le S}|u(t,R)|.\)
We combine this bound with \eqref{eq:FredholmBoundaryLeak} and the estimate for $v(S)$; we obtain \eqref{eq:righttailFredholm} after enlarging $C_1,C_2$ and weakening $e^{-3\delta S}$ to $e^{-2\delta S}$.  The same argument on $(-\infty,-R)$ proves the corresponding left-tail estimate.

We choose $L$ so large that the Gaussian term in \eqref{eq:righttailFredholm} is smaller than $1/8$.  We let $\chi_0\in C_c^\infty(\mathbb R)$ satisfy
\(\chi_0=1\ \text{on }[-R-L,R+L], \; \chi_0=0\ \text{outside }[-R-2L,R+2L].\)
The operator
\(K_c:=\mathcal U_c(S,0)M_{\chi_0}\)
is compact on $C_0(\mathbb R)$ for each fixed $c\in I$.  Indeed, the initial data are localized in a fixed compact interval; positive-time parabolic smoothing makes the image of the unit ball equicontinuous, and the Gaussian kernel bound makes its tails uniformly small.  Arzel\`a--Ascoli then provides relative compactness in $C_0$.  For the complementary part, we split $1-\chi_0$ into left and right cutoffs and use \eqref{eq:righttailFredholm} and its left analogue.  After replacing the cutoffs by a partition with values in $[0,1]$, we obtain
\(\mathcal P_c^m=K_c+R_c, \; \sup_{c\in I}\|R_c\|<\frac12.\)
Therefore
\(r_{\rm ess}(\mathcal P_c)^m =r_{\rm ess}(\mathcal P_c^m) \le\|R_c\|<\frac12,\)
which proves \eqref{eq:appendixessentialgap}, uniformly on $I$.

Since $r_{\rm ess}(\mathcal P_c)<1$, the element $I-[\mathcal P_c]$ is invertible in the Calkin algebra and $I-\mathcal P_c$ is Fredholm.  Moreover, for every $\theta\in[0,1]$,
\(r_{\rm ess}(\theta\mathcal P_c)=\theta r_{\rm ess}(\mathcal P_c)<1,\)
so $I-\theta\mathcal P_c$ is Fredholm along the norm-continuous path joining the identity to $I-\mathcal P_c$.  The Fredholm index is constant on this path and hence
\begin{equation}\label{eq:indexzeroMonodromy}
 \operatorname{ind}(I-\mathcal P_c)=0.
\end{equation}

It remains to pass from monodromy invertibility on the \emph{weighted $C_0$ trace lattice}
to the periodic H\"older equation.  The regularity of a periodic solution must be established on the entire time cylinder. Positive-time smoothing of a $C_0$ trace alone does not imply $C^{2+\alpha_w}$ regularity of the inhomogeneous terminal trace.

Let $g\in Y_\rho^P$, extend it periodically to $\mathbb R$ in time, and put
\[
 \mathcal G_cg:=\int_0^T\mathcal U_c(T,s)g(s)\,ds\in E_{\rho,0}.
\]
For the fixed weight, the conjugated evolution preserves $C_0(\mathbb R)$ and is strongly continuous there. Since $g\in Y_\rho^P$, the map $s\mapsto \omega_\rho g(s,\cdot)$ is a continuous $C_0$-valued function on the time circle. Consequently the Bochner integral is well defined in $E_{\rho,0}$; neither smooth-data density in a large H\"older norm nor regularity of the initial trace is required here.
The time-$T$ trace equation is
\begin{equation}\label{eq:periodicTraceEquationFredholm}
 (I-\mathcal P_c)h_0=\mathcal G_cg.
\end{equation}
If $1\notin\sigma(\mathcal P_c)$, it admits precisely one solution
$h_0=(I-\mathcal P_c)^{-1}\mathcal G_cg\in E_{\rho,0}$.
The mild solution on $[0,T]$ with this initial value satisfies
$h(T)=h(0)$ by variation of constants.  Repeating the same solution over
successive periods therefore produces a globally defined periodic mild
solution on $\mathbb R\times\mathbb R$.  Uniqueness of the initial-value
problem implies that the extensions agree across the joins; no assertion
about a priori $C^{2+\alpha_w}$ regularity of $h_0$ has been made.

We conjugate by the smooth weight, $z=\omega_\rho h$.  For the corresponding
parabolic equation the coefficients are uniformly bounded and uniformly
$C^{\alpha_w/2,\alpha_w}$ on all translated unit cylinders, with ellipticity
constant $a_0>0$.  The mild solution and the source are bounded, so the
interior parabolic Schauder estimate, applied on cylinders centered at
\emph{arbitrary} $(t,\xi)\in\mathbb R^2$, yields
\begin{equation}\label{eq:periodicSchauderFredholm}
 \|h\|_{\rho;1+\alpha_w/2,2+\alpha_w}
 \le C_I\big(\sup_{t\in[0,T]}\|h(t)\|_{E_{\rho,0}}+\|g\|_{Y_\rho^P}\big)
 \le C_I\big(1+\|(I-\mathcal P_c)^{-1}\|_{\mathcal L(E_{\rho,0})}\big)
 \|g\|_{Y_\rho^P}.
\end{equation}
The periodic extension across $t=0$ and $t=T$ precedes the interior estimate. The source and coefficients are periodic in the trace topology, and the extension is a bounded distributional solution on every cylinder intersecting an integer multiple of $T$. Thus the same interior estimate holds at the period joins; this is not an initial-boundary Schauder estimate. Here the interior Schauder theorem is used in the standard local regularity form for bounded distributional solutions with H\"older coefficients and source, not in a form requiring a preassigned $C^{2+\alpha_w}$ initial trace. This distinction is essential because the trace was constructed only in the weighted $C_0$ lattice.
The solution is continuous as a periodic $E_{\rho,0}$-valued function;
hence $z(t,\xi)\to0$ uniformly in $t$ as $|\xi|\to\infty$.
The global Schauder norm is bounded.  Set $\epsilon_R=\sup_{t\in[0,T],\,|\xi|\ge R-2}|z(t,\xi)|$. The weighted
$C^{2+\alpha_w}$ bound and the interpolation inequality on unit intervals yield
\[
 \sup_{t\in[0,T],\,|\xi|\ge R}
 (|z_\xi|+|z_{\xi\xi}|)
 \le C\big(\epsilon_R^{(1+\alpha_w)/(2+\alpha_w)}
 +\epsilon_R^{\alpha_w/(2+\alpha_w)}\big)\longrightarrow0.
\]
The constants are uniform for $c\in I$, apart from the inverse-monodromy bound already displayed.  Finally, the conjugated equation and the vanishing of $z$,
$z_\xi$, $z_{\xi\xi}$ and the source imply $z_t\to0$ uniformly at both ends.
Thus $h\in X_\rho^P$ \emph{with the precise space defined in the wave section}.
Conversely, every such periodic classical solution satisfies the trace equation;
therefore the periodic solution is unique.  This proves the isomorphism
$\mathcal A_c:X_\rho^P\to Y_\rho^P$ under $1\notin\sigma(\mathcal P_c)$.
The constants in the estimates are uniform for $c\in I$ except for the
explicit inverse-monodromy norm in \eqref{eq:periodicSchauderFredholm}.
\end{proof}
\section*{Data availability}
No datasets were generated or analyzed in this study.  The results are entirely analytical.

\section*{Conflict of interest}
The author declares no conflict of interest.

\begingroup
\emergencystretch=2em

\endgroup
\end{document}